\documentclass[11pt]{amsart}
\def\SBIMSMark#1#2#3{
 \font\SBF=cmss10 at 10 true pt
 \font\SBI=cmssi10 at 10 true pt
 \setbox0=\hbox{\SBF Stony Brook IMS Preprint \##1}
 \setbox2=\hbox to \wd0{\hfil \SBI #2}
 \setbox4=\hbox to \wd0{\hfil \SBI #3}
 \setbox6=\hbox to \wd0{\hss
             \vbox{\hsize=\wd0 \parskip=0pt \baselineskip=10 true pt
                   \copy0 \break%
                   \copy2 \break% 
                   \copy4 \break}}
 \dimen0=\ht6   \advance\dimen0 by \vsize \advance\dimen0 by 8 true pt
                \advance\dimen0 by -\pagetotal
 \dimen2=\hsize \advance\dimen2 by .25 true in
  \openin2=publishd.tex
  \ifeof2\setbox0=\hbox to 0pt{}
  \else 
     \setbox0=\hbox to 3.1 true in{
                \vbox to \ht6{\hsize=3 true in \parskip=0pt  \noindent  
                \input publishd.tex 
                \vfill}}
  \fi
  \closein2
  \ht0=0pt \dp0=0pt
 \ht6=0pt \dp6=0pt
 \setbox8=\vbox to \dimen0{\vfill \hbox to \dimen2{\copy0 \hss \copy6}}
 \ht8=0pt \dp8=0pt \wd8=0pt
 \copy8
 \message{*** Stony Brook IMS Preprint #1, #2 ***}
}

\usepackage{amssymb}
\usepackage{latexsym}
\usepackage{amscd}
\include{psfig}

\newcommand{\brkOK}{\discretionary{}{}{}}

\newtheorem{theorem}{Theorem}[section]
\newtheorem{introtheorem}{Theorem}
\newtheorem{corollary}[theorem]{Corollary}
\newtheorem{lemma}[theorem]{Lemma}
\newtheorem{definition}[theorem]{Definition}
\newtheorem{proposition}[theorem]{Proposition}
\def\C{{\Bbb C}}
\def\R{{\Bbb R}}
\def\Q{{\Bbb Q}}
\def\N{{\Bbb N}}

\def\D{{\Bbb D}}
\def\S{{{\Bbb R} / {\Bbb Z}}}
\def\T{{\Bbb T}}
\def\QS{{{\Bbb Q} / {\Bbb Z}}}
\def\CDC{{{\Bbb   C} \setminus \overline{\D}}}

\def\int{{\operatorname{Int}}}

\def\param{{{\mathcal P}_d}}
\def\ang{{\mathcal A}_d}
\def\conec{{{\mathcal C}_d}}
\def\shift{{{\mathcal S}_d}}
\def\vshift{{{\mathcal S}^{vis}_d}}

\def\diag{{\mathcal D}}

\def\orb{{\mathcal O}}
\def\ratrel{{\lambda_{\Q}}}
\def\ratcri{{\Lambda_{\Q}}}

\def\cp{{\Theta}}

\renewcommand{\newpage}{\vskip 2.5in \goodbreak\vskip -2.5in \bigskip\bigskip}
\begin{document}

\title{Rational rays and Critical portraits of complex polynomials}
\author{Jan Kiwi}
\maketitle
\thispagestyle{empty}

%\vfill
%\newpage
%$\mbox{}$
%\vfill
\tableofcontents
\vskip -1in
\SBIMSMark{1997/15}{October 1997}{}
\vfil\eject
%\input contents.tex
%\newpage

%\listoffigures
%\vfill
%\newpage
%$\mbox{}$
%\vfill

\vspace{0in}
\part*{Introduction}
\vspace{0cm}
\indent
Since external rays were introduced by
Douady and Hubbard~\cite{orsay-notes} they have played a key
role in the study of the dynamics of complex
polynomials.
The pattern in which external rays approach 
the Julia set allow us to investigate its topology
and to point out similarities and differences
between distinct polynomials.
This pattern can be organized in the form 
of ``combinatorial objects''. 
One of these combinatorial objects is the rational
lamination.

The rational lamination $\ratrel(f)$ of a polynomial
 $f$ with connected Julia set $J(f)$ captures how
 rational external rays land.
More precisely,
 the rational lamination $\ratrel(f)$  is an equivalence relation
 in $\QS$ which identifies two arguments $t$ and $t^\prime$
 if and only if the external rays with arguments $t$ and 
 $t^\prime$ land at a common point (compare~\cite{mcmullen-94}).

The aim of this work is to describe the equivalence
 relations in $\QS$ that arise as the rational lamination
 of  polynomials  with all cycles repelling.
We also describe where in parameter space one can
 find a polynomial with all cycles repelling and
 a given rational lamination.
At the same time we derive some consequences 
 that this study has regarding the topology of Julia sets.

\medskip
To simplify our discussion let us assume, 
 for the moment, that all the 
 polynomials in question are monic and that they have connected 
 Julia sets.

Now let us summarize our results.
A more detailed discussion can be found in the 
 introduction to each Chapter. 
We start with the results regarding the 
 topology of Julia sets.

Under the assumption that all the cycles
 of $f$ are repelling, the Julia set $J(f)$ 
 is a full, compact and connected set which might be
 locally connected or not.
If the Julia set $J(f)$ is locally connected then
 the rational lamination $\ratrel(f)$ completely determines
 the topology of $J(f)$ and the topological dynamics of $f$
 on $J(f)$ (see Proposition~\ref{lam-p}).
When the Julia set $J(f)$ is not locally connected
 it is meaningful to study its topology
 via prime end impressions.
We show that each point in $J(f)$ is contained
 in at least one and at most finitely many prime end
 impressions.
Also, we show that $J(f)$ is  locally connected
 at periodic and pre-periodic points of $f$:

\begin{introtheorem}
Consider a polynomial $f$ 
 with connected Julia set $J(f)$
and all cycles repelling.

(a) If $z$ is a periodic or pre-periodic point then $J(f)$
         is locally connected at z. 
    Moreover, if a prime end impression $Imp \subset J(f)$
     contains $z$ then the prime end impression $Imp$ 
         is the singleton
     $\{z\}$.

(b) For an arbitrary $z \in J(f)$, 
        $z$ is contained in at least one and 
        at most finitely many prime end impressions.
\end{introtheorem}

Roughly, part (a) of the previous Theorem
 is a consequence of the fact that
 the rational lamination $\ratrel(f)$ of a polynomial $f$
 with all cycles repelling is abundant in
 non-trivial equivalence classes. 
As a matter of fact we will show that 
 $\ratrel(f)$ is maximal with respect to some
 simple properties. 
Part (b) of the previous Theorem
 is ultimately a consequence of the following
 result which generalizes one by 
 Thurston for quadratic polynomials (see~\cite{thurston-85}):

\begin{introtheorem}
Consider a point $z$ in the connected Julia set $J(f)$ of
 a polynomial $f$ of degree $d$. 
Provided that $z$ has infinite forward orbit, there are
 at most $2^d$ external rays landing at $z$. 
Moreover, for $n$ sufficiently large, 
 there are at most $d$ external rays landing at $f^{\circ n}(z)$.
\end{introtheorem}

As mentioned before,
we describe which equivalence relations in $\QS$ 
 are the rational lamination $\ratrel(f)$
 of a polynomial $f$ with all
 cycles repelling.
The description will be in terms of critical portraits.
Critical portraits were introduced by Fisher in~\cite{fisher-89}
to capture the location of
the critical points of critically pre-repelling maps
and since then extensively used in the literature 
(see~\cite{bielefield-92,goldberg-93,poirier-93,goldberg-94}).
A critical portrait is a collection $\cp = \{\cp_1 , \dots
\cp_m \}$ of finite subsets of $\S$ that satisfy three properties:

$\bullet$ For every $j$, $|\Theta_j| \geq 2$ and 
$ |d \cdot \Theta_j| = 1 $,

$\bullet$ $\Theta_1 , \dots , \Theta_m $ are pairwise unlinked,

$\bullet$ $\sum (|\Theta_j| - 1) = d - 1$.

\goodbreak

Motivated by the work of Bielefield, Fisher and 
Hubbard~\cite{bielefield-92}, 
 each critical portrait $\cp$ 
 generates an equivalence relation $\ratcri(\cp)$ in $\QS$. 
The idea is that each critical portrait $\cp$ determines a
partition of  the circle 
 into $d$ subsets of lenght $1/d$ which we call $\cp$-unlinked classes. 
Symbolic dynamics of multiplication
 by $d$ in $\S$ gives rise to an equivalence relation 
$\ratcri(\cp)$ in $\QS$
 which is a natural candidate to be the rational
 lamination of a polynomial. 
A detailed discussion of this construction
 is given in Chapter~4 where we make a fundamental 
 distinction between critical portraits.
That is, we distinguish between critical portraits with
 ``periodic kneading'' and critical
 portraits with ``aperiodic kneading''.
We show that critical portraits with aperiodic 
 kneading  correspond to rational laminations of polynomials
 with all cycles repelling:

\begin{introtheorem}
Consider an equivalence relation $\ratrel$ in $\QS$.
$\ratrel$ is the rational lamination 
 $\ratrel(f)$ of  some polynomial $f$
 with connected Julia set and all cycles repelling if and only if
 $\ratrel = \ratcri(\cp)$ for some critical portrait
 $\cp$ with aperiodic kneading.

Moreover, when the above holds, there are at most finitely many
critical portraits $\cp$ such that $\ratrel = \ratcri (\cp)$.
\end{introtheorem}

In the previous Theorem, 
the fact that every rational lamination
is generated by a critical portrait is 
consequence of the study of rational laminations discussed
 in Chapter 3.
The existence of a polynomial with a given 
rational lamination relies on finding 
a polynomial in parameter space with the 
desired rational lamination.
Alternatively, the results of~\cite{bielefield-92} can be used
to give a simpler proof of Theorem 3. We will not do this here. 

In parameter space, following Branner and Hubbard~\cite{branner-88},
 we work in the space $\param$ of monic centered polynomials
 of degree $d$. 
That is, polynomials of the form:
 $$z^d  + a_{d-2} z^{d-2} + \cdots +a_0.$$
The set of polynomials $f$ in $\param$ with connected
 Julia set $J(f)$ is called the connectedness locus $\conec$.
We search for polynomials in $\conec$ by looking at $\conec$
 from outside. 
Of particular convenience for us is to explore the 
 shift locus $\shift$.
The shift locus $\shift$ is the open, connected 
 and unbounded set of polynomials $f$ such that
 all critical points of $f$ escape to infinity.
The shift locus $\shift$ is the unique hyperbolic component
 in parameter space
 formed by polynomials with all cycles repelling.
We will concentrate in the set $\partial \shift \cap \conec$ 
 where the shift locus $\shift$ 
 and the connectedness locus $\conec$ ``meet''. 
Conjecturally, every polynomial with all cycles
 repelling and connected Julia set lies in 
 $\partial \shift \cap \conec$.

To describe where in parameter space we can find
 a polynomial with a given rational lamination,
 in Chapter 2, we cover $\partial \shift \cap \conec$
 by smaller dynamically defined sets that we 
 call the ``impressions of  critical portraits''.
More precisely, inspired by Goldberg~\cite{goldberg-94}, we show 
 that each critical portrait $\cp$ naturally defines
 a direction to go from the shift locus $\shift$ to 
 the connected locus $\conec$. 
Loosely,  the set of polynomials
 in $\partial \shift \cap \conec$ reached by a given 
 direction $\cp$ is
 called the impression $I_{\conec}(\cp) \subset \conec$ 
 of the critical portrait $\cp$.
Each critical portrait impression is a closed connected subset 
 of $\partial \shift \cap \conec$ 
 and the set of all impressions covers 
 all of $\partial \shift \cap \conec$.
It is worth to point out that, for quadratic polynomials,
 there is a one to one correspondence between 
 prime end impression of the Mandelbrot
 set and impressions of quadratic critical portraits.

We characterize the impressions that contain polynomials
 with all cycles repelling and a given rational lamination:

\begin{introtheorem}
Consider a map $f$ in the impression 
 $I_{\conec}(\cp) \subset \conec$ of a critical portrait $\cp$.

If $\cp$ has aperiodic kneading then $ \ratrel(f) = \ratcri(\cp)$
and all the cycles of $f$ are repelling.

If $\cp$ has periodic kneading then  at least one
cycle of $f$ is non-repelling.
\end{introtheorem}

The previous Theorem 
 leads to a new proof of the Bielefield-Fisher-Hubbard 
 realization Theorem for critically pre-repelling 
 maps~\cite{bielefield-92}.
This proof replaces the use of Thurston's characterization
 of post-critically finite maps~\cite{douady-93a}
 with parameter space techniques.
The parameter space techniques shed some light on 
 how polynomials are distributed in $\partial \shift \cap \conec$
 according to their combinatorics.

For quadratic polynomials, the results previously stated
 are well known, sometimes in a different but equivalent
 language 
 (see~\cite{orsay-notes,thurston-85,lavaurs-86,atela-92,
douady-93,mcmullen-94,schleicher-94,milnor-95}). 
Also, the Mandelbrot Local Connectivity
 Conjecture states that quadratic impressions are singletons.
For higher degrees, we do not expect this to be true.
That is, we conjecture the existence non-trivial impressions of critical 
 portraits with aperiodic kneading.

\medskip
This work is organized as follows:

In Chapter 1, we study external rays that land 
at a common point.
Namely, following Goldberg and Milnor, the {\it type}
$A(z)$ of a point $z$ and the {\it orbit portrait} 
$A(\orb)$ of an orbit $\orb$ are introduced. 
Their basic properties are discussed and applied
to prove Theorem 2.

In Chapter 2,  we cover $\partial \shift \cap \conec$
 by critical portrait impressions. 
Here, the main issue is to overcome the nontrivial
 topology that the shift locus has for degrees greater
 than two~\cite{blanchard-91}. 
Inspired by Goldberg we do so by restricting 
 our attention to a dense subset of the shift locus that we call 
 the {\it visible shift locus}. In the visible shift locus
 one can introduce coordinates by means of {\it critical 
 portraits}. 
Then, after endowing  the
 set of critical portraits with the {\it compact-unlinked
 topology}, it is not difficult to
 define {\it critical portrait impressions}.
In this Chapter, we also discuss
 the pattern in which  external rays of polynomials
 in the visible shift locus land (see Section~\ref{dyn-s}).
This discussion, although elementary, plays
 a key role in the next two Chapters. 

In Chapter 3,
 we discuss the basic properties of rational
 laminations and at the same time 
 we proof Theorem 1. 
On one hand part (a)
 of this Theorem relies on working with {\it puzzle
 pieces}.  
On the other, we deduce part (b) of Theorem 1 from Theorem 2.
This argument makes use of an auxiliary
 polynomial in the visible shift locus. 
That is, 
 we profit from the fact that,
 for polynomials in the visible shift locus,
 the pattern in  which external rays land 
 is transparent.
 
In Chapter 4, we show how 
 each critical portrait $\cp$ gives rise to an equivalence
 relation $\ratcri(\cp)$ in $\QS$ and proceed to collect
 threads from Chapters~2 and~3 to prove Theorems 3 and 4.

\bigskip
\noindent
{\bf Acknowledgments:} This work is a copy of the author's Ph. D. Thesis.
I am grateful to Jack Milnor for his enlightening comments and suggestions
as well as for his support and good advise. 
I am grateful to Misha Lyubich 
for his insight and questions and to  
Alfredo Poirier for explaining me his work.
This work was partially supported by a ``Beca Presidente
de la Rep\'ublic'' (Chile).

\newpage
\vspace{0in}
\part*{Chapter 1: Orbit Portraits}
\vspace{-0.5cm}
\indent
\section{Introduction}
The main purpose of this chapter is to study
external rays that land at a common point in 
the Julia set $J(f)$ of a polynomial $f$.
The main result here is to give an upper 
bound on the number of external rays that 
can land at a point with infinite forward
orbit.

For quadratic polynomials, it follows from Thurston's work on quadratic
laminations 
that at most 4 rays can land at a point $z$ with
infinite forward orbit. Moreover, all but finitely many
forward orbit elements are the landing point of at most 2 rays
(see~\cite{thurston-85} Gaps eventually cycle).
Here we generalize this result:

\begin{theorem}
\label{con-th}
Let $f$ be a degree $d$ monic polynomial with connected
Julia set $J(f)$. If $z \in J(f)$ is 
a point with infinite forward orbit
then at most 
$2^d$ external rays can land at $z$. Moreover,  for $n$ sufficiently
large,
at most $d$ external rays can land at $f^{\circ n}(z)$. 
\end{theorem}

\begin{figure}[htb]
\centerline{\psfig{file=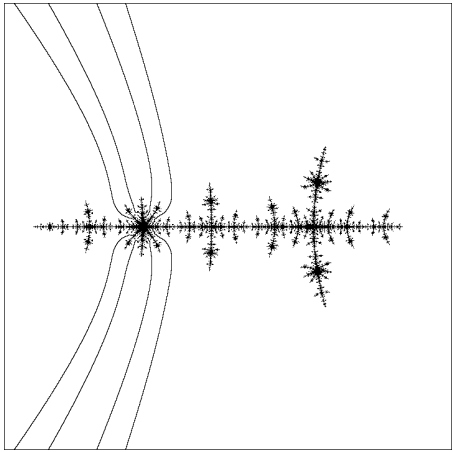,width=170pt}}
\label{8-ray-f}
\caption{The Julia set of the cubic polynomial
$f(z) = z^3 -1.743318 z + 0.50322$ with 
eight rays ``landing'' at the critical point
$c_- = -0.581106$.}
\end{figure}

Although our main interest is the finiteness
given by this result, we should comment
on the bounds obtained. For quadratic polynomials
both bounds are sharp and due to Thurston.
For higher degree polynomials, we 
expect $2^d$ to be optimal in the statement of
the Theorem (see Figure~\ref{8-ray-f}). But we do not know if 
there is an infinite 
orbit of a cubic polynomial
with exactly $3$ external rays landing at each
orbit element.

\medskip
The main ingredients in the proof of Theorem~\ref{con-th} 
are the ideas and techniques
introduced by Goldberg and 
Milnor~\cite{goldberg-93,milnor-95} to study external rays 
that land at a common point $z$ i.e. ``the type of $z$".

\smallskip
This Chapter is organized as follows:

In Section~\ref{preliminar}, we recall some results from polynomial
dynamics. For further reference see~\cite{milnor-90,carlesson-93}. 

In Sections~\ref{portraits} and~\ref{fin-por}, following Goldberg and
Milnor~\cite{goldberg-93}, orbit portraits are introduced and  
some of their basic properties are discussed. 

In Section~\ref{pop}, we apply these
properties to obtain  bounds on the number of cycles
participating in a periodic orbit portrait.
Also, this illustrates
 some of the ideas involved in the proof of Theorem~\ref{con-th}.

In Section~\ref{wop}, we prove Theorem~\ref{con-th}.

\section{Preliminaries}
\label{preliminar}

\noindent \indent
 Here we recall some facts about polynomial dynamics.
For more background material we refer the reader to~\cite{milnor-90}.

Consider a monic polynomial $f: \C \rightarrow \C$
of degree $d$.  Basic tools to understand the
dynamics of $f$ are the Green function $g_f$
and the B\"ottcher map $\phi_f$.

The Green function $g_f$ measures the escape rate of 
points to $\infty$:
$$ \begin{array}{lccl}
 g_f :& \C & \rightarrow & {\R}_{\geq 0} \\
    & z & \mapsto     & \lim \frac{\log_+ | f^{\circ n}(z)|}{d^n}  \\
\end{array}$$
It is a well defined continuous function which vanishes 
on the filled Julia set $K(f)$ and satisfies
the functional relation: 
$$g_f (f(z)) = d g_f(z).$$ 
Moreover, $g_f$ is positive and harmonic
in the basin of infinity $\Omega(f)$. In $\Omega(f)$, 
the derivative of $g_f$ vanishes  at $z$ if and only if 
$z$ is a pre-critical
point of $f$. In order to avoid confusions,
we say that $z$ is a {\bf singularity} of $g_f$. 

The B\"ottcher map $\phi_f$ conjugates 
$f$ with $z \mapsto z^d$ in a neighbourhood
of $\infty$.  The germ of $\phi_f$ at $\infty$
is  unique up to conjugacy 
by $z \mapsto \zeta z$ where $\zeta$ is a $(d-1)st$ root of unity.
Since $f$ is monic we can normalize $\phi_f$ to be asymptotic to the identity:
$$\frac{\phi_f(z)}{z} \rightarrow 1$$
as $z \rightarrow \infty$.
Observe that near infinity, $g_f(z) = \log | \phi_f(z) |$.

For the purpose of simplicity, 
we make a distinction according to whether the Julia
set $J(f)$ is connected or disconnected. 

\subsection{Connected Julia sets}
\noindent \indent
The Julia set $J(f)$ is connected if and only if
all the critical points of $f$ are non-escaping.
That is, the forward orbit of the critical points
remains bounded. Thus, $g_f$ has no singularities 
in $\Omega(f)$. Moreover, 
the B\"ottcher map extends to the basin of infinity $\Omega(f)$,
$$ \phi_f : \Omega(f) \rightarrow \CDC$$
and $\phi_f (f(z)) = (\phi_f(z))^d$ for $z \in \Omega(f)$.
Furthermore,  
$$g_f (z) = \log_+ | \phi_f (z)| \mbox{ for } z \in \Omega(f).$$

An {\bf external ray} $R^t_f$ is the pre-image of radial line
$(1, \infty) e^{2 \pi i t}$ under $\phi_f$, i.e.
$$ R^t_f = \phi^{-1}_f ( (1, \infty) e^{2 \pi i t}).$$
Thus, 
external rays
are curves that run from infinity
towards the Julia set $J(f)$. If $R^t_f$ has a well defined
limit $z \in J(f)$ as it approaches the Julia set $J(f)$
we say that $R^t_f$ {\bf lands} at $z$.

External rays 
are parameterized by the circle $\S$ and
$f$ acts on external rays as multiplication by $d$.
(i.e. $f(R^t_f) = R^{dt}_f$ ).
A ray $R^t_f$ is said to be rational if $t \in \QS$.
Rational rays can be either periodic or pre-periodic
according to whether $t$ is periodic or pre-periodic under 
$m_d : t \mapsto dt$.
A periodic ray always lands at a
repelling or parabolic periodic point.
A pre-periodic ray lands at a pre-repelling or
pre-parabolic point~\cite{orsay-notes}. Conversely,
putting together results of Douady, Hubbard, Sullivan and
Yoccoz~\cite{hubbard-93,milnor-90}, we have the following:

\begin{theorem} 
\label{land-th}
Let $z$ be a parabolic or repelling periodic
point  in a connected Julia set $J(f)$.
Then there exists at least one periodic ray landing at
$z$. Moreover, all the rays that land at $z$ are periodic
of the same period.
\end{theorem}

\subsection{Disconnected Julia sets}
\label{dis-jul}
\noindent \indent
A polynomial $f$ has disconnected Julia set $J(f)$
if and only if  some critical point of $f$ lies
in the basin of infinity $\Omega(f)$. 
In this case, the B\"ottcher map does not extend to all of $\Omega(f)$.
It extends, along flow lines,  to the basin of infinity under
the gradient flow $\mathrm{grad} g_f$.
Following Levin and Sodin~\cite{levin-91}, the {\bf reduced basin 
of infinity} $\Omega^*(f)$ is the basin of
infinity under the gradient flow $\mathrm{grad} g_f$.
Now 
$$\phi_f : \Omega^*(f) \rightarrow U_f \subset \CDC$$
is a conformal isomorphism from
 $\Omega^*(f)$ onto a starlike (around $\infty$) domain $U_f$.
A flow line of $\mathrm{grad} g_f$ in $\Omega^*(f)$
is an {\bf external radius}. 
An external radius maps into
an external radius by $f$. Thus, $f (\Omega^*(f)) \subset \Omega^*(f)$.

External radii are 
parameterized by $\S$. 
More precisely,  for $t \in \S$ let
$(r, \infty) e^{2 \pi i t}$ be the maximal
portion of $(1, \infty) e^{2 \pi i t}$ contained
in $U_f$. The external radius $R^{*t}_f$ with argument $t$ 
is 
$$ R^{*t}_f = \phi^{-1}_f ( (r, \infty)e^{2 \pi i t}).$$
As one follows the external radius $R^{*t}_f$ from $\infty$ towards
the Julia set $J(f)$ one might hit a singularity
$z$ of $g_f$ or not. In the first case, $r > 1$ and
we say that $R^{*t}_f$ {\bf terminates} at $z$. In the second case,
$r =1$ and $R^{*t}_f$ is in fact the smooth external ray $R^{t}_f$
with argument $t$.
Notice that, from the point of view of the gradient flow, an
 external radius  which terminates at a singularity $z$
 is an unstable manifold of $z$. 

Under iterations of $f$, each point $z$ in 
 the basin of infinity $\Omega(f)$ eventually maps
 to a point in the reduced basin of infinity $\Omega^*(f)$.
Say that $f^{\circ n}(z) \in \Omega^*(f)$
 and that the local degree of $f^{\circ n}$ at $z$ is $k$.
After a conformal change of coordinates, the gradient flow lines 
 nearby $z$ are the pre-image under $w \mapsto w^k$ of the 
 horizontal flow lines near the origin. 
Thus, at a singularity $z$
 of $g_f$, there are exactly $k$ local unstable and
 $k$ local stable manifolds which alternate as one goes
 around $z$. 
A local unstable manifold is contained
 in $\Omega^*(f)$ if and only if it is part of an external radius
 that terminates at $z$.

Now let $\theta_1, \dots, \theta_l$ be the arguments
 of the external radii that terminate at critical points
 of $f$. 
Since every pre-critical point of $f$ is a singularity
 of $g_f$, the external radii with arguments in 
 $$ \Sigma = \bigcup_{n \geq 0} m^{-n}_d (\{ \theta_1, \dots, \theta_l \})$$
 also terminate at a singularities.
Since every singularity is a pre-critical point
 we have smooth external rays defined for arguments
 in $\S \setminus \Sigma$.
Following Goldberg and Milnor, for $t \in \S$ let
 $$R^{t^\pm}_f = \lim_{s \rightarrow t^{\pm}} R^s_f.$$   
If $t \notin \Sigma$ then $R^{t^\pm}_f$ coincide, and
 we say that $R^t_f$ is a smooth external ray.
If $t \in \Sigma$ then $R^{t^\pm}_f$ do not agree,
and we say that they are {\bf non-smooth or bouncing rays}
with argument $t$.

Notice that $f(R^{t^\epsilon}_f) = R^{dt^\epsilon}_f$.
We say that $R^{t^\epsilon}_f$ is periodic or pre-periodic
if $t$ is periodic or pre-periodic  under $m_d : t \mapsto dt$.
Here, we also have that periodic rays land at 
repelling or parabolic periodic points.
But, there might be rays, which are not
periodic, landing at a periodic point $z$.
Following Levin and Przyticky~\cite{levin-95}, 
the landing Theorem stated for connected Julia
sets generalizes to:

\begin{theorem}
\label{dis-lan-th}
Let $z$ be a repelling or parabolic periodic point.
 Then there exists at least one external
ray landing at $z$. Moreover, 

Either all the external rays, smooth and non-smooth,
landing at $z$ are periodic of the same period,

Or, the arguments of the external rays, smooth and non-smooth,
landing at $z$ are irrational and form a Cantor set.
Furthermore, $\{ z \}$ is a 
connected component of $J(f)$ and there are non-smooth
rays landing at $z$.
\end{theorem}

\section{Orbit Portraits}
\label{portraits}

\noindent \indent
We fix, for this Chapter, a monic polynomial
$f$ of degree $d$ with Julia set $J(f)$
(possibly disconnected). Our  goal is to study 
external rays that land at a common point:

\begin{definition}
Consider a point $z \in J(f)$. Suppose that at least 
one external ray lands at $z$ and that all the external
rays which land at $z$ are smooth. We say that 
$$A(z) =\{ t \in \S : R^t_f \mbox{ lands at } z \}$$
is the {\bf type} of $z$.

\newpage
Let $\orb = \{ z , f(z) , \dots \}$
be the forward orbit of $z$, we say that
$$A(\orb) = \{ A(w) : w \in \orb \}$$
is the 
{\bf orbit portrait} of $\orb$.
In particular, when $\orb$ is a periodic cycle 
we say that $A(\orb)$ is a {\bf periodic orbit portrait}.
\end{definition}

\begin{figure}[hbt]
\centerline{\psfig{figure=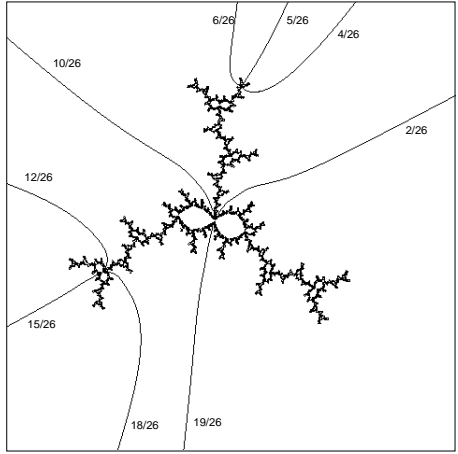,width=170pt}}
\medskip
\caption{External rays landing at a period 3 orbit of a cubic polynomial.}
\label{3-cycle-fig}
\end{figure}

Figure \ref{3-cycle-fig} shows the external rays landing at a period 
3 orbit $\orb$ of a cubic polynomial with orbit portrait
$A(\orb)$:
$$\{ \{2/26,10/26,19/26\}, \{4/26,5/26,6/26\},\{12/26,15/26,18/26\} \}.$$

\noindent
{\bf Remark:}
Below we will see that if the type $A(z)$ of $z$ is well
defined then the type $A(f(z))$ is also well defined
(Lemma~\ref{inv-por}).

\smallskip
Theorem~\ref{con-th} follows from 
the slightly more general:

\begin{theorem}
\label{wop-th}
Consider a monic polynomial $f$ of degree $d$ with
Julia set $J(f)$. If $A(z)$ is the type of 
a Julia set element $z$ with infinite forward orbit
then the cardinality of $A(z)$ is at most $2^d$.
 Moreover, for $n$ sufficiently large, the cardinality of
$A(f^{\circ n} (z))$ is at most $d$.
\end{theorem}

\noindent
{\bf Remark:} Above, in the statement of
Theorem, we do not assume that the Julia set is connected
because, in Chapter 3, we will need to apply this result
for polynomials with disconnected Julia set.

\medskip
Now we list the basic properties of types, 
proofs are provided at the end of this section.
Recall that $m_d: t \rightarrow dt$ denotes multiplication
by $d$ modulo $1$.

\smallskip 
Types are invariant under dynamics:

\begin{lemma}
\label{inv-por}
If $A(z)$ is the type of
$z$  then $A(f(z)) = m_d (A(z))$.
Moreover, ${m_d}_{|A(z)}$ is a $k$ to $1$ map where
$k$ is the local degree of $f$ at $z$.
\end{lemma}

\smallskip
Provided that $z$ is not a critical point,  
the transition from the type
of $z$ to that of its image $f(z)$ 
is cyclic order preserving:

\begin{lemma}
\label{mon-por}
If $z$ is not a critical point of $f$ then
$${{m_d}_|}_{A(z)}: A(z) \rightarrow A(f(z))$$
is a cyclic order preserving bijection.
\end{lemma}

\smallskip
Often we study types of several points
at the same time.
Since smooth external rays are disjoint,
types of distinct points embed in $\S$ in 
an ``unlinked'' fashion.

\begin{lemma}
\label{unl-por}
If $A(z)$ and $A(\hat{z})$ are distinct types
then $A(z)$ is contained in a connected component 
of $\S \setminus A(\hat{z})$.
\end{lemma}

\begin{definition}
\label{unl-def}
We say that two subsets $A, A^{\prime} \subset \S$
are {\bf unlinked} if and only if $A$ is contained in a connected
component of $\S \setminus A^{\prime}$.
\end{definition}

While types live in $\S$, external rays are contained
in the complex plane $\C$. It is convenient 
to have both objects in the same topological space:

\begin{definition}
The circled plane $\copyright$ is the closed topological
disk obtained by adding to $\C$ a circle of points 
$\lim_{r \rightarrow \infty} re^{2 \pi i t}$ at infinity.
The boundary $\partial \copyright$ is canonically identified
with $\S$. 
\end{definition}

Thus, a type $A(z)$ can be  considered as a subset of $\S \cong \partial
\copyright$ and the external rays landing at $z$ are arcs
that join $A(z) \subset \partial \copyright$ with $z$.  

Now we proceed to prove the Lemmas stated above.
But before, let us fix the standard orientation in $\S$ and
use interval notation accordingly with the agreement that
the interval $(t,t)$ represents the circle $\S$ with the point $t$ removed. 

\smallskip
\noindent
{\sc Proof of Lemma~\ref{inv-por}:}
If $t \in A(z)$ then the external ray $R^t_f$ lands
at $z$. Continuity of $f$ plus the fact that $f(R^t_f) = R^{dt}_f$
assures that $dt \in A(f(z))$. Conversely, if $s \in A(f(z))$
then, locally around $z$, the preimage of $R^s_f$ is formed by $k$
arcs. Each of these arcs must belong to a smooth external ray
because in the definition of $A(z)$ we assume that all the
external rays landing at $z$ are smooth. 
\hfill $\Box$

\smallskip
\noindent
{\sc Proof of Lemma~\ref{mon-por}:} 
Since $f$ is locally orientation preserving around
$z$ it must preserve the cyclic order of the rays 
landing at $z$.
\hfill $\Box$

\smallskip
\noindent
{\sc Proof of Lemma~\ref{unl-por}:}
By contradiction, assume that $\{\hat{s}, \hat{s}^{\prime} \} \subset 
A(\hat{z})$ and $\{s, s^{\prime} \} \subset A(z)$ are
such that $\hat{s} \in (s,s^{\prime})$ and $\hat{s}^{\prime} \in
(s^{\prime},s)$. Then the rays $R^s_f$, $R^{s^\prime}_f$ together
with $z$ chop the complex plane $\C$ into two connected components. One which contains
$R^{\hat{s}}_f$ and another which contains $R^{\hat{s}^{\prime}}_f$.
Thus $\hat{z}$ lies in two different sets. Contradiction.
\hfill $\Box$.

\section{Sectors}
\label{fin-por}

\noindent \indent
We want to count the number of external 
rays that participate in some  types.
Following Goldberg and Milnor~\cite{goldberg-93,milnor-95},
several counting problems can be tackled
by a detailed study of the partitions 
of $\C$ and $\S$ which arise from a given 
type. 
To obtain useful information we work 
 under the assumption that there are finitely many
 elements participating in a given type $A(z)$.
Although we have not proved that almost all
 types are finite this will follow from
Theorems~\ref{land-th},~\ref{dis-lan-th},~\ref{wop-th}. 
More precisely, the only types 
 that have a chance of being infinite are
 types of Cremer points and pre-Cremer points.
But it is not known if there exists a Cremer
 point with a ray landing at it.

\begin{definition} 
Let $A(z)$ be a type with finitely many elements.
A connected component of 
$$\C \setminus \bigcup_{t \in A(z)} \overline{R^t_f}$$
is called a {\bf sector}  with basepoint $z$.
A sector $S$ lies in a connected component
of 
$$\copyright \setminus \bigcup_{t \in A(z)} \overline{R^t_f}$$
which intersects $\partial \copyright \cong \S$
in 
an open interval $\pi_{\infty}S \subset \S$.
We say that
the length of $\pi_{\infty}S$ is 
the {\bf angular length} $\alpha (S)$ of $S$.
\end{definition}

For us the circle $\S$ has total length one.
Note that, 
each connected component of 
$\S \setminus A(z)$ corresponds to 
$\pi_\infty S$  for some sector 
$S$ based at $z$.  

\begin{figure}[hbt]
\centerline{
\psfig{file=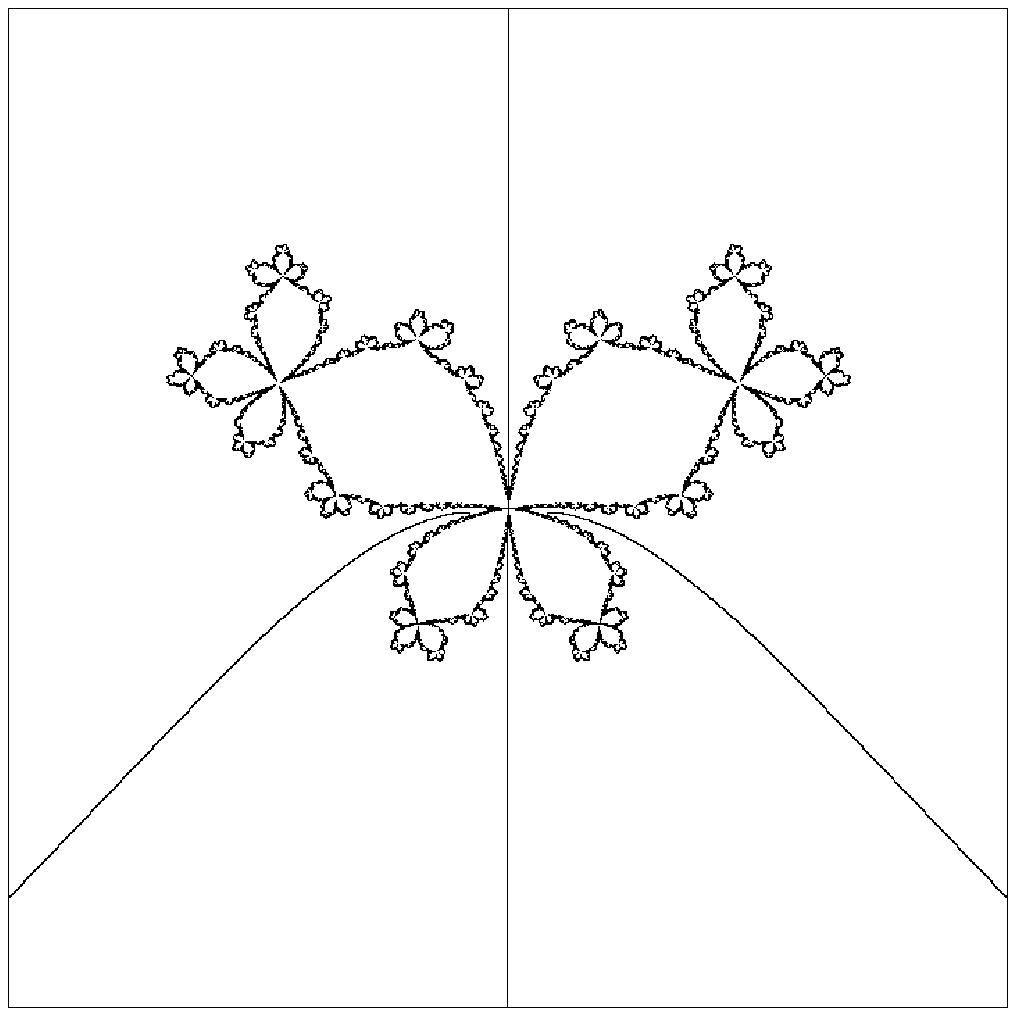,width=170pt}
\hspace{1cm}
\raisebox{0cm}{\psfig{file=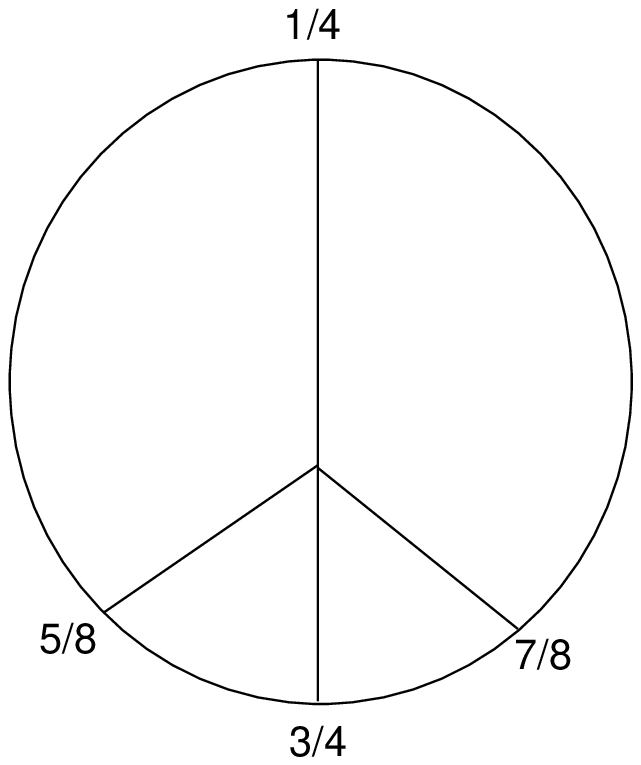,width=140pt}}
}
\caption{At the left, the Julia set of the
 cubic polynomial $f(z) = -1.1z - i z^2 + z^3$
and the external rays landing at the repelling fixed point $0$. 
At the right, the diagram of the type $A(0) = \{1/4,5/8,3/4,7/8\} $.}
\label{anti-butterfly-diagram}
\end{figure}

\smallskip
``Diagrams'' will help us  illustrate,
in the closed unit disk, the partitions 
which arise from a type with finite cardinality.
In the circled plane $\copyright$ consider the union $\Gamma$
of the external rays landing at $z$, the type $A(z) \subset \partial 
\copyright$, and the Julia set element $z$. Let
$h :\copyright \rightarrow \overline{\D}$
be a homeomorphism that fixes the points in 
the circle $\S \cong \partial \copyright \cong \partial \overline{\D}$.
The image $\diag(A(z))$ of $\Gamma$ under $h$ is a {\bf diagram}
of $A(z)$ (See figure~\ref{anti-butterfly-diagram}).

For example, a diagram of $A(z) =\{t_1 , \dots , t_n \}$
can be obtained as follows ($n \geq 2$).
Denote by  
$\zeta$  the center of gravity 
of $\{ e^{ 2 \pi i t_1}, \dots , e^{ 2 \pi i t_n} \}$ and
draw $n$ line segments in $\overline{\D}$ joining $\zeta$ to 
$ e^{ 2 \pi i t_1}, \dots , e^{ 2 \pi i t_n}$.
The resulting graph  $\diag(A(z))$ is a  diagram of $A(z)$.

\smallskip
A first question is to establish how many
critical points or values does a sector 
contain.

\begin{definition}
Let $S$ be a sector. We say that the
{\bf critical weight} $w(S)$ is the number of critical 
points (counting multiplicity) of $f$ contained in 
the open set $S$.
The {\bf critical value weight} $v(S)$ is the number of critical
values of $f$ contained in the open set $S$.  
\end{definition}

\smallskip
In order to detect the presence of critical points 
and critical values in a given sector we 
have to understand how sectors behave under 
iterations of $f$.
Although the global image under $f$ of a sector 
based at $z$ is not necessarily
a sector based at $f(z)$,
locally around $z$ sectors 
map to sectors.

\begin{definition}[Sector map]
\label{tau-d}
For a type $A(z)$ with finite cardinality,
we define a  map $\tau$
which assigns to each sector based at $z$ a 
sector based at $f(z)$ as follows.
Given a sector $S$ based at $z$ 
let $\tau (S)$ be the unique sector based
at $f(z)$ such that $f(S \cap V) \subset \tau(S)$
for some neighborhood $V$ of $z$.
We call $\tau$ the {\bf sector map} at $z$. 
In general, 
for an orbit $\orb$ we introduce as above
the {\bf sector map} $\tau$ at $\orb$ that
takes sectors based at the points of $\orb$
to sectors based at the points of  $f(\orb)$.
\end{definition}

\begin{figure}[htb]
\centerline{
\hspace{-122pt}
\psfig{file=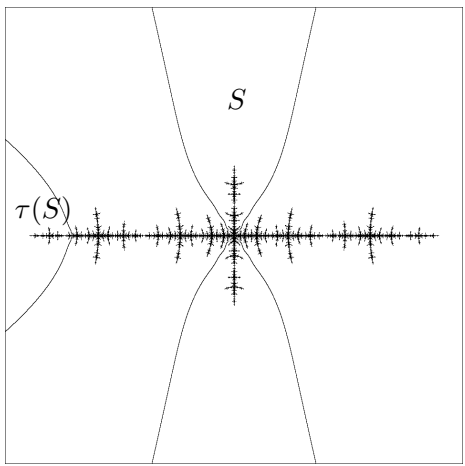,width=174pt}
\hspace{-97pt}
\raisebox{136pt}{$S$}
\hspace{-95pt}
\raisebox{95pt}{$\tau(S)$}
}
\caption{ This figure illustrates
Definition~\ref{tau-d} and Lemma~\ref{op-l} when the sector $S$
is based at the critical point $0$  of the 
Coullet-Feigenbaum-Tresser quadratic polynomial
$z \mapsto z^2 - 1.4101155..$. 
The sector $S$ is bounded
by the external rays with arguments $t_1 = 0.206227..$
and $t_2 = 0.293773..$. 
The sector $\tau(S)$ is bounded
by the external rays with arguments $2t_1 = 0.412454..$
and $2t_2 = 0.587546..$.}
\label{feig2-f}
\end{figure}

It is convenient to understand the action of the sector map
in the circle at infinity. 
If $S$ is a sector based at
$z$ and ``bounded" by the external rays 
with arguments $t_1$ and $t_2$ then, in a neighbourhood of
$z$, the sector $S$ maps to the sector ``bounded"
by $dt_1$ and $dt_2$ (see Figure~\ref{feig2-f}):

\begin{lemma}
\label{op-l}
If $S$ is  a sector such that $\pi_{\infty} S = (t_1,t_2)$ then
$\pi_{\infty} \tau(S) = (dt_1,dt_2)$. 
\end{lemma}

\smallskip
\noindent
{\bf Remark}: From the Lemma above, it follows that
if $(t_1,t_2)$ is a connected component of $\S \setminus A(z)$
then $(dt_1, dt_2)$ is a connected component of $\S \setminus A(f(z))$.

\smallskip
\noindent
{\sc Proof:} 
Pick a neighbourhood $V$ of $z$ such that 
$f(S \cap V) \subset \tau(S)$.
Consider the graph $\Gamma$ formed by the union of the
external rays landing at $z$ and the point $z$. 
If $\pi_{\infty} S = (t_1, t_2)$ then there exists
a connected component $P$ of $\C \setminus f^{-1}(f(\Gamma))$
which contains $R^{t_1 + \epsilon}_f$ for $\epsilon$ small
enough. Now $P \subset S$ and the boundary of $P$ contains
the rays $R^{t_1}_f$ and $R^{t_2}_f$. Moreover, we may assume
that $R^{t_1 + \epsilon}_f \cap V \neq \emptyset$. 
Hence $f_{|P}$ maps $P$  onto a connected component
of $\C \setminus f(\Gamma)$ which is the sector $\tau(S)$ based
at $f(z)$. This sector $\tau(S)$ has in its boundary $R^{dt_1}_f, R^{dt_2}_f$
and contains $R^{dt_1 + d\epsilon}_f$. It follows that
$\pi_{\infty}\tau(S) = (dt_1, dt_2)$.
\hfill $\Box$

\medskip
Following Goldberg and Milnor (see~\cite{goldberg-93} Lemma 2.5 and
Remark 2.6) we state the basic relations between the maps and quantities
introduced above:

\begin{lemma}[Properties]
\label{pro-por}
Let $S$ be a sector of a type $A(z)$ with finitely many
elements, then:

(a) $w(S)$ is the largest integer strictly less than  $d \alpha (S)$.

(b) $\alpha(\tau(S)) = d \alpha(S) - w(S)$.

(c) If $w(S) > 0$ then $v(\tau(S))>0$.

(d) If $\alpha ( \tau(S) ) \leq \alpha (S)$ then $v(\tau(S)) >0$. 

\end{lemma}

\noindent
{\sc Proof:}
For simplicity let us assume that $d \alpha(S)$
is not an integer. The general proof is a 
small variation of the one below.
Let $n$ be the largest integer strictly
less than $d \alpha(S)$.
Consider the loop $\gamma$ in the circled 
plane $\copyright$ that goes from  $t_1$ to $t_2$ 
along the interval $[t_1, t_2]  \subset \partial \copyright$,
it continues along the ray $R^{t_2}_f$ until it reaches $z$
and it goes back to $t_1$ along the ray $R^{t_1}_f$.
Now $f$ acts in $\gamma$ taking $t_1$
to $dt_1$ then it goes $n$ times around
the circle $\partial \copyright$ up to $dt_2$, afterwards
it goes to $z$ along $R^{dt_2}_f$ and back up
to $dt_1$ along $R^{dt_1}_f$. 
Push $\gamma$ to a smooth path $\tilde{\gamma} \subset \C$
and notice that the winding number of the tangent vector
to $\tilde{\gamma}$ around zero is $n+1$. By the Argument
Principle it follows that $f^{\prime}$ has $n$ zeros
in the region enclosed by $\tilde{\gamma}$. Hence,
$w(S) = n$ and (a) of the Lemma follows.
Part (b) is a direct consequence
of (a) and the previous Lemma.

For (c), if no critical value lies in $\tau(S)$ then
consider a branch of the inverse map $f^{-1}$
which takes $\tau(S)$ to $S$. It follows that
$S$ cannot contain critical points of $f$.
Now  (b) and (c) imply (d).
\hfill $\Box$

\smallskip
Observe that part (d) of the previous Lemma says that if a sector $S$
decreases 
in angular length then its ``image'' $\tau (S)$
contains a critical value. 

\begin{figure}[hbt]
\centerline{\psfig{file=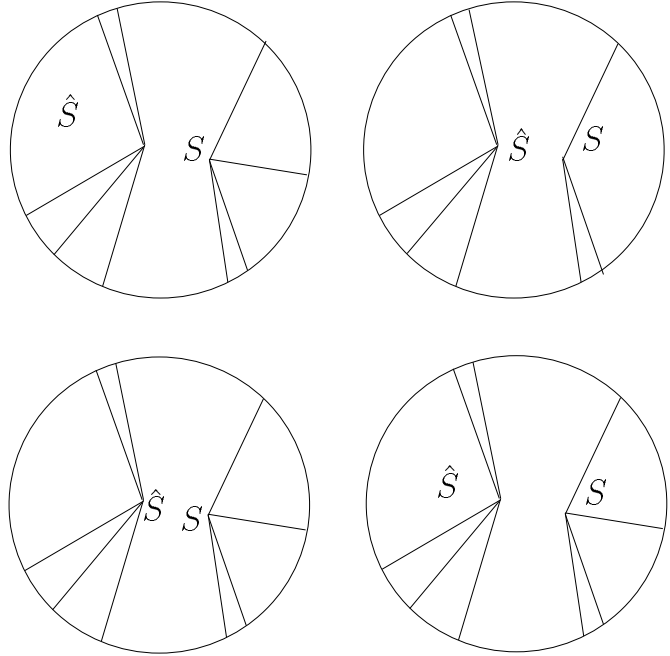,width=250pt,clip=}}
% \vskip -2.0cm
\caption{The four possibilities for the relative position of 
two sectors based at different points}
\label{two-sectors}
\end{figure}

\medskip
Sectors of distinct types are organized in an ``almost"
nested or disjoint fashion. 
As illustrated in figure~\ref{two-sectors},
 there are exactly four alternatives for the relative position
of two sectors based at different points. That is, two sector $S$ and
$\hat{S}$ are either nested or disjoint or each sector contains
the complement of the other.

 For later reference, 
let us record several immediate consequences of this picture in 
the three Lemmas below:

\begin{lemma}
\label{nes-dis2}
Let $S$ and $ \hat{S} $ be  sectors of distinct types
of finite cardinality.
Then one and only one of the following holds:

(i) $S \cap \hat{S} \neq \emptyset$,
 
(ii) $S \subset \hat{S}$,

(iii) $ \hat{S} \subset S$,

(iv) $\C \setminus S \subset \hat{S}$ and $\C \setminus \hat{S} \subset S$.
\end{lemma}

\begin{lemma}
\label{nes-dis0}
Let $S$ and $ \hat{S} $ be  sectors of distinct finite types.
Then:

(a) $S \cap \hat{S} \neq \emptyset$ if and only if 
$\pi_{\infty}S \cap \pi_{\infty}\hat{S} \neq 
\emptyset$.

(b) $S \subset \hat{S}$ if and only if 
$\pi_{\infty}S \subset \pi_{\infty}\hat{S}$.
\end{lemma}

\begin{lemma}
\label{nes-dis1}
Consider two distinct finite types $A(z)$ and $A(\hat{z})$  and 
let $S$ (resp. $ \hat{S} $) be a sector with basepoint
$z$ (resp. $\hat{z}$) then:

(a) If $S \cap \hat{S} \neq \emptyset$ then 
$\hat{z} \in S$ or $z \in \hat{S}$. 

(b) If $\hat{z} \in S$ then
 $S$ contains all but one of the sectors based at $\hat{z}$

(c) If $\hat{z} \notin S$ then  $S$ is contained in exactly one sector 
based at $\hat{z}$. 

\end{lemma}

\section{Periodic Orbit Portraits}
\label{pop}
\noindent \indent
In this section we study types of periodic points
which are the landing point of smooth periodic rays. We
give an upper bound on the number of cycles of rays
that can land at a periodic orbit.
The results here are an immediate generalization of
the ones obtained by Milnor for quadratic polynomials
(see~\cite{milnor-95}). 
In the next Section we will apply similar ideas to
give to find an upper bound on the number of external
rays landing at a point with infinite forward orbit.

We assign to a periodic orbit portrait a {\bf rotation number} as follows.
Let $\orb = \{ z_0,\brkOK f(z_0),\brkOK \dots,\brkOK  f^{\circ (p-1)} (z_0) \}
 \subset J(f)$ 
be a periodic cycle of period $p$
with portrait $A(\orb)$ formed by periodic arguments.
That is, $\orb$ is a parabolic or repelling cycle and each
point in $\orb$ is the landing point of smooth periodic 
rays of the same period.
By Lemma~\ref{mon-por} the return map
$$ {{m_d}^{\circ p}_|}_{A (z_0)} : A(z_0)
\rightarrow A(z_0)$$
is cyclic order preserving.
Thus, the return map has a well defined rotation number 
$\mbox{rot} A (\orb) \in \QS$ which does not depend on the choice of
$z_0 \in \orb$.  The rotation number of $A(\orb)$ 
is also called the {\bf combinatorial rotation number} of $\orb$.

The {\bf number of periodic cycles} of $A (\orb)$ is
the number of cycles of $m_d$ that participate in
$A(z_0) \cup \cdots \cup A(f^{\circ (p-1)}(z_0))$.
For a quadratic polynomial, Milnor~\cite{milnor-95} showed that
the number of cycles of a periodic
orbit portrait $A (\orb)$  is at most $2$.
Moreover, if the number of cycles is $2$ then
$A (\orb)$ has zero rotation number.
We generalize this result:

\begin{theorem}
\label{pop-th1}
The number of cycles of $A (\orb)$ is at most
$d$. Moreover, if the number of cycles is 
$d$ then $A(\orb)$ has zero rotation number.
\end{theorem}

The bounds above are obtained by showing that 
the number of cycles of a portrait gives rise to a lower bound
on the number of critical values of the polynomial
in question. 
Hence, we prefer to state the result as
follows:
 
\medskip
\begin{theorem}
\label{pop-th2}
If  $f$ has exactly  $k$ distinct critical values then
the number of cycles of $A ( \orb )$
is at most $k+1$. Moreover, if ${A (\orb)}$ has $k+1$
cycles then $A(\orb)$ has zero rotation number.
\end{theorem}

\noindent
{\bf Remark:} The bounds on the number of cycles
are sharp. In fact, every parabolic periodic
orbit $\orb$ with $d-1$ immediate basins and
multiplier distinct from $1$ has exactly
$d-1$ cycles of rays participating
in $A(\orb)$. In this case $A(\orb)$ has nonzero rotation number
(also compare with Figure~\ref{anti-butterfly-diagram}).

Figure~\ref{3-cycle-fig} shows a cubic periodic orbit portrait
with $3$ cycles and zero rotation number.

\smallskip
Notice that the sector map $\tau$
at $\orb$ is a well defined permutation of the sectors
based at $\orb$.
Observe that the number of cycles of sectors
under $\tau$ coincides with the number of 
cycles of external rays participating in $A (\orb)$.

\smallskip
The following Lemma (see \cite{milnor-95}) shows that the smallest sector
in a cycle contains a critical value:

\begin{lemma}
\label{pop-pro}
If $\alpha(S) = \min \{ \alpha (\tau^{\circ n}(S)) : n \in {\N} \}$
then $v(S)>0$.
\end{lemma}

\noindent
{\sc Proof:}
Consider a sector $S$ with minimal angular size in its
cycle. If $v(S) = 0$ then Lemma~\ref{pro-por}~(c)
  shows that $w ( \tau^{-1} (S)) =0$. By
Lemma~\ref{pro-por}~(b) we have that $ \alpha( \tau^{-1}(S)) = \alpha(S) /d$,
which contradicts minimality of the angular length.
\hfill $\Box$

\smallskip

\noindent
{\sc Proof of Theorem~\ref{pop-th2}:} 
By contradiction, suppose 
that  $\mbox{rot}{A (\orb)} \neq 0$ and that  
 ${A ( \orb)}$ has more than  $k$ cycles.
Select $k+1$ sectors $S_1, \dots , S_{k+1}$ 
such that:

(a) If $\tau^{ \circ n} (S_i) = S_j $ then $i = j$.
(i.e. $S_1, \dots , S_{k+1}$  belong to different cycles of sectors).

(b) The angular length of $S_i$ is minimal in its cycle
of sectors.

By  Lemma~\ref{pop-pro}
we have that  $v(S_i) > 0$. Thus, in order to obtain a
contradiction it is enough to show that $S_1, \dots , S_{k+1}$ are
pairwise disjoint.

If $S_i \cap S_j \neq \emptyset$ then 
$S_i$ contains the basepoint $z$ of  $S_j$ or  $S_i$ contains the 
basepoint $w$ of $S_j$
(Lemma~\ref{nes-dis1} (a)). In the first case $S_i$ contains all the sectors
based at $z$ with the exception of one
(Lemma~\ref{nes-dis1} (b)). Since the cycle of $S_i$ has 
at least 2 sectors based at $z$ (nonzero rotation number), it 
follows that $S_i$ properly contains a sector in its cycle. This
contradicts (b), i.e. minimality of the angular length. The  same reasoning
gives a  contradiction in the second case.

\smallskip
Now suppose that  $\mbox{ rot}{A (\orb)} = 0$ and that ${A(\orb)}$ has 
more than $k+1$ cycles. Consider a minimal angular length  sector in each 
cycle of sectors  and select  $k+1$ amongst them with smaller angular 
length. Denote these sectors by $S_1 \dots S_{k+1}$. Again we obtain
a contradiction by showing that they are pairwise disjoint.

If $S_i \cap S_j \neq \emptyset$ then 
$S_i$ contains the basepoint $z$ of  $S_j$ or  $S_i$ contains the 
basepoint $w$ of $S_j$. In the first case, $S_i$ contains all the sectors
based at $z$ with the exception of one, which has to be in the cycle of
$S_i$. Hence $S_i$ properly contains  at least $k+1$
sectors of different cycles, this contradicts
the choice of $S_i$. The second case is identical.
\hfill $\Box$

\section{Wandering Orbit Portraits}
\label{wop}
\noindent \indent
A priori we do not know that the
type $A(z)$ of a point $z$ with infinite orbit has finite
cardinality.
 In order to 
apply the results obtained in Section~\ref{fin-por} for  types
with finitely many elements 
we restrict our attention
to finite subsets of $A(z)$. Accordingly
we restrict to finite subsets along the 
forward orbit of $z$:

\newpage
\begin{definition}
Consider an infinite orbit $\orb$ that
does not contain critical values.
We say that 
$$A^* (\orb) = \{ A^* (z) : z \in \orb \}$$
is an {\bf orbit sub-portrait} of $A (\orb )$ if:

 $A^* (z) \subset A (z)$,

 $A^* (z)$  is finite and

 $m_d ( A^* (z) ) = A^* ( f(z))$.
\end{definition}

For orbit sub-portraits we can introduce
sectors, angular length and the sector map $\tau$ 
just as we did in Section~\ref{fin-por}.
It is not difficult to check that the results obtained
in Section~\ref{fin-por} remain valid for sub-portraits.

Observe that if $\orb = \{ z , f(z) , \dots \}$ does
not contain critical values and $A^* (z)$ is a 
finite subset of $A(z)$ then $\{ A^* (z) , m_d(A^*(z)) , \dots \}$
is an orbit sub-portrait of $A(\orb)$.

\smallskip
\noindent
{\bf Remark:} In the definition of sub-portraits 
we avoid orbits containing a critical value
for two reasons. The first one is that we exclude the
special case in which a critical value does not
belong to any of the sectors based at a given point.
The second reason is that since an orbit without critical
values is also free of critical points we have that
the cardinality of $A^* (z)$ is independent of $z \in \orb$.   

\medskip
\noindent
{\sc Proof of Theorem~\ref{wop-th}: }
Given $z \in J(f)$, as in the statement of the Theorem, 
pick $N$ such that 
the forward orbit $\orb$ of $z_0 = f^{\circ N}(z)$
does not contain a critical value. 
First we show that $A(z_0)$ contains at most $d$ elements.

Consider an   orbit
sub-portrait ${A^*}(\orb) $  and 
assume that the cardinality of $A^* (z_0)$ is $d+1$.
After some work we obtain a contradiction.

 Let $z_n = f^{\circ n} (z_0)$
and enumerate the sectors of $A^*(z_n)$ based at $z_n$ 
by $S^1 ( n ) , \dots , S^{d+1} (n)$ according to their angular 
length:$$\alpha ( S^1 (n)) \leq \alpha(S^2 (n)) \leq
\cdots \leq \alpha(S^{d+1} (n)).$$

\smallskip
Intuitively we interpret the angular length of a sector
as its size. 
Under the sector map $\tau$,
sectors of angular length beneath $1/d$ increase
their size.
In contrast,   
for $n$ large, at most 2 sectors based at $z_n$
are not arbitrarily small:

\smallskip
\noindent
{\it Claim 1:} $\lim_{n \rightarrow \infty} \alpha(S^{d-1} (n)) = 0$.
 
\noindent
{\it Proof of Claim 1:} 
By contradiction, suppose that 
$$\lim \sup \alpha(S^{d-1} (n)) = a >0$$
and let $n_k$ be a subsequence such that:
$$\frac{2 a}{3} < \alpha(S^{d-1}(n_k)) < \frac{4 a}{3}.$$
The sectors $S^{d-1} (n_k)$ cannot be disjoint because
otherwise 
we would have infinitely many
disjoint intervals
of length greater than $2 a /3$ contained in $\S$ (Lemma~\ref{nes-dis0}).
Let $n_{k_0}$ and $n_{k_1}$ be such that 
$$S^{d-1} (n_{k_0}) \cap S^{d-1} (n_{k_1}) \neq \emptyset.$$
From
Lemma~\ref{nes-dis1}~(a) we conclude that   
$$z_{n_{k_0}} \in S^{d-1} (n_{k_1}) \mbox{ or }  
z_{n_{k_1}} \in S^{d-1} (n_{k_0}). $$
Without loss of generality,
$z_{n_{k_0}} \in S^{d-1} (n_{k_1})$. This implies that 
all the sectors based at $z_{n_{k_0}}$ with the exception of one 
  are contained in $S^{d-1} (n_{k_1})$.
Since there are 3 sectors
$$S^{d-1} (n_{k_0}),S^{d} (n_{k_0}),S^{d+1} (n_{k_0})$$
based at $z_{n_{k_0}}$  of angular length greater than $2 a /3$
it follows that at least 2 of these must be contained in 
$S^{d-1} (n_{k_1})$.
Thus, the angular length of $S^{d-1} (n_{k_1}) $ 
is greater than $ 4 a / 3$
which is impossible.

\smallskip 
For our purposes we do not distinguish between critical 
values that lie in the same sector based at $z_n$ for
all $n$.
That is, critical values $v$ and $v^\prime$ 
such that $v \in S^k (n)$ if and only if 
$v^\prime \in S^k(n)$ are regarded as ONE
critical value of $f$.
With this in mind, for each critical value $v$ let 
$$\delta(v) = \inf \{ \alpha (S^{k}(n)) : v \in S^k (n) \}.$$

Observe that $\delta(v)=0$ if and only if $v$ is contained in an
arbitrarily small sector.

\smallskip
Loosely speaking, 
we want to show that there
is a correspondence between 
sectors that become arbitrarily 
small and critical values $v$ such that $\delta (v) =0$.
In order to establish this correspondence we need
to ``isolate'' each critical value $v$ such that
$\delta (v) = 0$ from the rest of the critical 
values of $f$. 

Let  
$$\epsilon = \min \{ \delta(v) : \delta(v) \neq 0 \} \cup \{1/d \}.$$

\smallskip
\noindent
{\it Claim 2:} For each critical 
value $v$ such that $\delta(v) = 0$
there exists $n(v)$ such that:

(a) The sector $S(n(v))$ based at 
$z_{n(v)}$ containing $v$ has angular
length $\alpha (S (n(v))) < \epsilon$.

(b) $v^\prime \notin S(n(v))$
for all critical values $v^\prime \neq v$.

(c) $S(n(v))\cap S(n(v^\prime)) = \emptyset$
for all critical values $\delta (v^\prime) =0$. 

\noindent
{\it Proof of Claim 2:}
For parts (a) and (b) enumerate by $v_1 \dots v_m$ the critical
values of $f$ distinct from $v$. 
We already identified the critical values that
always belong to the same sector so there exists
sectors $S_{v_1}, \dots ,S_{v_m}$ such that 
$v \in S_{v_k}$ and $v_k \notin S_{v_k}$.
Now the critical value $v$ is contained in arbitrarily
small sectors ($\delta(v) =0$), thus there exists
an integer $n(v)$ such that the sector $S(n(v))$ 
based at $z_{n(v)}$ containing $v$ has angular length:
$$\alpha (S(n(v))) < \min \{ \alpha(S_{v_k}) ,1- \alpha(S_{v_k})
: k=1 , \dots ,m \} \cup \{ \epsilon \}.$$
The sectors $S(n(v))$ and $S_{v_k}$ are not
disjoint because both contain $v$.
By Lemma~\ref{nes-dis2} we know  that 
one of the following holds:
$$ S(n(v)) \subset S_{v_k}, $$
$$ S_{v_k} \subset S(n(v)), $$
$$ \C \setminus S_{v_k} \subset S(n(v)). $$
The upper bound on $\alpha(S(n(v)))$ says that only the first possibility
can hold. It follows that $v_k \notin S(n(v))$ for all $k$.

For part (c), if $S(n(v)) \cap S(n(v^\prime)) \neq \emptyset$
then one of the following holds
$$S(n(v)) \subset S(n(v^\prime))$$
$$S(n(v^\prime)) \subset S(n(v))$$
$$\C \setminus S(n(v^\prime)) \subset S(n(v))$$
Part (b) of this Claim rules out
the first and second  possibility.
The third one implies that
$$\alpha(S(n(v))) > 1 - \alpha(S(n(v^\prime)))
 \geq 1 - \epsilon \geq \epsilon,$$
which contradicts part (a) and finishes the proof of
Claim 2.

\smallskip
In the next claim we start to establish the correspondence
between small sectors and critical values:

\noindent
{\it Claim 3:} There exists $N_0$ such that
$S^1 (N_0)$ contains a critical value and 

(a) $\alpha(S^1(N_0)) \leq \alpha(S^1(n(v)))$ for all
critical values $v$ such that $\delta(v) =0$.

(b) $\alpha(S^1(N_0)) < \epsilon$.

\noindent
{\it Proof of Claim 3:}
From Claim 1 we know that there exists $M$ such
that for all $n \geq M$:

 $\alpha(S^1(n)) \leq \alpha(S(n(v)))$ for all $v$ such that $\delta(v)=0$ and,

 $\alpha(S^1(n)) < \epsilon$.

 Now, for some $N_0  \geq M$,
the sector $S^1 (N_0)$ must contain a critical value,
otherwise $\alpha(S^1(n))$ would be increasing for $n \geq M$.

\smallskip
We think of $\alpha(S^1(N_0))$ as the threshold for a
sector to be considered ``big'' or ``small''. That is,
if a sector has angular length greater (resp. less)
than $\alpha(S^1(N_0))$ is thought as
being ``big'' (resp. ``small''). 
Now we show that in the transition from sectors
based at $z_n$ to the sectors based at $z_{n+1}$
at most one sector that is ``big'' can ``become small''.

\smallskip
\noindent
{\it Claim 4:} If $\alpha (S^j(n)) \geq \alpha (S^1 (N_0))$
then  $\alpha (\tau(S^{j+1}(n+1))) \geq \alpha (S^1 (N_0))$.

\noindent
{\it Proof of Claim 4:} By contradiction, if 
  $\alpha (\tau(S^{j+1}(n+1))) < \alpha (S^1 (N_0))$
then there are at least 2 sectors $S$ and $S^\prime$
based at $z_n$ such that:
$$\alpha(S) \geq \alpha(S^1(N_0)) > \alpha(\tau(S)) \mbox{ and}$$
$$\alpha(S^\prime) \geq \alpha(S^1(N_0)) > \alpha(\tau(S^\prime))$$
Hence,  $\tau(S)$ (resp. $\tau(S^\prime)$) 
contains a critical value $v$ (resp. $v^\prime$).
Since $\alpha(\tau(S)) < \alpha(S^1(N_0)) \leq \alpha(S(n(v)))$
it follows that $\tau(S) \subset S(n(v))$. Similarly
$\tau(S^\prime) \subset S(n(v^\prime))$. This implies 
that the common basepoint $z_{n+1}$ of the sectors 
$\tau(S)$ and $\tau (S^\prime)$ is contained both
in $S(n(v))$ and in $S(n(v^\prime))$ which contradicts
Claim 2 part (c).

\smallskip
For $1 \leq k \leq d-1$, let $N_k$ be the smallest integer
greater than $N_0$ such that $\alpha(S^k(N_k)) < \alpha(S^1(N_0))$.
That is, $z_{N_k}$ is the first iterate after $z_{N_0}$ for
which $k$ of the sectors based at $z_{N_k}$ are ``small''.
Observe that $N_0 \leq N_1 \leq \dots \leq N_{d-1}$
and that the existence of such integers $N_k$ is guaranteed by Claim 1.
We need to show that there are at least two ``big" sectors
based at each $z_{N_k}$:

\noindent
{\it Claim 5:}  $S^d (N_{k}) \geq 
S^1(N_0)$ for $0 \leq k \leq d-1$.

\noindent
{\it Proof of Claim 5:}
For $k=0$ the claim is trivial.
Given $k \geq 1$ we have that $\alpha(S^k (N_k -1)) \geq 
\alpha(S^1(N_0))$ and by the previous Claim 
we conclude that  $\alpha(S^{k+1}(N_k)) \geq 
\alpha(S^1(N_0))$. Since $\alpha(S^d(N_k)) \geq 
\alpha(S^{k+1}(N_k))$ we are done.

\smallskip
For $1 \leq k \leq d-1$, let $l_k$ be such that 
$$\alpha (S^{l_k} (N_k -1)) \geq \alpha (S^1 (N_0)) >
\alpha (\tau (S^{l_k} (N_k -1))).$$

\noindent
{\it Claim 6:} $S^1 (N_0), \tau (S^{l_1} (N_1 -1)),
\dots \tau (S^{l_{d-1}} (N_{d-1} -1))$ are disjoint and each contains
a critical value.

\noindent
{\it Proof of Claim 6:} 
By Claim 4 and Lemma~\ref{pro-por}~(d),
 we know that these sectors contain critical values.
Now we have to show that they are disjoint.
If $ \tau (S^{l_{k_0}} (N_{k_0} -1))
\cap   \tau (S^{l_{k_1}} (N_{k_1} -1)) \neq \emptyset $ then
without loss of generality 
we may assume that $z_{N_{k_0}} \in \tau (S^{l_{k_1}} (N_{k_1} -1))$.
Lemma~\ref{nes-dis1}~(b) implies that all but one of the sectors based at 
$z_{N_{k_0}}$ are contained in $ \tau (S^{l_{k_1}} (N_{k_1} -1)) $.
From Claim 5, there is at least one sector of angular length
greater then $\alpha(S^1(N_1))$ contained in
$\tau (S^{l_{k_1}} (N_{k_1} -1))$. 
That is, 
$$\alpha (\tau (S^{l_k} (N_k -1))) \geq \alpha(S^1(N_0))$$
which is a contradiction. This finishes 
the proof of Claim 6. 

\smallskip
\noindent
{\bf Remark:}
If a critical value $v$ is contained in one of
the sectors\ \   $S^1 (N_0)$,\ \brkOK $\tau (S^{l_1} (N_1 -1)),\ \brkOK
\dots \ \brkOK \tau (S^{l_{d-1}} (N_{d-1} -1))$ then 
$\delta(v) = 0$ (the angular length of each
of these sectors is less than $\epsilon$).

\smallskip
It follows from Claim 6 and the fact that 
polynomials of degree $d$ have at most $d-1$ critical values
that the cardinality of $A(f^{\circ N}(z))$ is at most $d$.

\smallskip
We modify the arguments above in order to show that 
the cardinality of $A(z)$ is at most $2^d$.

Let $b$ be the number of critical values
that are not in the forward orbit of $z$.
First suppose that $b \geq 1$,
and replace $d$ by $b+1$ in all the statements
from the beginning of the proof up to the end of 
Claim 6. That is, suppose that there are
$b+2$ rays landing at $f^{\circ N}(z)$ and
obtain $b+1$ critical values $v$ such that
$\delta(v) =0$. This is a contradiction
because all the critical values in 
the forward orbit of $z$ cannot be contained 
in arbitrarily small sectors, hence there are
at most $b$ critical values $v$ with $\delta(v) =0$.
Now that we know that at most $b+1$ rays
land at $f^{\circ N}(z)$ let $m_1, \dots, m_a$ be the multiplicity
of the critical points in forward orbit of $z$.
Hence, $A(z)$ has cardinality at most 
$$(m_1 + 1) \cdots (m_a+1) \cdot (b+1).$$ 
Since the
sum $(m_1 +1) + \cdots + (m_a+1) \leq d - 1 - b + a$
and $m_k +1 \geq 2$, it is not difficult to
show that $(m_1 + 1) \cdots (m_a+1) \leq 2^{d-1-b}$.
Then 
$$(m_1 + 1) \cdots (m_a+1) \cdot (b+1) \leq 2^{d-1-b}(b+1) \leq 2^d.$$

Now suppose that $b=0$ and replace 
$d$ by $2$, starting at the beginning of
the proof up to the end of Claim 3.
That is, assume that $3$ rays land at $f^{\circ N}(z)$
and obtain a critical value $v$ such that 
$\delta(v)=0$, this is impossible because
all the critical values are in the orbit of
$z$. Hence, the cardinality
of $A(f^{\circ N}(z))$ is at most 2. The product 
of the local degree of $f$ at the critical points is
at most $2^{d-1}$. Therefore,
when $b=0$ we also have that 
$A(z)$ contains  at most $2^d$ elements.
\hfill $\Box$

\newpage
\part*{Chapter 2: The Shift Locus}
\vspace{-0.8cm}
\indent
\section{Introduction}
In parameter space, following
 Branner and Hubbard~\cite{branner-88},
 we work in the set $\param \cong \C^{d-1}$ of
 monic centered polynomials of degree $d$.
Namely, polynomials of the form:
 $$z^d + a_{d-2}z^{d-2} + \cdots + a_0.$$

Parameter space $\param$ is stratified according to how many
 critical points escape to $\infty$. 
One extreme is the connectedness locus $\conec \subset \param$,
 which is
 the set of polynomials $f$ that have connected Julia
 set $J(f)$. 
Equivalently, all the critical points of $f$ are
 non-escaping.  
The other extreme is the shift locus $\shift \subset \param$,
 formed by the polynomials $f$ that have all
 their critical points escaping. 
In this Chapter we prepare ourselves to explore the
 set $\partial \shift \cap \conec$ where these two 
 extremes meet.

The connectedness locus $\conec$ is compact, connected and
 cellular
 (see~\cite{douady-82,branner-88,lavaurs-89}).
For $d \geq 3$, $\conec$ is known not to be locally connected~\cite{lavaurs-89}.
In contrast, the quadratic connectedness locus,
 better known as the Mandelbrot set $\mathcal{M}$, is conjectured
 to be locally connected (see~\cite{orsay-notes}).

\medskip
The dynamics of a polynomial $f$ in the shift locus
 $\shift$  is completely understood. 
In fact, $f$ has a Cantor set as Julia set $J(f)$ and, $f$ acts on
 $J(f)$ as a hyperbolic dynamical system which is topologically
 conjugate to the one sided shift in $d$ symbols.
The shift locus $\shift$ is open, connected and unbounded.
For $d \geq 3$, $\shift$ has a highly non-trivial topology. 
More precisely,
 its fundamental group is infinitely generated~\cite{blanchard-91}.
In contrast, after Douady and Hubbard~\cite{douady-82},
 the quadratic shift locus 
 $\mathcal{S}_2 = \C \setminus \mathcal{M}$ 
 is conformally isomorphic
 to the complement $\CDC$ of the unit disk.

In the dynamical plane, we describe the location
 of points in the Julia set, which is the boundary of the  
 basin of infinity, by means of external rays and prime end 
 impressions.
In parameter space, we introduce objects
 that will allow us to explore the 
 portion of $\partial \shift$ contained in $\conec$.
More precisely, we define  what it means to go from the shift
 locus $\shift$ towards the connectedness locus 
 $\conec$ in a given direction. 
Each direction will be specified by a ``critical portrait''
 and will determine an ``impression'' in the connectedness
 locus $\conec$. 

In Chapter 4, we are going to show that the ``combinatorics''
 of a polynomial $f$ in $\partial \shift \cap \conec$
 is completely determined by the
 ``impression(s)'' to which $f$ belongs, provided that
 $f$ has all its cycles repelling. 

\medskip
For quadratic polynomials, we have a
 dynamically defined conformal isomorphism
 from $\mathcal{S}_2=\C \setminus \mathcal{M}$ onto
 $\CDC$ (see~\cite{douady-82,orsay-notes}).
This map provides us with parameter rays and a dynamical
 parameterization of the prime end impressions of 
 $\mathcal{S}_2$ in $\partial \mathcal{M}$.
For higher degrees, we need to overcome
 the difficulties that stem from the non-trivial
 topology of the shift locus. 
Motivated by Goldberg~\cite{goldberg-94},
 it is better to work with a dense subset
 of $\shift$ where the critical points 
 are easily located by the B\"ottcher coordinates.
That is, the polynomials $f$ such that each 
 critical point of $f$ is ``visible'' from $\infty$, in the sense
 defined below.
Recall that an external radius is a gradient 
 $\mathrm{grad} g_f$ flow line that reaches $\infty$ 
 (see~\ref{preliminar}).

\begin{definition}[Visible Shift Locus]
Consider a  polynomial $f$ which belongs to the shift locus $ \shift$.
We say that $f$ belongs to the {\bf visible
shift locus} $\vshift$ if for each critical point
$c$ of $f$:

(a)  there are exactly $k$ external radii
terminating at $c$, where $k$ is the local degree
of $f$ at $c$;

(b)  the critical value
$f(c)$ belongs to an external radius.
\end{definition}

Our definition has a slight difference
with Goldberg's definition of the ``generic shift locus''.
Thus, although we use a different name, there is a strong
overlap with the ideas found in~\cite{goldberg-94}.

\smallskip
The quadratic shift locus coincides with the quadratic
 visible shift locus. 
In fact, for
 a quadratic polynomial $f(z) = z^2 + a_0$ in the 
 shift locus $\mathcal{S}_2 = \C \setminus \mathcal{M}$
 there are two external radii $R^{*\theta}_f$ and
 $R^{*\theta + 1/2}_f$ which terminate at
 the unique critical point $c =0$. 
Both of these external radii
 map into $R^{*2 \theta}_f$ which contains the critical value 
 $f(c) = a_0$. 
Similarly, a polynomial, of any degree, 
 with a unique escaping critical point always lies
 in the visible shift locus (see Corollary~\ref{rate-c}).

\begin{figure}[htb]
\label{3-pos-f}
\centerline{
\psfig{file=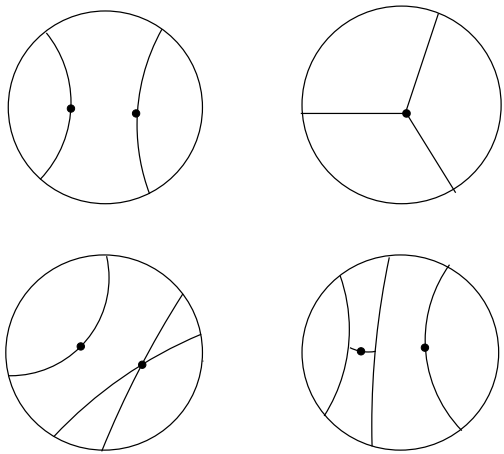,width=187pt}
}
\caption{The possible configurations of gradient flow
lines connecting the critical point(s) of a cubic polynomial
to $\infty$. Only the upper two are allowed in
the visible shift locus $\mathcal{S}^{vis}_3$}
\end{figure}

\smallskip
For a cubic polynomial $f \in \mathcal{ S}_3$
 with two distinct critical points, there are three cases. 
Namely, two external radii might terminate
 at each critical point, or two external radii terminate at one and
 four at the other, or two external radii terminate at
 one and none at the other (see Figure~\ref{3-pos-f}). 
The first case is the
 only one allowed in the visible shift locus $\mathcal{ S}_3$, 
 here, external radii with arguments $\{ \theta_1, \theta_1 + 1/3 \}$
 terminate at one critical point and external radii with
 arguments $\{ \theta_2, \theta_2 + 1/3 \}$ terminate at the other. 
In the second case, one critical point eventually maps to the other.  
In the third case, one critical point lies on external rays
 that bounce off some iterated pre-image of the other 
 critical point. 

\medskip
We keep track of the external radii that terminate 
 at  the critical points:

\begin{definition}
Let $f$ be a polynomial in the visible shift locus $\vshift$
with critical points $c_1 , \dots , c_m$
and $\Theta_i \subset \S$
 be the set formed by the arguments of the external radii 
that terminate at $c_i$. We say that 
$\Theta (f) = \{ \Theta_1 , \dots , \Theta_m \}$
 is the {\bf critical portrait of $f$}.
\end{definition}

The  main properties of $\cp(f)$ are (see Lemma~\ref{cp-l}):

(CP1) For every $j$, $|\Theta_j| \geq 2$ and $ |m_d (\Theta_j)| = 1 $,

(CP2) $\Theta_1 , \dots , \Theta_m $ are pairwise unlinked,

(CP3) $\sum (|\Theta_j| - 1) = d - 1$.

\begin{definition}[Critical Portraits]
A collection $\Theta = \{ \Theta_1 , \dots , \Theta_m \}$
of finite subsets of $\S$ is called a {\bf critical portrait}
of degree $d$ if (CP1), (CP2) and (CP3) hold.
\end{definition}

Critical portraits were introduced by Fisher~\cite{fisher-89} 
to study critically pre-repelling maps and,
since then, widely used in the literature to capture the
location of the critical points 
(e.g.~\cite{bielefield-94, poirier-93, goldberg-93,goldberg-94}).

\smallskip
A result, due to Goldberg~\cite{goldberg-94}, says that for each
 critical portrait $\cp$ there exists a map $f \in \vshift$
 such that $\cp(f) = \cp$.

For $f \in \vshift$, the external radii which 
 terminate at the critical points cut the plane into
 $d$ components. 
In order to capture this situation in the circle at infinity, 
 we define $\cp$-unlinked classes.

\begin{definition}
\label{unl-d}
We say that $t,t^{\prime} \in \S$ are 
$\cp =\{ \cp_1, \dots, \cp_m \}$-{\bf unlinked equivalent}
if $\{t, t^\prime \},\brkOK \cp_1, \dots , \cp_m $ are pairwise
unlinked. 
\end{definition}

Given a degree $d$ critical portrait $\Theta$, there are exactly
$d$  $\Theta$-unlinked classes
$L_1, \dots , L_d$. Moreover, each unlinked class $L_j$ is the union of open
intervals with total length $1/d$.
Intuitively, for polynomials in $\vshift$ close 
to $f$, this partition does not change to much.
Formally, we introduce a topology on the set of
all critical portraits:

\begin{definition}
Let $\ang$ be the set formed 
by all critical portraits endowed with the
{\bf compact-unlinked} topology which is generated by the
subbasis formed by
$$V_{X} = \{ \Theta \in \ang : X \subset L_{\Theta} \}$$
where $X$ is a closed subset of $\S$ and $L_{\cp}$
is a $\cp$-unlinked class.
\end{definition}

\noindent
{\bf Remark:} For ``low'' degrees, a critical portrait $\cp$
 is  uniquely determined by the set $\cp^\cup$ of 
angles which participate in $\cp$ and
the compact-unlinked topology in $\ang$ coincides with the 
Hausdorff topology on subsets $\cp^\cup$ of $\S$. 
For ``high'' degrees, this is not
true. In fact, consider
the degree six critical portraits 
$$\{ \{ 1/12, 1/4 ,7/12, 3/4 \}, \{ 1/3 , 1/2 \} , \{ 5/6 , 0 \} \},$$
$$\{ \{ 1/12, 1/4 \} , \{7/12, 3/4 \}, \{ 1/3 , 1/2 , 5/6 , 0 \} \}.$$

\medskip
The set of quadratic critical portraits $\mathcal{ A}_2$ is
 homeomorphic to $\S$. 
The homeomorphism
 is given by $\{ \theta, \theta + 1/2 \} \mapsto 2 \theta$. 
The set of cubic critical portraits $\mathcal{A}_3$ can be 
 obtained from a M\"obius band $M$ as follows. 
Parameterize the boundary of $M$ by $\S$ and identify $\beta, \beta + 1/3$
 and $\beta + 2/3$. 
The resulting topological space is homeomorphic to $\mathcal{A}_3$.
In general, for $d \geq 3$, 
 $\ang$ is compact and connected but it is not a manifold
 (Lemma~\ref{ang-l}).
The set of critical portraits
 $\ang$ is homeomorphic to the subset $\mathcal{E} \subset 
 \shift$ of polynomials $f$ such that all the critical points $c$
 of $f$ escape to $\infty$ at a fixed rate $\rho = g_f(c)$
 (see Lemma~\ref{ang-l}).  
The topology of $\mathcal{E}$, for cubic polynomials,
 has been previously described by Branner and Hubbard 
 in~\cite{branner-88}.

\smallskip
Now, the topology in the set of critical
portraits  allows us to introduce the impression
of a critical portrait $\cp$ in the connectedness locus
$\conec$:

\begin{definition}
Let $\cp$ be a critical portrait.
We say that $f$ belongs to the {\bf impression}
$I_{\conec}(\cp)$ of the critical portrait $\cp$ if
there exists a sequence of maps $f_n \in \vshift$ 
converging to $f$ such that
the corresponding critical portraits $\cp (f_n)$
converge to $\cp$.
\end{definition}

For quadratic polynomials, $I_{\mathcal{M}}(\{\theta, \theta + 1/2 \})$
is a prime end impression. More precisely, it is
the prime end impression corresponding to $2 \theta$ under
the Douady-Hubbard map 
$\Phi : \C \setminus \mathcal{M} \rightarrow \CDC$.

In order to show that  impressions of critical portraits
are connected and cover 
all of $\partial \shift \cap \conec$ we 
study the basic properties of the map $\Pi$
from  $\vshift$ onto 
the set of critical portraits $\ang$.
The following Theorem asserts that critical portraits
depend continuously on $f \in \vshift$. Also, the set 
$S_{\cp}$ of polynomials in $\vshift$ which share a common 
critical portrait $\cp$ form a sub-manifold of 
$\shift$ parameterized by the escape rates of the critical
points:

\begin{theorem}
\label{coor-th} The subset
$\vshift$ is dense in $\shift$, and the map 
$$\begin{array}{lccl}
 {\Pi}: & \vshift & \rightarrow & \ang \\
                    & f & \mapsto     &
                         \Theta(f)\\
\end{array}$$
is continuous and onto.

Moreover, for any critical portrait
$\cp =\{ \cp_1 , \dots , \cp_m \}$, 
the preimage $S_{\Theta} = \mathbf{\Pi}^{-1}(\cp)$ 
is a $m$-real dimensional manifold. In fact, let
$$\begin{array}{lccl}
     G:   & S_{\cp} & \rightarrow & \R^{m}_{>0} \\
                    & f & \mapsto     &
                                       (g_f(c_1), \dots , g_f(c_m))\\
\end{array}$$
where $c_i$ is the critical point corresponding to $\Theta_i$.
Then $G$ is injective and 
$$ G(S_{\Theta}) = \{ (r_1, \dots ,r_m): d^n \cdot \cp_i \in \cp_j
\Rightarrow d^n r_i > r_j \}.$$
\end{theorem}

The proof of this Theorem appears in Section~\ref{coor-s}.
Afterwards, in Section~\ref{imp-s}, we deduce the following:

\begin{corollary}
The impression $ I_{\conec} (\cp)$ of a critical portrait $\cp$ is a
non empty and connected subset of $\partial \shift \cap \conec$.  
Moreover, 
$$\bigcup_{\cp \in \ang} I_{\conec}(\cp) = \partial \shift \cap \conec.$$
\end{corollary}

\noindent
{\bf Remark:} For quadratic polynomials,
 the sub-manifolds $S_{\cp}$ are the parameter rays introduced
 by Douady and Hubbard in~\cite{orsay-notes}.
For cubic polynomial, if $\cp = \{ \cp_1 , \cp_2 \}$ then
 $S_{\cp}$ is an interval worth
 of Branner-Hubbard ``stretching'' rays (see~\cite{branner-88}).
If $ \cp = \{ \theta, \theta+1/3, \theta+ 2/3 \}$ then
 $S_{\cp}$  is a parameter ray in the parameter plane 
 of the family $z^3 + a_0$.

\section{Dynamical Plane}
\label{dyn-s}
\noindent \indent
For $f$ in the shift 
locus $\shift$,  the Julia set $J(f)$ 
is a measure zero Cantor set
and, on $J(f)$, the map $f$ is 
topologically conjugate to the one sided
shift on $d$ symbols (see~\cite{blanchard-84}). 

In Section~\ref{preliminar}, we summarized some results about
polynomials with disconnected Julia set.
Here, we go into more details about polynomials in 
the visible shift locus $\vshift$. 

\smallskip
In the introduction, we defined $\vshift$ by imposing
conditions on the external radii that terminate at critical
points. Sometimes it is easier to look at the gradient
flow nearby the critical points. Recall that at the critical
points $\mathrm{grad} g_f$ vanishes and that the reduced
basin of infinity $\Omega^*(f)$ is the basin of 
infinity under the gradient flow (Section~\ref{preliminar}).

\begin{lemma}
A polynomial $f \in \shift$ lies in $\vshift$
if and only if for each critical point $c$ of $f$:

(a) there are exactly $k$ local unstable manifolds of
the gradient flow at $c$, where $k$ is the local degree
of $f$ at $c$,

(b) each local unstable manifold of $\mathrm{grad} g_f$
at $c$ is contained in $\Omega^*(f)$.
\end{lemma}

\noindent
{\sc Proof:} It is not difficult to show
that if  $f \in \vshift$ then conditions (a) and (b)
hold. Conversely, (b) implies that each of the local unstable 
manifold of $c$
must lie in an external radius which terminates at 
$c$. These $k$ external radii must map, under $f$, into the same external
radius $R^{*t}_f$. Hence, either $R^{*t}_f$ terminates at
$f(c)$ or $R^{*t}_f$ contains $f(c)$.  In the first case
we have that $f(c)$ must be a pre-critical point. Thus,
the number of unstable manifolds around $c$ would be greater
then $k$, which contradicts (a). Therefore,
$R^{*t}_f$ contains $f(c)$ and $f \in \vshift$.
\hfill
$\Box$

\smallskip
As an immediate consequence we have that:
\begin{corollary}
\label{rate-c}
If $f \in \shift$ is such that all the critical points $c$
of $f$ have the same escape rate $\rho =g_f(c) $ then $f$ lies
in the visible shift locus $\vshift$.
\end{corollary}

In particular, any polynomial of the form $z \mapsto z^d + a_0$
which belongs to  $\shift$ also belongs to $\vshift$.

\smallskip
For the rest of this section, unless otherwise stated,
$f$ is a polynomial in the visible shift locus $\vshift$.
The basic properties of the critical portrait of $f$
are stated below:

\begin{lemma}
\label{cp-l}
Let $f$ be a polynomial in the visible shift locus $\vshift$
with critical points $c_1 , \dots , c_m$
and critical portrait $\Theta (f) = \{ \Theta_1 , \dots , \Theta_m \}$,
where $\Theta_i$ is formed by the arguments of the external radii 
that terminate at $c_i$. Then

(CP1) For every $j$, $|\Theta_j| \geq 2$ and $ |m_d (\Theta_j)| = 1 $,

(CP2) $\Theta_1 , \dots , \Theta_m $ are pairwise unlinked,

(CP3) $\sum (|\Theta_j| - 1) = d - 1$.
\end{lemma}

\smallskip
\noindent
{\sc Proof:}
For (CP1), observe that the external radii that terminate at $c_i$
must map into the unique external radius or ray which contains
the critical value $f(c_i)$.
For (CP2), just notice that external radii are disjoint.
By counting multiplicities (CP3) follows.
\hfill $\Box$

\begin{figure}[htb]
\centerline{
\psfig{file=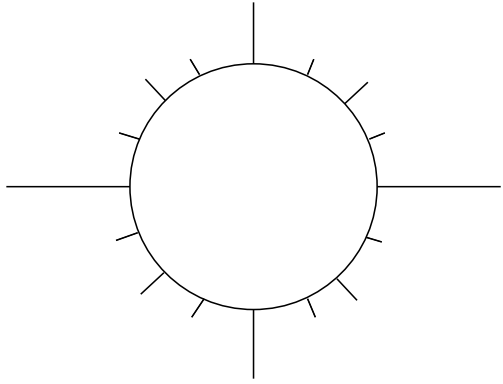,width=187pt}
}
\label{hedgehog}
\caption{Hedgehog}
\end{figure}

\smallskip
From the critical portrait $\cp(f)$
and the escape rates of the critical points of $f$
 we can describe the image $U_f$ of the B\"ottcher map:
$$\phi_f : \Omega^*(f) \rightarrow U_f.$$ 
In fact, assume that  $\cp(f) = \{ \cp_1, \dots, \cp_m \}$
is the critical portrait of $f$ and $g_f(c_1), \dots, g_f(c_m)$
are the escape rates of the corresponding critical points.
Following Levin and Sodin~\cite{levin-91},
for each $\theta \in \cp_i$, let $I_\theta \subset \C \setminus
\D$ be the ``needle'' based at $e^{2 \pi i \theta}$
of height $e^{g_f(c_i)}$:
$$I_{\theta} = [1, e^{g_f(c_i)}]e^{2 \pi i \theta}.$$
Now consider all the iterated preimages of 
$$\bigcup_{\theta \in \cp_1 \cup \dots \cup \cp_m} I_{\theta}$$
under the map $z \mapsto z^d$ to obtain a ``comb''
$C \subset \C \setminus \D$. It follows that
the ``hedgehog'' $\overline{\D} \cup C$ is the complement of 
$U_f$. Equivalently, $U_f = \C \setminus \overline{\D} \cup C$
(see Figure~\ref{hedgehog}).

\smallskip
\noindent
{\bf Example 1:} Consider a quadratic polynomial
$f_v : z \mapsto z^2 + v$ where  $v$ is real and $v > 1/4$.
The external radii with arguments $0$ and $1/2$ terminate
at the critical point $0$ and $\cp(f_v) = \{ 0, 1/2 \}$.
 Say that the escape rate of 
$0$ is $\log r$. Then the ``hedgehog'' for $f_v$ is 
the closed unit disk $\overline{\D}$ union a comb of
needles based at every point of the form $e^{2 \pi i p/2^{n}}$.
For $p$ odd, at each point $p/2^n$ the needle has height
$r^{1/2^{n-1}}$  and at $1 \in \partial \D$
the needle has height $r$.

\smallskip
A point $e^{\rho + 2 \pi i t}$ belongs to $ U_f$
if either $d^n t \notin \cp_1 \cup \dots \cup \cp_m$,
or $d^n t \in \cp_i$ and $d^n \rho > g_f(c_i)$.
Since each critical value $f(c_i)$ belongs to $\Omega^*(f)$
and $\phi_f (f(c_i)) = e^{g_f(c_i) + 2 \pi i d \cp_i}$
we have that:
$$d^n \cdot \cp_i \in \cp_j \Rightarrow d^n g_f(c_i) > g_f(c_j).$$
This explains why in the statement of Theorem~\ref{coor-th}
the image of $G$ is contained in 
$$ \{ (r_1, \dots ,r_m): d^n \cdot \cp_i \in \cp_j
\Rightarrow d^n r_i > r_j \}.$$

\begin{lemma}
\label{subset-l}
In the notation of Theorem~\ref{coor-th}, 
$$G(S_\cp) \subset \{ (r_1, \dots ,r_m): d^n \cdot \cp_i \in \cp_j
\Rightarrow d^n r_i > r_j \}.$$
\end{lemma}

Also note that for 
$$t \notin \Sigma = \bigcup_{n \geq 0} m_d^{-n}(\cp_1 \cup \dots \cup \cp_m)$$
the external ray $R^t_f$ is smooth.
For $t \in \Sigma$ we have two non-smooth external rays
$R^{t^+}_f$ and $R^{t^-}_f$ which bounce off some
pre-critical point(s).

\smallskip
\noindent
{\bf Example 1:} (continued)
 For a quadratic polynomial $f_v: z \mapsto z^2 +v$
where $v > 1/4$, the external rays with arguments
of the form $p/2^n$ eventually map to one of the
fixed non-smooth external rays $R^{0^\pm}_{f_v}$ which contain
the critical point. It follows that the external rays
with argument $p/2^n$ are not smooth.

\smallskip
\noindent
{\bf Example 2:} Consider a cubic polynomial 
$f \in \mathcal{ S}^{vis}_3$ with
critical portrait $\{ \cp_1 = \{ 1/3, 2/3 \},\brkOK
 \cp_2 = \{1/9, 7/9 \} \}$.
Then $g_f(c_2) > g_f(c_1)/3$. The external rays with arguments
of the form $p/ 3^q$, where $p \neq 0$,  are not smooth 
because $m^{\circ q-1}_3 (p/3^q) = 1/3 \mbox{ or }2/3$.

\begin{figure}[htb]
\label{ex3}
\centerline{
\psfig{file=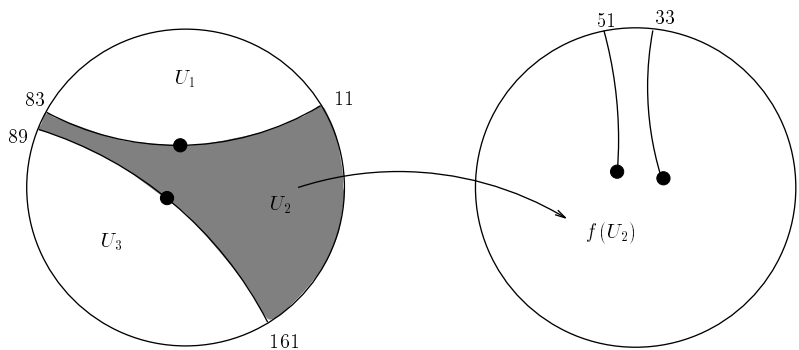,width=306pt}
}
\caption{Schematic picture of the external radii
terminating at the critical points of a cubic polynomial 
$f$ with critical
portrait $\{ \{11/216,83/216\} , \{ 89/216,161/216\} \}$.
Also we illustrate the image of these external radii
and of the region $U_2$.
Units are in $1/216$.}
\end{figure}

\smallskip
The  external radii with arguments in $\cp_1 \cup \dots \cup \cp_m$
together with the critical points chop the complex plane into
$d$ connected components $U_1, \dots, U_d$. 
The boundary
of $U_i$ is formed by pairs of external radii
that terminate at a common critical point.
Each of these pairs is mapped onto an
arc which joins a critical value to $\infty$.
Moreover, $f$ maps $U_i$ homeomorphically onto a slited complex plane
and $\overline{U}_i$ onto $\C$ (see Figure~\ref{ex3}).

In the circle at infinity, each connected
component $U_i$ spans a $\Theta$-unlinked
class $L_i$. Each $\cp$-unlinked class $L_i$
is a finite union of intervals with total length
$1/d$. The boundary points of $L_i$
are mapped two to one by $m_d$ and $L_i$ is mapped
injectively onto its image.

\smallskip
\noindent
{\bf Example 3:} Consider a cubic polynomial $f$
with critical portrait 
$$\cp = \{ \{ 11/216 , 83/216 \}, \{ 89/216 , 161/216 \} \}.$$
The $\cp$-unlinked classes are $L_1 = (11/216,83/216)$,
$$L_2 = (83/216,89/216) \cup (161/216,11/216)$$ and
$L_3 = (89/216,161/216)$. The schematic situation 
is represented in Figure~\ref{ex3}.

\smallskip
Since the Julia set $J(f)$ is a Cantor set, every external
ray lands. 
The symbolic dynamics induced on $J(f)$
 by the connected components $U_1, \dots, U_d$
corresponds to the symbolic dynamics induced on the arguments
of the external rays by the $\cp(f)$-unlinked classes.

\begin{definition}
\label{itin-d}
Given  a critical portrait $\cp$ of degree $d$ with $\cp$-unlinked classes
$L_1, \dots, L_d$,  let
$$\begin{array}{lccl}
\mathrm{itin}^{\pm}_{\Theta} : & \S & \rightarrow & 
        \{ 1, \dots , d \}^{{\N}\cup \{0\}} \\
                    & t  & \mapsto     &
                        ( j_0 , j_1 , \dots )\\
\end{array}$$
if, for each $n \geq 0$, there exists $\epsilon >0$ such that
$(d^n t, d^n t\pm \epsilon) ) \subset L_{j_n}$.
\end{definition}

Now we have the following: 

\begin{lemma}
\label{lan-vshift-l}
Consider  $f$ in the visible shift 
locus $\vshift$ with critical 
portrait $\Theta (f)$. Two external rays
$R^{t^\epsilon}_f$ and $R^{s^\delta}_f$
land at a common point if and only if
$\mathrm{itin}^{\epsilon}_{\Theta}(t) =
\mathrm{itin}^{\delta}_{\Theta}(t)$ where $\epsilon, \delta = \pm$.
\end{lemma}

\begin{figure}[htb]
\label{ex3-c}
\vspace{1cm}
\centerline{
\psfig{file=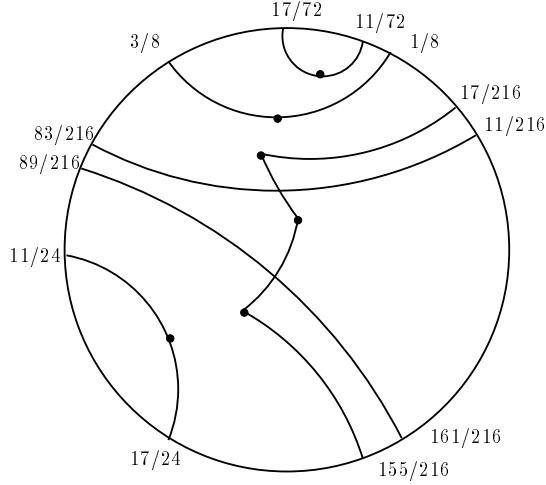,width=204pt}
}
\smallskip
\caption{The pattern in which some external rays land for a 
cubic polynomial $f$ with critical portrait 
$\cp(f) = \{ \{ 11/216 , 83/216 \}, \{ 89/216 , 161/216 \} \}$.
Notice that the rays with arguments that participate in $\cp(f)$
bounce off critical points. Dots represent the landing points}
\end{figure}

Before we prove the Lemma let us discuss an example:

\smallskip
\noindent
{\bf Example 3} (continued): Since
$$\mathrm{itin}^+_{\cp}(t=161/216) = \mathrm{itin}^-_{\cp}(s=11/216) = 213111111...,$$
it follows that $R^{t^+}_f$ and $R^{s^-}_f$ land at a common
point $z$. The external rays with arguments $3t = 17/72$ and
$3s = 11/72$ are smooth and land at the same point $f(z)$.
See Figure~\ref{ex3-c} for a schematic picture which illustrates 
how these and other rays land.

\bigskip
\noindent
{\sc Proof of Lemma~\ref{lan-vshift-l}:}
Consider the forward invariant closed set
formed by the iterates of the external radii
which terminate at critical points:
$$ X = \bigcup_{t \in \cp_1 \cup \dots \cup \cp_m} 
\bigcup_{n \geq 1} f^{\circ n}(\overline{R}^{*t}_f)$$
The inverse image of $\C \setminus X$ has $d$ components
$V_1 \subset U_1, \dots, V_d \subset U_d$. In the Julia set
$J(f)$ each branch of the inverse is a strict contraction
with respect to the hyperbolic metric in $\C \setminus X$.
Thus, a point $z \in J(f)$ is completely determined
by its itinerary in $V_1, \dots, V_d$. Hence, the landing
point of $R^{t^{\delta}}_f$ is completely determined
by $\mathrm{itin}^{\delta}_{\cp}(t)$.
\hfill $\Box$

\bigskip
When no periodic argument $\theta$
participates of $\cp(f)$ all the periodic rays
are smooth. 
In fact, given a periodic
argument $t_1$, we have that
$(1, \infty) e^{2 \pi i t_1} \subset U_f$.
Moreover, the next Lemma, due to Levin and
Sodin~\cite{levin-91}, shows that there exists a definite 
``triangular'' neighbourhood of $(1, \infty) e^{2 \pi i t_1}$
contained in $U_f$. We will need this result in Chapter 4.

\begin{lemma}
\label{ls-l}
Consider $f \in \vshift$ such that $\cp(f) =
 \{ \cp_1, \dots, \cp_m \}$ and let
$$\mu = \max g_f(c)$$
be the maximal escape rate of the critical points.
Consider the exponential map
$$\begin{array}{lccl}
     exp:   & H = \{z = x + iy : x > 0 \} & \rightarrow & \CDC \\
            &  z                       & \mapsto     &
                                              e^z\\
\end{array}$$ 
and let $\tilde{U}_f = exp^{-1}(U_f)$. Let 
$t_1 \in \QS$ be a periodic argument with 
orbit $\{t_1, \dots, t_p \}$ under $m_d$.
Denote by $\delta$ the angular distance
between $\{ t_1, \dots, t_p \}$ and
$\cp_1 \cup \cdots \cup \cp_m$.
Let 
$$ \tilde{V} = 2 \pi i t_1 + 
\{ z \in H :  |\arg (z)| < arctan(\delta/2 \pi \mu) \}.$$
Then $\tilde{V} \subset \tilde{U_f}$.
\end{lemma}

\section{Coordinates}
\label{coor-s}
\noindent \indent
In this section we prove Theorem~\ref{coor-th}.
First we need some facts about how the B\"ottcher map
$\phi_f$ and the Green function $g_f$ depend on $f \in \param \cong \C^{d-1}$.

\begin{lemma} Consider the open sets:
$$\mho = \{(f,z) \in \param \times \C : z \in \Omega(f) \},$$
$$\mho^* = \{(f,z) \in \param \times \C : z \in \Omega^*(f) \}.$$
The Green function:
$$\begin{array}{lccl}
     g:   & \mho & \rightarrow & \R_{>0} \\
          & (f,z)& \mapsto     &
                                       g_f(z)\\
\end{array}$$ 
is real analytic.
The B\"ottcher map:
$$\begin{array}{lccl}
     \phi:   & \mho^* & \rightarrow & \CDC \\
          & (f,z)  & \mapsto   & \phi_f(z)\\
\end{array}$$
is holomorphic.
\end{lemma}

\noindent
{\sc Proof:} For $z$ in a neighbourhood of infinity,
$\phi_f(z)$ depends holomorphically both on $f$ and $z$
(see~\cite{branner-88} I.1).
Since 
$$g_f(z) = \frac{g_f(f^{\circ n}(z))}{d^n} = \frac{\log_{+} |\phi_f(z)|}{d^n}$$
we have that $g$ is real analytic in $\mho$. By continuous dependence
on paramaters of the gradient flow lines, $\mho^*$ is open.
Spreading, along flow lines,
the holomorphic dependence of $\phi_f(z)$, for $z$
near infinity, to the domain $\mho^*$ the Lemma follows.
\hfill
$\Box$

\smallskip
The ``visibility condition'' imposes restrictions 
on the unstable manifold 
of the critical points under the gradient flow. 
To rule out certain situations
we work with broken flow lines.

\begin{definition}
For $[a,b] \subset (0, \infty]$, a broken flow line
is a path:
$$\gamma_f: [a,b] \rightarrow \Omega(f) \cup \{\infty\}$$
such that 
for $r \in (a,b)$
$$g_f(\gamma_f(r)) = r$$
and $\gamma_f(r)$ is either a singularity or $\mathrm{grad}g_f$ is tangent
to $\gamma_f$ at $\gamma_f(r)$.
\end{definition}

\begin{lemma}
\label{broken-l}
Consider a sequence $ \{ f_n \} \subset \param$ 
which converges to $f \in \param$. Let 
$$\gamma_{f_n}: [a_n , b_n] \rightarrow \Omega(f_n) \cup {\infty}$$
 be broken flow lines
such that $a_n \rightarrow a \in (0, \infty]$ and
$b_n \rightarrow b \in (0, \infty]$.
 Then, there exists a
subsequence $\{ \gamma_{f_{n_i}} \}$ which converges
to a broken flow line 
$$\gamma_f : [a, b] \rightarrow \Omega(f) \cup {\infty}.$$
 Here by convergence
we mean that if $s_{n_i} \rightarrow s \in [a,b]$ then 
$\gamma_{f_{n_i}} (s_{n_i}) \rightarrow \gamma_{f}(s)$.
\end{lemma}

\noindent
{\sc Proof:} By passing to 
a subsequence we may assume that $\gamma_f(a_n)
\rightarrow z_0$. There are only finitely many broken flow
lines starting at $z_0$ and ending at the equipotential
$g_f = b$ or at $\infty$
(in the case $b = \infty$).
 By continuous dependence of the gradient flow,
a subsequence of $\{ \gamma_{f_n} \}$ must converge to one of these
broken flow lines.
\hfill
$\Box$

\begin{lemma}
The visible shift locus 
$\vshift$ is dense in shift locus $\shift$.
\end{lemma}

\noindent
{\sc Proof:}
By contradiction, suppose that $\shift \setminus \vshift$
has non-empty interior. Under this assumption, we
restrict to an open set where ``visibility" fails 
in a controlled manner:

\smallskip
\noindent
{\it Claim 1:} There exists an open set $V \subset \shift \setminus
\vshift$ and holomorphic functions $c : V \rightarrow \C$,
$\tilde{c} : V \rightarrow \C$ and
$s : V \rightarrow \C$ such that:

(a) Each $f$ in $V$ has $d-1$ distinct critical points.

(b) $c(f)$ and $\tilde{c}(f)$ are critical points
$f$ and $s(f)$ is a singularity of $g_f$.

(c) There exists a broken flow line of 
$g_f$ from $c(f)$ to $s(f)$.

(d) There exists $k$ such that $f^{\circ k}(s(f)) = \tilde{c}(f)$.

\noindent
{\it Proof of Claim 1:} 
Condition (a) is open and dense and implies
that the critical points depend holomorphically
on $f$. 
There can be only finitely many singularities 
between the slowest escaping critical point
and the fastest one. We can assume that 
these singularities also depend holomorphically
on $f$ in an open dense set of $\shift$. Now since we suppose 
that $\shift \setminus \vshift$ has non-empty interior
there exists an open set $W$ where there is a broken 
flow line between a critical point and a singularity.
Locally in $W$, there are finitely many possible combinations
and each occurs in a closed set (Lemma~\ref{broken-l}).
Hence, there must exist an open set $V \subset W$
such that for $f \in V$ there exists a broken flow line
between a critical point $c(f)$ and a singularity $s(f)$
which depends holomorphically on $f$.
Thus, $s(f)$  must
map onto a critical point $\tilde{c}(f)$ after a fixed 
number of iterates $k$ (i.e. $f^{\circ k}(z) = \tilde{c}(f)$)
and the Claim follows.

\medskip
For $n$ sufficiently large, $f^{\circ n}(c(f))$ is 
close to $\infty$ and  $f^{\circ n}(c(f))$ 
belongs to the domain $\Omega^*(f)$
of the B\"ottcher function $\phi_f$.
Now $f^{\circ n}(s(f))$ also lies in $\Omega^*(f)$, 
closer to $\infty$, along the same external radius which contains 
$f^{\circ n} (c(f))$.
Furthermore, for $m = n -k$ we
have that $f^{\circ n}(s(f)) = f^{\circ m}(\tilde{c}(f))$.
Hence the quotient
$$
\frac{\phi_f \circ f^{\circ m}(\tilde{c}(f))}{\phi_f \circ f^{\circ n}({c}(f))}
                                                \in \R $$
and depends holomorphically on $f \in V$. 
It follows that for some $K > 1$,
$$\phi_f \circ f^{\circ m}(\tilde{c}(f)) = K \phi_f \circ f^{\circ n}({c}(f))$$
for all $f \in V$.

To show that this situation cannot occur we perturb $f_1 \in V$ using
Branner and Hubbard's wringing construction (see~\cite{branner-88}).
We will only need the {\it stretching} part of this construction
that we briefly summarize below. 
For $s > 0$, consider the quasiconformal map
$$\begin{array}{lccl}
     l_s :  & \CDC & \rightarrow & \CDC \\
           & re^{2 \pi i \theta} & \mapsto &
                                       r^s e^{2 \pi i \theta}\\
\end{array}$$
which commutes with $z \mapsto z^d$.
The pull-back $\mu_s = l^*_s \mu_0$ of the standard conformal
structure $\mu_0$ is a Beltrami differential which depends 
smoothly on $s$.

From $\mu_s$ one obtains a conformal structure $\nu_s$ invariant under $f_1$
as follows. Let $R \geq 1$ be large enough so that  
$U_R  = {\phi_{f_1}}^{-1} (\CDC_R)$ is well defined.
We may assume that 
$ f^{\circ n}_1 (c(f_1)) \in U_R.$ Let
$$ \nu_s (z) =  \phi^*_{f_1} \mu_s (z) \mbox{  for  } z \in U_R$$
and extend $\nu_s$ to the basin of
infinite $\Omega(f_1)$ by successive pull-backs of $\mu_s$ under $f_1$.
Finally, let $\nu_s (z) = 0$ for $z \in J(f_1)$. 

Apply the Measurable Riemann Mapping Theorem (~\cite{ahlfors-66}
ch. V) with parameters
to obtain a continuous family of quasiconformal maps
$h_s$ such that  $h^*_s \mu_0 = \nu_s$ where
$h_s$ is normalized to fix $0$, $1$ and  $\infty$. It follows
that $h_s \circ f_1 \circ h^{-1}_s$  is a 
family of polynomials, but a priori we do not know 
if they are monic and centered. Following
Branner and Hubbard, we adjust $h_s$ in order to meet
the required properties:

\smallskip
\noindent
{\it Claim 2:} There exists a continuous 
family $\tilde{h}_s: \C \rightarrow \C$ of quasiconformal
maps such that:

(a) $\tilde{h}_1 (z) = z$ for $z \in \C$,

(b) $f_s = \tilde{h}_s \circ f_1 \circ \tilde{h}^{-1}_s$
is a continuous family of monic centered polynomials.

(c) $\phi_{f_s} (z) = l_s \circ \phi_{f_1} \circ \tilde{h}^{-1}_s (z)$
for $z \in \tilde{h}_s(U_R)$ where $\phi_{f_s}$ is the B\"ottcher map
of $f_s$.

\noindent
\smallskip
{\it Proof of Claim 2:}
With $h_s$ as above we 
have that 
$$h_s \circ f_1 \circ h^{-1}_s 
= a_{d}(s) (z^{d} + a_{d-1}(s) z^{d-1} + \cdots + a_0(s)).$$ 
Notice that  $h_1$ is the identity, hence $a_d(1) =1$
and $a_{d-1} (1) = 0$.

To check that $a_0(s), \dots , a_d(s)$ are 
continuous observe that 
the critical points of $h_s \circ f_1 \circ h^{-1}_s $
vary continuously with $s$ 
because they are the image under $h^{-1}_s$ of the critical points of $f_1$.
The coefficients
$a_1(s), \dots, a_{d-1}(s)$ are continuous functions of
the critical points of  $h_s \circ f_1 \circ h^{-1}_s $
and hence of $s$. 
Since $h_s$ fixes $0$ and $1$ it follows that
$a_0(s)$ and $a_d(s)$ also depend continuously on 
$s$. 

Choose a continuous branch of $a_d(s)^{1/d-1}$
such that $a_d(1)^{1/d-1} =1$ and
let  
$$\tilde{h}_s (z) = a_d(s)^{1/d-1} ( h_s(z) + a_{d-1}(s)/d). $$
Now $f_s = \tilde{h}_s \circ f_1 \circ \tilde{h}^{-1}_s$
is a continuous family of monic centered polynomials.
By construction 
$$ l_s \circ \phi_{f_1} \circ \tilde{h}^{-1}_s : \tilde{h}_s (U_R)
\rightarrow  \CDC_{R^s}$$
is a conformal isomorphism which conjugates 
$f_s$ and $z \mapsto z^d$. Hence, it must be the B\"ottcher
map of $f_s$ up to a $(d-1)st$ root of unity. But for 
$s=1$ we have that 
$l_1 \circ \phi_{f_1} \circ \tilde{h}^{-1}_1 = \phi_{f_1}$.
Thus, by continuity,
$l_s \circ \phi_{f_1} \circ \tilde{h}^{-1}_s$ is tangent to the 
identity at infinity for all $s$.
Uniqueness of $\phi_{f_s}$ finishes the proof of the Claim.

\medskip
Since $\tilde{h}_s$ is a conjugacy, it maps
critical points
to critical points and their iterates also correspond.
In particular,
 $\tilde{h}_s \circ f^{\circ n}_1 (c(f_1)) = f^{\circ n}_s (c(f_s))$
and $\tilde{h}_s \circ f^{\circ m}_1 (\tilde{c}(f_1)) 
= f^{\circ m}_s (\tilde{c}(f_s))$.
After replacing in part (c) of the Claim:
$$| \phi_{f_s} \circ f^{\circ m}_s (\tilde{c}(f_s)) |
  = | \phi_{f_1} \circ f^{\circ m}_1 (\tilde{c}(f_1)) |^s$$
and
$$ | \phi_{f_s} \circ f^{\circ n}_s ({c}(f_s)) |
= | \phi_{f_1} \circ f^{\circ n}_1 ({c}(f_1))|^s$$
which gives us the desired contradiction.
\hfill $\Box$

\smallskip
Recall that the map $\Pi$ assigns to each polynomial
$f \in \vshift$ its critical portrait $\cp(f) \in \ang$.

\begin{lemma}
$\Pi$ is continuous.
\end{lemma}

\noindent
{\sc Proof:} Consider a closed subset $X \subset \S$ 
and the corresponding element $V_X$ of the subbasis that
generates the compact-unlinked topology in $\ang$.
We must show that $\Pi^{-1}(V_X)$ is open or equivalently
that $\vshift \setminus \Pi^{-1}(V_X)$ is closed. 
Take a sequence $\{f_n\} \subset \vshift$ such that 
$f_n \rightarrow f \in \vshift$ 
and $X$ is not contained in a $\cp(f_n)$-unlinked 
class. Thus, there exists
two external radii of $f_n$ with arguments $t_n$ and $t^{\prime}_n$
which terminate at a common critical point $c_n$ of $f_n$ and
$X$ is not contained in a connected component of 
$\S \setminus \{t_n , t^{\prime}_n \}$. By passing to 
a subsequence we may assume that $t_n \rightarrow t$,
$t^{\prime}_n \rightarrow t^{\prime}$ and
$c_n \rightarrow c$ where $c$ is a critical
point of $f$. 
In view of Lemma~\ref{broken-l},
 by passing to a further subsequence, the closure 
of the external radii $\overline{R^{*t_n}_{f_n}} $ converge
to a broken flow line that connects 
a critical point $c$ of $f$ to infinity. 
Near infinity, this broken flow line
coincides with $R^{*t}_f$. Since $f$ lies  
in $\vshift$ the broken flow lines connecting $c$ to $\infty$ are
the closure of external radii. 
Hence, $R^{*t}_f$ terminates at $c$ and,
similarly, $R^{*t^{\prime}}_f$ also terminates at $c$.
In the limit we also have that the closed set $X$ is not contained 
in a connected component of $\S \setminus \{ t, t^{\prime} \}$.
Therefore, $X$ is not contained in a $\Theta(f)$-unlinked class
and $\vshift \setminus \Pi^{-1}(V_X)$ is closed. 
\hfill
$\Box$

\medskip
The fact that $\Pi$ is onto relies on a result of 
Goldberg (see~\cite{goldberg-94} Proposition 3.8):

\begin{proposition}
\label{goldberg-p}
Let $\cp$ be a critical portrait. Then there
exists a polynomial $f \in \vshift$ such that 
$\cp = \cp(f)$.
\end{proposition}

Recall that, given a critical portrait
  $\cp = \{ \cp_1 , \dots , \cp_m \}$,
$G$ assigns $$(g_f(c_1), \dots, g_f(c_m))$$ to each polynomial 
$f$ with critical portrait $\cp$,
where the external radii with arguments in $\cp_i$
terminate at the critical point $c_i$.

\begin{lemma}
\label{inj-l}
G is injective.
\end{lemma}

\noindent
{\sc Proof:} Assume that $f_1$ and $f_2$ are polynomials
in the visible shift locus $\vshift$ with the same critical
portrait $\cp = \{ \cp_1, \dots , \cp_m \}$ and such that
$G(f_1) = G(f_2) = ( \rho_1 , \dots , \rho_m )$.
We must show that $f_1 = f_2$. The idea is to use 
the  ``pull-back argument" to construct a 
quasiconformal conjugacy
$\hat{h}$ between $f_1$ and $f_2$ which is conformal
in the basin of infinity $\Omega(f_1)$. 
Then, we can argue that 
$\hat{h}: \C \rightarrow \C$ is actually conformal because
the Julia set $J(f_1)$ has zero Lebesgue measure
(see~\cite{lehto-73} ch. V.3).

\smallskip
For $i = 1,2$ consider the sets $X_i$ formed by the union of:

(a) The region outside an high enough equipotential:
$$\{ z / g_{f_i} (z) > \rho \}$$
where $\rho = 2d \max \{ \rho_1 , \dots , \rho_m \}$,

(b) The portion of the external radii that run down
from infinity up to a point in the forward orbit of
a critical value:
$$\bigcup_{n \geq 1} f^{\circ n}_i (R^{*t}_{f_i})$$
where $t \in \cp_1 \cup \cdots \cup \cp_m$.

(c) The forward orbit of the critical values.

Observe that $X_i$ is completely contained in
the domain $\Omega^*(f_i)$ of the B\"ottcher map $\phi_{f_i}$.
Also notice that, in part (b), although we take the union
over infinitely many sets, all but finitely of these
are outside the equipotential of level $\rho$.

In $Y_i = f^{-1}_i (X_i) \supset X_i$ the only singularities
of the gradient flow are the critical points.
By analytic continuation along flow lines
of $\mathrm{grad} g_{f_1}$ extend 
$$\phi^{-1}_{f_2} \circ \phi_{f_1} : X_1 \rightarrow X_2$$
to a conformal isomorphism $h_0$ from a 
connected neighbourhood $N(Y_1)$ of $Y_1$ 
onto a neighbourhood  $N(Y_2)$ of $Y_2$.
This is possible because $\cp(f_1) = \cp(f_2)$
and $G(f_1) = G(f_2)$.
In fact, 
$\phi^{-1}_{f_2} \circ \phi_{f_1}$ around a critical
value $f_1(c)$ can be lifted to a map around the corresponding
critical point $c$ in order to agree with the analytic
continuation along the external radii that terminate at
$c$.

The complement of $Y_i$ are $d$ topological disks,
therefore after shrinking
$N(Y_i)$ (if necessary) $h_0$ extends to a $K$-quasiconformal map 
$$\hat{h}_0 : \C \rightarrow \C.$$

So far we have a $K$-quasiconformal map $\hat{h}_0$
which is a conformal conjugacy in $N(Y_1)$ (i.e. $f_2 \circ \hat{h}_0 (z) =
\hat{h}_0 \circ f_1 (z)$ for $z \in N(Y_1)$).
The region  $N(Y_1)$ is connected and contains all the critical
values of $f_1$.
The critical values of $f_1$ are taken onto the critical values of
$f_2$ by $\hat{h}_0$. A similar situation occurs 
with the pre-image  of the critical values. 
It follows that $\hat{h}_0$ lifts to a unique $K$-quasiconformal
map $\hat{h}_1$ which agrees with $\hat{h}_0$ in $N(Y_1)$:

%\begin{equation*}
$$\begin{CD}
\C @>\hat{h}_1>>  \C \\
@Vf_1VV    @VVf_2V \\
\C @>>\hat{h}_0> \C 
\end{CD}$$
%\end{equation*}

Now we have a conjugacy in a larger set:
$$\hat{h}_1 \circ f_1 (z) = f_2 \circ \hat{h}_1 (z) \mbox{ for } z \in 
f^{-1}(N(Y_1))$$
which is, in particular, conformal in $\{ z : g_{f_1}(z) > \rho/d^2 \}$.

Continue inductively to obtain a sequence $\{ \hat{h}_n \}$ of
$K$-quasiconformal maps such that 
$\hat{h}_n$ is a conformal conjugacy in 
$\{ z : g_{f_1}(z) > \rho/d^{n+1} \}$.
All of the maps $\hat{h}_n$ agree with $\hat{h}_0$ in a neighbourhood
of $\infty$. By passing to the limit of a subsequence
we obtain a  $K$-quasiconformal 
conjugacy $\hat{h}$ (see~\cite{lehto-73} ch. II.5).
which is conformal in the basin of infinity $\Omega(f_1)$ and
asymptotic to the identity at $\infty$.
Since $J(f_1)$
has measure zero, the conjugacy $\hat{h}: \C \rightarrow \C$ must be in fact
an affine translation. 
But $f_1$ and $f_2$ are monic
and centered  so we conclude that $f_1 = f_2$.
\hfill
$\Box$

\smallskip
To show that the set $S_{\cp}$ of polynomials
in $\vshift$ sharing a critical portrait $\cp$ is
a sub-manifold parameterized by the escape rates of 
the critical points we need the following result
of Branner and Hubbard (see~\cite{branner-88} ch. I.3):

\begin{lemma}
\label{compact-l}
Given $\rho >0$, let $\mathcal{B} \subset \param$ be the set
formed by polynomials $f$ such that:
$$ \max g_f(c) \leq \rho$$
where the maximum is taken over the critical
points $c$ of $f$. Then $\mathcal{B}$ is compact.
\end{lemma}

\begin{lemma}
The map $G$ is onto and 
the set $S_{\cp}$ is a $m$-dimensional real analytic sub-manifold.
\end{lemma}

\noindent
{\sc Proof:}
Given a critical portrait $\cp = \{ \cp_1 , \dots , \cp_m \}$
we want to show that $G(S_{\cp})$ is
$$W = \{  (r_1, \dots ,r_m): d^n \cdot \cp_i \in \cp_j
\Rightarrow d^n r_i > r_j \}.$$
Proposition~\ref{goldberg-p} says that $G(S_{\cp})$ is
not empty and Lemma~\ref{subset-l}
guarantees that $G(S_\cp) \subset W$. Note that $W$ is convex, 
in particular connected. So it is enough to
show that $G(S_{\cp})$ is both closed and open.

\smallskip
To show that $G(S_{\cp})$ is closed let
$f_n \in S_{\cp} \subset \vshift$  be such that
$$G(f_n) \rightarrow G_0 = (\rho_1 , \dots , \rho_m) \in W$$
The set of polynomials $f$
such that:
$$g_f(c) \leq \max \{ \rho_1, \dots, \rho_m \}$$ 
is compact (Lemma~\ref{compact-l}). Therefore, by passing to a subsequence,
we may assume that $f_n \rightarrow f$. 
Label the critical points
$c_1(f_n), \dots , c_m(f_n)$ of $f_n$ so that 
the external radii of $f_n$ with arguments in $\cp_i$
terminate at $c_i(f_n)$. The critical points
$c_i(f_n)$ converge to a critical point $c_i(f)$ of $f$.
Moreover, $c_1(f), \dots, c_m(f)$ is a list of all
the critical points of $f$. A priori we do not know
whether there are any repetitions in this list or not.
By continuity of
the Green function we know that $g_f(c_i(f)) = \rho_i > 0$.
Hence, $f$ lies in the shift 
locus $\shift$. We must show that $f$ actually lies 
in the visible shift locus $\vshift$ which automatically
implies that $f \in S_{\cp}$ (continuity of $\Pi$) and
$G(f) = G_0$ (continuity of $g$).

For $t \in \cp_i$, consider the broken flow lines 
$$\overline{R}^{*t}_{f_n}: [g_{f_n}(c_i(f_n)), \infty] \rightarrow \hat{\C}$$
which go from $c_i(f_n)$
to $\infty$. By passing to a subsequence,
 $\overline{R}^{*t}_{f_n}$ converge
to a broken flow line 
$$ \gamma_f : [\rho_i , \infty] \rightarrow \hat{\C}$$
connecting $c_i(f)$ to $\infty$.
This broken flow line $\gamma_f$, near infinity,
coincides with $R^{*t}_f$.

\smallskip
We claim that $\gamma_f$ is the
external radius $R^{*t}_f$ union $c_i(f)$.
In fact, consider the B\"ottcher maps
$\phi_{f_n} : \Omega^*(f) \rightarrow U_{f_n}$.
From Section~\ref{dyn-s},
it is not difficult to conclude that, given $\epsilon >0$, there 
exists a definite neighbourhood $V$ of
$[e^{\rho_i} + \epsilon , \infty) e^{2 \pi i t}$
contained in $U_{f_n}$ 
(i.e. $V$ is independent of $n$).
Thus, $(e^{\rho_i} , \infty) e^{2 \pi i t}$ is contained
in the domain of $\psi_f = \phi^{-1}_f: U_f \rightarrow \Omega^*(f)$.
Hence, for $t \in \cp_i$,
the external radius $R^{*t}_f$ terminates at $c_i(f)$.

The critical value $f(c_i(f))$
belongs to $\Omega^*(f)$ because
the same argument used above shows that $e^{d \rho_i + 2 \pi i dt} \in U_f$ 
and  $f(c_i(f)) = \psi_f (e^{d \rho_i + 2 \pi i dt} )$.

We need to show that $c_1(f), \dots ,c_m(f)$
are distinct. For this we apply a counting argument.
The local degree of $f_n$ at $c_i(f_n)$ is $d_i = |\cp_i|$.
If $c = c_{i_1}(f) = \dots = c_{i_k}(f)$ then
the local degree of $f$ at $c$ is $d_{i_1} + \cdots + d_{i_k} - k+1$.
But there are $d_{i_1} + \cdots + d_{i_k}$ external
radii terminating at $c$. Moreover, the critical value
$f(c)$ belongs to only one external radius. Thus $k=1$,
and $c_1(f), \dots , c_m(f)$ are distinct. 
Hence, $f \in \vshift$ and $G(S_{\cp})$ is closed.

\bigskip
We proceed to show that $S_\cp$ is a $m$-real dimensional
manifold at the same time that we show that $G(S_{\cp})$ is open in $W$.
 Consider
$f_0 \in S_{\cp}$ and observe that
$f_0$ belongs to a $m$-complex dimensional sub-manifold $M$ of
$\param$ formed by polynomials that have
$m$ distinct critical points $c_1(f), \dots, c_m(f)$ with
corresponding local degrees
$d_1 , \dots , d_m$. In $M$, the
critical points vary holomorphically with $f$.
The visible shift locus $\vshift$ contains an open neighbourhood $V$
of $f_0$ {\bf in} $M$. 
Consider the holomorphic map
$$\begin{array}{lccl}
     \Phi:  & V & \rightarrow & (\CDC)^m \\
           & f & \mapsto & (\phi_f \circ f(c_1(f)), \dots , 
                        \phi_f \circ f(c_m(f))\\
\end{array}$$
Since the external radii that terminate at 
$c_i(f)$ vary continuously with $f \in V$ we have that $\Phi(f)$
completely determines $\cp(f)$ (after shrinking
$V$ if necessary). Moreover,
$g_f(c_i(f))$ is also determined by
$\Phi(f)$. By Lemma~\ref{inj-l}, it follows that $\Phi$ is injective
and hence a biholomorphic isomorphism between $V$
and its image.
Thus, 
$$ (\phi_{f_0} \circ f_0 (c_1(f_0)), \dots , 
\phi_{f_0} \circ f_0 (c_m (f_0))$$
is a regular value.
\hfill
$\Box$

\section{Impressions}
\label{imp-s}
\noindent \indent
Here we prove that critical portrait
impressions are connected and
that their union is the portion of
$\partial \shift$ contained in $\conec$.
We start with the basic properties of $\ang$.

\begin{lemma}
\label{ang-l}
Given $r_0 >0$, let $\mathcal{ E} \subset \vshift$
be the set of polynomials $f$ such that all the critical 
points $c$ of $f$ have escape rate $g_f(c) = r_0$.
Then:
$$\Pi_{|\mathcal{ E}} : {\mathcal{ E}} \rightarrow \ang$$
is a homeomorphism. Furthermore, $\ang$ is compact and connected.
\end{lemma} 

\noindent
{\sc Proof:} 
First we show that $\ang$ is Hausdorff.
Consider two distinct critical portraits
$\cp$, $\cp^{\prime}$ and observe that 
the $\cp$-unlinked classes $L_1, \dots , L_d$ must be distinct from
the $\cp^{\prime}$-unlinked classes $L^{\prime}_1,
\dots, L^{\prime}_d$. Hence,
the  union of a $\cp$-unlinked class
and a $\cp^{\prime}$-unlinked class has total length 
strictly greater than $1/d$. Pick closed sets 
$$X_1 \subset L_1, \dots , X_d \subset L_d$$
$$X^{\prime}_1 \subset L^{\prime}_1, \dots , 
X^{\prime}_d \subset L^{\prime}_d$$
such that, for $1 \leq i,j \leq d$,
 the union $X_i \cup X^{\prime}_j$ has 
measure greater than $1/d$.
It follows that the 
neighbourhood $V = V_{X_1} \cap \dots \cap V_{X_d}$ of
$\cp$ and the neighbourhood 
$V^{\prime} = V^{\prime}_{X^{\prime}_1} \cap \dots 
\cap V^{\prime}_{X^{\prime}_d}$ of $\cp^{\prime}$ 
are disjoint.

\smallskip
The set $\mathcal{ E}$ is compact (Lemma~\ref{compact-l}). 
$\Pi|_{\mathcal{ E}}$ is one to one and onto from a
compact to a Hausdorff space. Thus, 
$\Pi|_{\mathcal{ E}}$ is a homeomorphism and $\ang$ is compact.

\smallskip
To show that $\ang$ is connected,  consider
the subset $S \subset \ang$ formed by critical portraits
of the form:
$$\cp = \{ \{ \theta, \theta + 1/d, \dots, \theta+ (d-1)/d \} \}.$$
These are the critical portraits corresponding to polynomials
of the form $z^d + a_0$. Notice that $S \cong \S$, in particular
$S$ is connected.
Pick a critical portrait which is a collection
$\cp = \{ \cp_1, \dots , \cp_m \}$ of $m \geq 2$ sets.
It is enough to show that $\cp$ lies in the same connected component 
of $\ang$ than a critical portrait formed by $m-1$ sets.
In fact, assume that the angular distance $\epsilon$ between
the pair $\cp_1$ and $\cp_2$ is minimal amongst
all possible pairs. Without loss of generality, there exists 
$\theta_1 \in \cp_1$ and $\theta_2 \in \cp_2$ such
that $\theta_2 = \theta_1 + \epsilon$.
Therefore, $\{\cp_1 + \epsilon s, \cp_2 , \dots , \cp_m \}$
where $ 0 \leq s < 1$ is a path between $\cp$
and $\{ (\cp_1 + \epsilon) \cup \cp_2, \dots, \cp_m \}$. 
\hfill
$\Box$

\smallskip
\begin{corollary}
The impression $I_{\conec}(\cp)$ of a critical portrait is 
a non-empty connected subset of $\partial \shift$.
Moreover, 
$$\bigcup_{\cp \in \ang} I_{\conec}(\cp) = \partial \shift \cap \conec$$
\end{corollary}

\noindent
{\sc Proof:} 
Take a basis of connected neighbourhoods $V_n \subset \ang$ around 
$\cp$ and let 
$$W_n = \{ f \in \vshift : 
\Theta(f) \in V_n \mbox{ and } \sum_{f^{\prime}(c) = 0} g_f(c) \leq 1/n \}.$$
From Theorem~\ref{coor-th}, it follows that $W_n$ is connected. 
Thus,
$$I_{\conec}(\cp) = \bigcap_{n \geq 1} \overline{W}_n$$
is connected. To see that impressions cover all 
of $\partial \shift \cap \conec$ just observe that
$\vshift$ is dense and $\ang$ is compact.
\hfill
$\Box$

\newpage
\part*{Chapter 3: Rational Laminations}
\vspace{-0.5cm}
\indent
\section{Introduction}
In this Chapter we study some topological features of
 the Julia set of polynomials with all cycles repelling.
Important objects that help us understand the 
 topology of a connected Julia set are prime end impressions
 and external rays. 

Let us briefly recall the definition of a prime end 
impression:

\begin{definition}
Consider a polynomial $f$ with connected Julia set $J(f)$ and
 denote the 
 inverse of the B\"ottcher map 
 by $\psi_f: \CDC \rightarrow \Omega(f)$. 
Given $t \in \S$, we say that $z \in J(f)$
 belongs to the {\bf prime end impression} $Imp(t)$ if
 there exists a sequence $\zeta_n \in \CDC$ 
 converging to $e^{2 \pi i t}$
 such that the points $\psi_f (\zeta_n)$ converge
 to $z$.
\end{definition}

Note
 that if the external ray $R^t_f$ lands at $z$ then
 $z$ belongs to the prime end impression $Imp(t)$. 
In particular, for $t \in \QS$
 the impression $Imp(t)$ contains a 
 pre-periodic or periodic point. 

A prime end impression $Imp(t)$ is a singleton if and only if 
 $\psi_f$ extends continuously to $e^{2 \pi i t}$. 
A result, due to Carath\'eodory, says that every impression
 is a singleton if and only if $J(f)$ is locally connected.
In this case, $\psi_f$ extends continuously to the boundary
 $\partial \D \cong \S$ and establishes a semiconjugacy
 between 
 $m_d: \S \rightarrow \S$ and the map $f: J(f) \rightarrow J(f)$. 
Recall that $m_d(t) = d \cdot t \mbox{ (mod $1$)}.$

The Julia set is not always locally connected.
But, under the assumption that all the cycles of $f$
 are repelling, we show that $J(f)$ is locally connected
 at every pre-periodic and periodic point (see 
 Theorem~\ref{par-imp-th} below). 
Moreover,
 $\psi_f$ extends continuously at every rational
 point $t$ in $\QS \subset \S \cong \partial \D$.  
That is, the topology of $J(f)$ is rather ``tame''
 at periodic and pre-periodic points and,
 the boundary behavior of $\psi_f$ is also ``tame''
 in the rational directions.

Another closely related issue 
 is to know how many impressions contain 
 a given point $z \in J(f)$. 
We apply the results from
 Chapter 1 and 2 to show that $z$ is contained in 
 at most finitely many impressions provided that 
 all  cycles of $f$  are repelling. 
Observe that while 
 there might be no external ray landing at $z$ there
 are always impressions which contain $z$.

\begin{theorem}[Impressions]
\label{par-imp-th}
Consider a monic polynomial $f$ with connected Julia set
$J(f)$ and all   cycles repelling.
Let $\mbox{Imp}(t) \subset J(f)$ be the prime end impression
corresponding to $t \in \S$ under the B\"ottcher map. 

(a) If $t \in \QS$ then $\mbox{Imp}(t) = \{ z \}$
where $z$ is a periodic or pre-periodic point.

(b) If $t \notin \QS$ then $\mbox{Imp}(t)$ does not 
contain periodic or pre-periodic points.

(c) If $z \in J(f)$ is a periodic or pre-periodic point
then $J(f)$ is locally connected at $z$.

(d) Every $z \in J(f)$  is contained in at least one and
at most finitely many impressions.
\end{theorem}

Loosely, a polynomial $f$ with all cycles
 repelling has ``a lot'' of periodic and 
 pre-periodic orbit portraits which are non-trivial.
We will make this more precise later.
Now let us observe that there is at least one fixed
 point $z$ with more than one ray landing at it.
This is so because there are $d$ repelling fixed points
 and only $d-1$ fixed rays.

Roughly speaking,
 the abundance of nontrivial periodic and pre-periodic 
 orbit portraits gives rise to a wealth of possible partitions of 
 the complex plane
 into ``Yoccoz puzzle pieces''.
The proof of parts (a), (b) and (c) of the previous Theorem
 relies on finding an appropriate puzzle piece for each
 periodic or pre-periodic point $z$. 
As mentioned above, the proof of part (d) uses  
 results from the previous Chapters.

\smallskip
Under the assumption that all cycles are repelling,
 every pre-periodic or periodic
 of $f$ is the landing point of rational rays 
 (see Theorem~\ref{land-th}).
The pattern in which rational external rays land is captured by
 the rational lamination $\ratrel(f)$ of $f$:

\begin{definition}
Consider a polynomial $f$ with connected Julia set $J(f)$.
The equivalence relation
 $\ratrel(f)$ in $\QS$ that identifies $t, t^\prime \in \QS$
 if the external rays $R^t_f$ and $R^{t^\prime}_f$ land
 at a common point is called the {\bf rational lamination} of $f$.
\end{definition}

\noindent
{\bf Remark:} We work with the definition of rational 
 lamination which appears in~\cite{mcmullen-94}. 
The word ``lamination'' corresponds to the usual
 representation of this equivalence relation in the unit 
 disk $\D$. 
That is, each equivalence class $A$ is represented
 as the convex hull (with respect to the Poincar\'e metric)
 of  $A \subset \S \cong \partial \D$.
The use of ``laminations'' to represent the 
 pattern in which external rays of a polynomial land or can land 
 was introduce by Thurston in~\cite{thurston-85}
 (also see~\cite{douady-93}).

\smallskip
We explore the basic properties that
 will allow us, in Chapter 4, to describe the equivalence
 relations in $\QS$ that arise as the rational 
 lamination of a polynomial with all cycles repelling.
Recall that we fix the standard orientation in $\S$ and
 use interval notation accordingly. 

\begin{proposition}
\label{lam-p}
Let $\ratrel(f)$ be the rational lamination of a polynomial
 $f$ with all   cycles repelling and connected Julia set $J(f)$.
Then:

\noindent
(R1) $\ratrel(f)$ is a closed equivalence relation in $\QS$

\noindent
(R2) Every $\ratrel(f)$-equivalence class $A$ 
is a finite set.

\noindent
(R3) If $A_1$ and $A_2$ are distinct equivalence classes
then $A_1$ and $A_2$ are unlinked.

\noindent
(R4) If $A$ is an equivalence class then $m_d(A)$ is an 
equivalence class.

\noindent
(R5) If $(t_1, t_2)$ is a connected component of 
$\S \setminus A$ where $A$ is an equivalence class
then $(dt_1, dt_2)$ is a connected component of $\S \setminus m_d(A)$.

Moreover, $\ratrel(f)$ is maximal with respect to properties
(R2) and (R3). 
Furthermore, there exists a unique closed equivalence
 relation $\lambda_{\R}$ in $\S$ which agrees with $\ratrel(f)$ in
 $\QS$ such that $\lambda_{\R}$-classes are unlinked.
Also, the equivalence classes of $\lambda_{\R}$ satisfy properties
 (R2) through (R5) above.
\end{proposition}

Recall that an equivalence relation $\ratrel$ 
(resp. $\lambda_\R$) is closed
if it is a closed subset of $\QS \times \QS$
(resp. $\S \times \S$). 
Also, by a ``maximal'' equivalence relation in $\QS$
we mean an equivalence relation which is maximal with
respect to the partial order determined by inclusion
in subsets of $\QS \times \QS$.

\smallskip
Each $\ratrel(f)$-equivalence class $A$ is the 
 type of a periodic or pre-periodic point. 
We refer to  $A$ as a {\bf rational
 type} to emphasize that $A$ is the  type of a point.
In particular,  $\ratrel(f)$-equivalence classes 
 inherit the basic properties, (R2) through (R5) above,
 of types discussed in Chapter~1.
The property (R1), i.e. $\ratrel(f)$ is closed, is more
 delicate.  Again, our proof of (R1) will rely on constructing a 
 puzzle piece around each periodic or pre-periodic point.

\smallskip
It is worth pointing out that,
 from the above Proposition, it follows that 
 when the Julia set $J(f)$ is locally connected 
 $J(f)$ must be homeomorphic
 to $(\S)/ \lambda_{\R}$.
Moreover, $m_d$ projects to a map from
 $(\S)/ \lambda_{\R}$ onto itself.
Thus, $\lambda_{\R}$ gives rise to an
 ``ideal''  model for the topological dynamics of $f$.

It is also worth mentioning that, for quadratic polynomials
with all cycles repelling, the Mandelbrot local
connectivity Conjecture implies that
the rational lamination uniquely 
determines the quadratic polynomial
in the family $z \mapsto z^2 + a_0$.
We do not expect this to be true 
for cubic and higher degree polynomials.

\medskip
This Chapter is organized as follows:

\smallskip
In Section~\ref{puz-s} we fix notation and summarize
 basic facts about Yoccoz puzzle pieces. 
Our notational approach is  slightly nonstandard 
 because we need some  flexibility to work, at the
 same time, with all the possible
 puzzles for a given polynomial.

\smallskip
Section~\ref{pai-s} contains the proofs of the 
results discussed above.

\section{Yoccoz Puzzle}
\label{puz-s}

\noindent \indent
For this section, unless otherwise stated, we let $f$ be a monic
 polynomial with connected Julia set $J(f)$ and, possibly,
 with non-repelling cycles.

\bigskip
Every point $z \in J(f)$ which eventually maps onto a repelling
 or parabolic periodic point is the landing point of
 a finite number of rational rays. 
The arguments of 
 these external rays form the type $A(z)$ of $z$.
For short, we simply say that $A(z)$ is a
 {\bf rational type} for $f$.
Note that the rational types for $f$ 
 are the equivalence classes
 of the rational lamination $\ratrel(f)$. 
We will usually prefer to call them ``rational types'' 
 to emphasize that we are talking about rational rays
 landing at a common point rather than an abstract equivalence
 class of $\ratrel(f)$.

\bigskip
Rational external rays landing
 at a finite collection of points chop the complex plane 
 into puzzle pieces:

\newpage
\begin{definition}
Let $G=  \{ A(z_1) , \dots , A(z_p) \}$ be a collection
of rational types.
The union $\Gamma$ of the external rays with arguments in 
$A(z_1) \cup \cdots \cup A(z_p)$ together with their
landing points $\{z_1, \dots, z_p \}$ cuts the complex plane
into one or more connected components.
A connected component $U$ of $\C \setminus \Gamma$ is called an
{\bf unbounded  $G$-puzzle piece}. 
The portion of an unbounded $G$-puzzle piece
contained inside an equipotential is called
a {\bf bounded  $G$-puzzle piece}.
\end{definition}

When irrelevant or clear from the context we will 
not specify whether a given puzzle piece is bounded
or unbounded. 

\smallskip
Usually, it is convenient to start with a forward invariant
 puzzle.
That is, a puzzle $G^{\prime} = \{ A(z_1), \dots, A(z_p) \}$
 such that 
 $$\{ A(f(z_1)), \dots , A(f(z_p)) \} \subset G^{\prime}.$$
Then we consider the collection
 $G$  formed by all the rational types that
 map onto one in $G^\prime$ (i.e. $G$ is the
 pre-image of $G^\prime$).
In this case, 
 if $U$ is an unbounded $G$-puzzle piece then $f$ maps 
 $U$ onto a $G^\prime$-puzzle piece $U^\prime$.
Moreover, if $k$ is the number of critical points
 in $U$ counted with multiplicities then
 $f : U \rightarrow U^{\prime}$ is
 a degree $k+1$ proper holomorphic map.
Furthermore, $U \subset U^\prime$ or $U$ and $U^\prime$ are
 disjoint.
A similar situation occurs when $U$ is bounded by the
 equipotential $g_f = \rho$ and $U^{\prime}$ is bounded by the
 equipotential $g_f = d \rho$.

\smallskip
Puzzle pieces are useful to construct a basis of
connected neighborhoods around a point in the
Julia set $J(f)$ because of the following:

\begin{lemma}
\label{lc-l}
Let $G$ be a collection of rational types and
 $U$ be a $G$-puzzle piece. 
Then $\overline{U} \cap J(f)$ is connected.
\end{lemma}

\noindent
{\sc Proof:} 
We proceed by induction on
 the number of rational types in $G$ (see~\cite{hubbard-93}).

Let $G_n = \{ A(z_1), \dots , A(z_n) \}$. 
The statement is true, by hypothesis, when $n=0$ (taking
the associated puzzle piece to be the entire plane).
For $n \geq 1$, denote by $V$ the  $G_{n-1}$-puzzle piece
which contains $z_{n}$. 
We assume that $X = \overline{V} \cap J(f)$ is
connected and show that the same holds for
$G_{n}$-puzzle pieces.
 
Denote by $S_1 , \dots , S_q$ the sectors based at 
$z_{n}$. The $G_{n}$-puzzle pieces that are not 
$G_{n-1}$-puzzle pieces are 
$U_1 = S_1 \cap V, \dots , U_q = S_q \cap V$.

For each $i = 1, \dots , q$,
we must show that 
$$X_i =  \overline{U_i} \cap J(f) =  \overline{U_i} \cap X = 
\overline{S_i} \cap X$$
is connected. 
Without loss of generality we show that $X_1$ is connected.
 
Let $W_1$ and $W_2$ be two disjoint open sets 
such that:
$$ X_1 = (W_1 \cap X_1) \cup (W_2 \cap X_1). $$
Since $z_{n} \in X_1$ we may assume that $z_{n} \in W_1$.
It follows that 
$$(W_1 \cap X_1) \cup (S_2 \cap X)\cup \dots \cup(S_q \cap X)$$
and $W_2 \cap X_1$ are two disjoint open sets (in $X$) whose union is 
the connected set $X$. Hence, $W_2 \cap X_1$ is empty and 
$X_1$ is connected. 
\hfill $\Box$

\smallskip
Given a puzzle piece $U$,
we keep track of the situation in the circle at infinity
by considering:
$$\pi_{\infty} U = \{ t \in \S : R^t_f \cap U \neq \emptyset \}.$$
Notice that if $t \in \pi_{\infty} U$ then the impression $Imp(t)$
is contained in $\overline{U} \cap J(f)$. Moreover, if
$t \notin \overline{\pi_{\infty}U}$ then $Imp(t)$ is
contained in $(\C \setminus U) \cap J(f)$.

\smallskip
Douady's Lemma, below, will enable us to show that 
 certain subsets of the circle
 are finite (see~\cite{milnor-90}):

\begin{lemma}
If $E \subset \S$ is closed and $m_d$ maps $E$ homeomorphically
onto itself then $E$ is finite.
\end{lemma}

\smallskip
Puzzles will allow us to extract polynomial-like maps
from $f$.
Following Douady and Hubbard, we say that  $g : V \rightarrow V^{\prime}$ is 
{\bf polynomial-like map} if $V$ and $V^{\prime}$ are Jordan
domains in $\C$ with smooth boundary such that
$V$ is compactly contained in $V^{\prime}$  (i.e. $\overline{V} \subset V$) and 
$g$ is a degree $k >1$ proper holomorphic map.

A polynomial-like map $g$ has a {\bf filled Julia set} $K(g)$ and 
a {\bf Julia set} $J(g)$ just as polynomials do: 
$$K(g) = \bigcap_{n \geq 1} f^{-n}(\overline{V})$$
$$J(g) = \partial K(g)$$
Moreover, a polynomial-like map can be extended to a map from  $\C$ 
onto itself which is quasiconformally conjugate to
a polynomial: 

\begin{theorem}[Straightening]
If $g: V \rightarrow V^{\prime}$ is a degree $k$ polynomial-like map then there exists a quasiconformal map
$h : \C \rightarrow \C$ and a degree $k$ polynomial $f$
such that $h \circ g = f \circ h$ on a neighbourhood of $K(g)$.
\end{theorem}

In particular, the quasiconformal map $h$ of the Straightening
Theorem takes the, possibly disconnected,
 Julia set  of $g$ onto the
 Julia set  of $f$.

\smallskip
Under certain conditions we can apply the
{\it thickening procedure} to extract a polynomial-like map from a polynomial and a puzzle:

\begin{lemma}[Thickening]
\label{thick-l}
Consider a collection of repelling periodic orbits $\orb_1 , \dots , \orb_k$ 
 of $f$ and, for some $l > 0$, let:
$$Z = \bigcup^{l}_{i= 0} f^{-i}(\orb_1 \cup \dots \orb_k).$$ 
Assume that $Z$ does not contain critical points.

Consider 
$G= \{ A(z) : z \in Z \}$ and $G^\prime = \{ A(f(z)) : z \in Z \}$.
Suppose that $U$ (resp. $U^{\prime}$) is a bounded $G$-puzzle piece
(resp. $G^{\prime}$-puzzle piece) such that
$U \subset U^\prime$  and $f : U \rightarrow U^\prime$ is a degree
$k>1$ proper map.

Then there exists Jordan domains $V$ and $V^\prime$ 
with smooth boundary such that
$U \subset V$, $\overline{V} \subset V^\prime$ and
$$f : V \rightarrow V^{\prime}$$ is a degree $k$
proper map (i.e. a polynomial-like map).
\end{lemma}

\noindent
{\sc Proof:} We restrict to the case in which
there is only one periodic orbit 
$\orb = \{ z_0, \dots, z_{p-1} \}$ involved.
The construction generalizes easily.
In order to fix notation let $g_f = \rho$ be the equipotential
inside which $U$ lies.

For each point $z \in Z$, let $\Gamma_{z}$ be the graph
formed by the union of the external rays landing
at $z$ and the point $z$. Since $f$ maps $\Gamma_{z}$
onto $\Gamma_{f(z)}$ injectively we can choose 
neighborhoods $W_z$ of $\Gamma_z$ such that:

$\bullet$ $W_z \cap W_w = \emptyset$ for $z \neq w$,

$\bullet$ $W_{f(z)} \subset f(W_z)$ and,

$\bullet$ $f_{|W_z}$ is injective.

\smallskip
Inside $W_z$ we will {\it thicken} the graph $\Gamma_z$.
First we construct open disks $D^{(l)}_z \subset W_z$ around
$z \in Z$ and $D^{(l-1)}_z \subset W_z$ around $z \in f(Z)$.
These disks will have several properties:

(a) For $z \in f(Z)$, $\overline{D}^{(l-1)}_z \subset D^{(l)}_z$.

(b) For  $z \in Z$,  $f(D^{(l)}_z) = D^{(l-1)}_{f(z)}$.

(c) For $z \in Z$, the portion of the external rays
landing at $z$ contained in $D^{(l)}_z$ is equal
to the portion of the external rays contained inside
the equipotential $g_f = r_0/d^l$ where $r_0 < \rho $ is independent
of $z$ and the external ray.

(d) $\partial D^{(l)}_z$ is smooth. 

(e) If $R^t_f$ lands at $z \in Z$
then there exists a small open arc
in $\partial D^{(l)}_z$ around 
$R^t_{f} \cap \partial D^{(l)}_z$ which is
contained in the equipotential $g_f = r_0/d^l$.

\smallskip
To construct these disks we start  by finding
$p+1$ nested disks around the periodic point $z_0$.
The construction only relies on the fact the 
$z_0$ is a repelling periodic point.

Pick an open
topological disk $D^{(0)}_{z_0} \subset W_{z_0}$ around $z_0$ such that:

(i) $\overline{D}^{(0)}_{z_0} \subset D^{(-p)}_{z_0}
= f^{\circ p} (D^{(0)}_{z_0})$.

(ii) $ \overline{D}^{(-p)}_{z_0} \subset W_{z_0}$.

(iii) $\partial D^{(0)}_{z_0}$ is smooth.

(iv) The portion of each external ray 
landing at $z_0$ inside $D^{(0)}_{z_0}$ is connected.

(v) For each point in $\partial D^{(0)}_{z_0} \cap \Gamma_{z_0}$ 
there exist a small open arc of $\partial D^{(0)}_{z_0}$ that
is contained the equipotential $ g_f = r_0 $.

Now choose a nested collection of $p-2$ disks between 
$D^{(0)}_{z_0}$ and $D^{(-p)}_{z_0}$:
$$D^{(0)}_{z_0} \subset
D^{(-1)}_{z_0} \subset  \cdots  \subset D^{(-p + 1)}_{z_0}
\subset D^{(-p)}_{z_0}$$
such that the closure of each disk is contained in the next disk
and for $n= -1, \dots , -p+1$:

(vi) $\partial D^{(n)}_{z_0}$ is smooth.

(vii) The portion of each external ray landing
landing at $z_0$ inside $D^{(n)}_{z_0}$ is connected.

(viii) For each point in $\partial D^{(n)}_{z_0} \cap \Gamma_{z_0}$ 
there exist a small open arc of $\partial D^{(0)}_{z_0}$ that
is contained the equipotential $g_f = r_0 /d^n $.

\smallskip
To build a disk $D^{(0)}_{z_i}$ around each periodic
point $z_i$ observe that $f^{\circ p-i}(z_i) = z_0$.
Let $D^{(0)}_{z_i}$ be the connected component 
of $f^{-(p-i)}(D^{(-(p-i))}_{z_0})$ which contains
$z_i$. It follows that $D^{(0)}_{z_i}$ has
the properties (iii), (iv), (v) above. Moreover,
$D^{(0)}_{z_{i+1}}$ is compactly contained in $f(D^{(0)}_{z_i})$
(subscripts mod $p$).

\smallskip
For $n \leq l-1$, define inductively disks around each point in $z \in Z$
as follows. If $z \in f^{-1}(w)$ then let $D^{(n+1)}_z$
be the connected component of $f^{-1}(D^{(n)}_w)$ that contains $z$.
Hence, $D^{(l)}_z$ is  defined for every
point $z \in Z$ and $D^{(l-1)}_z$ is defined 
for all $z \in f(Z)$. By construction, these disks have the
desired properties (a) through (e).

\smallskip
The second step is to {\it thicken} the rays landing at 
points in $Z$.
Choose  $\delta >0$ small enough so  that
 the following conditions hold:

(a) $$T^{(l)}_{t} = \{ z : g_f(z) \geq r_0/d^{l} \mbox{  and } 
\arg \phi_f(z) \in (t-\delta , t + \delta) \} \subset W_z.$$
where $z$ is the landing point of $R^t_f$.

(b) $T^{(l-1)}_{t} = f(T^{(l)}_{t}) \subset W_{f(z)}$.

(c) The portion of $T^{(l)}_t$ contained
in the equipotential $g_f = r_0/d^l$ lies in $\partial D^{(l)}_z$.

\smallskip
Finally, {\it thicken} the puzzle pieces $U$ 
contained inside the 
equipotential $g_f = \rho$:

$$V = (U \cup \bigcup_{z \in \partial U} D^{(l)}_z \cup 
                \bigcup_{R^t_f \cap \partial U \neq \emptyset} 
                        T^{(l)}_t) \cap \{ z : g_f (z) < \rho \}.$$
It follows that 
$$f(V) = V^{\prime} = 
          (U^{\prime} \cup \bigcup_{z \in \partial U^{\prime}} D^{(l-1)}_z 
                          \cup 
                \bigcup_{R^t_f \cap \partial U^{\prime} \neq \emptyset} 
                        T^{(l-1)}_t) \cap \{ z : g_f (z) < d \rho \}$$
is compactly contained in $V^{\prime}$.
After rounding of corners of $\partial V$ we obtained the
desired polynomial-like map.
\hfill $\Box$

\section{Puzzles and Impressions}
\label{pai-s}

\begin{figure}[htb]
\centerline{
\psfig{file=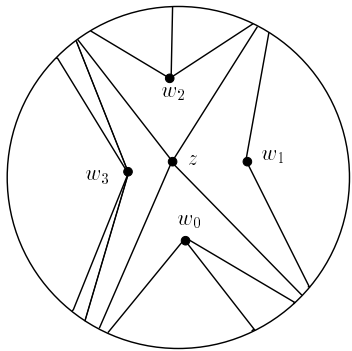,width=136pt}
}
\caption{A diagram of the situation 
described in Lemma~\ref{approx-l}}
\label{approx-f}
\end{figure}

\noindent \indent
In this section we prove Theorem~\ref{par-imp-th}. The proof 
of parts (a), (b), (c) relies on constructing a puzzle piece 
around each pre-periodic or periodic point $z$ of 
a polynomial $f$ with all   cycles repelling.
Simultaneously, in the Lemma below we show that the type $A(z)$ 
of $z$ is well {\it approximated} by other rational types
(see Figure~\ref{approx-f}). 
This result will be useful to prove Proposition~\ref{lam-p} and
in the next Chapter.

\begin{lemma}
\label{approx-l}
Consider a polynomial $f$ with all   cycles repelling
and connected Julia set $J(f)$. Let 
$A(z) = \{ t_0 , \dots , t_{p-1} \} \subset \QS$
be a rational type
(subscripts respecting cyclic order and mod $p$).
Given $\epsilon > 0$, 
there exists rational types
 $A(w_0), \dots , A(w_{p-1})$ such that $A(w_i)$ has elements
both in $(t_i, t_i + \epsilon)$ and $(t_{i+1} - \epsilon, t_{i+1})$.
Moreover, $w_0, \dots, w_{p-1}$ can be chosen so that
they do not belong to the grand orbit of a critical
point.
\end{lemma}

Recall that the grand orbit of $z$ is the set formed by the points 
$z^\prime$ such that $f^{\circ n}(z) = f^{\circ m}(z^\prime)$
for some $n, m \geq 0$.

\smallskip
\noindent
{\sc Proof of Lemma~\ref{approx-l}  and Theorem~\ref{par-imp-th}}
 (a),(b),(c):
First, consider a periodic point $z \in J(f)$. We
pass to an iterate of $f$ such that $z$ is a fixed point
and every ray landing at $z$ is also fixed. 
If necessary, we pass to an even higher iterate
of $f$ so that each periodic
point $\zeta$ in the post-critical orbit  is fixed
and the rays landing at $\zeta$ are also fixed.

Consider the collection $\mathcal{G}$ formed by the rational types 
$A(w)$ such that:

(i) $w$ is not in the grand orbit of a critical point.

(ii) $w$ is not in the grand orbit of $z$.

(iii) $A(w)$ has more than one element.

\smallskip
We saturate $\mathcal{G}$ by an increasing sequence 
 of finite collections $G_l$ of rational types.
That is, for $l \geq 1$, let $G_l \subset \mathcal{G}$ be the 
 collection formed by the 
 types of the points $w$ such that $f^{\circ l}(w) $ 
 is periodic of period less or equal to $l$.
Notice that $G^{\prime}_l = \{ A(f(w)) : A(w) \in G_l \}$
 is contained in $G_l$ and that $G_l$ is the preimage 
 of $G^{\prime}_l$.

Let $U_l(z)$ be the bounded $G_l$-puzzle piece which contains
the fixed point $z$ (inside the equipotential $g_f = \rho >0$). 
Notice that $U_l(z)$ maps onto the $G^{\prime}_l$-puzzle piece
$U^{\prime}_l(z)$, inside $g_f = d \rho$, which contains $z$.
Also, observe that $\{ U_l(z) \}$ is a decreasing sequence 
of puzzle pieces.

\smallskip
Since every fixed point type
with nonzero rotation number  belongs to
 $G_1$, it is not difficult to show that 
every fixed point contained in $\overline{U}_2(z)$
is the landing point of fixed external rays.

\smallskip
\noindent
{\it Claim 1:}
For $l$ large enough, $f : U_l(z) \rightarrow U^{\prime}_l(z)$ is
one to one.

\smallskip
\noindent
{\it Proof of Claim 1:}
The degree of $f$ restricted to $U_l(z)$ is non-increasing as $l$ 
increases. 
Thus, we proceed by contradiction and suppose that $f$
 restricted to $U_l(z)$ has degree $k>1$, 
 for all $l \geq l_0 \geq 2$.
Recall that each fixed point in 
 $\overline{U}_{l_0}(z) \subset \overline{U}_2 (z)$ is 
 the landing point of fixed rays. 
Observe that $f$ has $k-1$ critical points
(counting multiplicities) in $U_{l_0}(z)$.
The forward orbit of these critical points 
must be contained in $U_{l_0}(z)$, otherwise 
for some $l > l_0$ the degree of $f$ restricted to the puzzle piece
$U_l(z)$ would be less than $k$.
Apply the thickening procedure to 
 $f : U_{l_0}(z) \rightarrow U^{\prime}_{l_0}(z)$ to extract 
 a polynomial-like map $g$ of degree $k$ with 
 connected Julia set $J(g)$ and all its fixed points repelling.
After
 straightening we obtain a degree $k$ polynomial $P$ such that
 each of its  $k$ fixed points is repelling.
Each fixed point of $g$ is accessible through 
 a fixed ray in the complement of 
 the filled Julia set $K(g) = J(g) \subset J(f)$.
After straightening, each fixed point 
 of $P$ is accessible through a fixed arc in the complement
 of $K(P) = J(P)$. 
Hence, every one of the $k$ fixed points of $P$ is the 
 landing point of a fixed
 ray $R^t_P$. Since there are only $k-1$ fixed rays of $P$,
 this is impossible.

\smallskip
In the circle at infinity,
recall that  
$$\pi_{\infty} U_l(z) = 
\{ t \in \S : R^t_f \cap U_l(z) \neq \emptyset \},$$
and let $E = \cap_{l \geq 1} \overline{\pi_{\infty} U_l(z)}$.

By the previous Claim, for $l$ large enough, 
$f$ maps $\overline{U}_l(z)$ homeomorphically
onto its image $\overline{U}^{\prime}_l(z)$. 
It follows that 
$m_d$ is a cyclic order preserving
bijection from $\overline{\pi_{\infty} U_l (z)}$ onto
its image $\overline{\pi_{\infty} U^\prime_l (z)}$.
Since, $$U_l(z) \subset U^\prime_l(z) \subset U_{l-1}(z) $$
we have that 
$$\pi_\infty U_l(z) \subset \pi_\infty U^\prime_l(z) 
        \subset \pi_\infty U_{l-1}(z) $$
and we conclude that $m_d$ leaves $E$ invariant.
 By Douady's Lemma, $E$ is finite.

\smallskip
Recall that the rays landing at $z$ are fixed.
Thus, $E$ contains a fixed point
 of $m_d$.  Since ${m_d}_{|E}$ is cyclic order preserving
 it follows that $E = \{ s_0 , \dots, s_{m-1} \}$ 
 is a collection of 
 arguments fixed under $m_d$
 (subscripts respecting cyclic order and mod $m$).

\smallskip
We claim that $A(z) = E$. In fact, take $\epsilon$ small enough
so that an $\epsilon$ neighborhood of $E$ in $\S$
is mapped by $m_d$ injectively onto its image.
By construction of $E$ there exists types 
$A(w_0), \dots, A(w_{m-1})$
in $\mathcal{G}$ such that
$A(w_i)$ has elements both in $(s_i , s_i + \epsilon)$ and
$(s_{i+1} - \epsilon, s_{i+1})$. 

Now consider a bounded
$\{ A(w_0), \dots, A(w_{m-1}) \}$-puzzle piece $U$ which contains
$z$. It follows that
$f$ maps $U$ onto a  domain $U^{\prime}$ which 
compactly contains $U$. Moreover, $f$ is univalent in $U$.
The inverse branch of $f$ that takes $U^{\prime}$ onto $U$
is a strict contraction in $\overline{U}$ with respect to 
the hyperbolic metric on $U^{\prime}$. 
Hence:
$$ \cap_{k \geq 1} {f_{|U}}^{-k}(\overline{U}) = \{ z \}.$$
The external rays with
arguments in $E$ must land at $z$ and  $A(z) = E$.
Therefore, the Lemma is proved for periodic
$z$. By Lemma~\ref{lc-l},  
${f_{|U}}^{-k}(\overline{U}) \cap J(f)$ is connected.
Thus, $J(f)$ is locally connected 
at $z$. Also, for each
$t \in A(z)$ we have that $Imp(t) = \{z \}$.
Moreover, if $t \notin A(z)$ then $Imp(t)$ cannot contain $z$.

\smallskip
A similar situation occurs
 at a pre-periodic point $\tilde{z} \in f^{-k}(z)$.
In fact, for $l \geq 1$, let
 $\tilde{G}_l$ be the  collection of rational types
that map under $m^{\circ k}_d$ onto a rational type in $G_l$.
Consider the
 $\tilde{G}_l$-puzzle piece  $\tilde{U}_l(\tilde{z})$,
 inside $g_f = \rho/d^k$,
 which contains $\tilde{z}$ and we note that 
$$f^{\circ k} ( \overline{\tilde{U}_l(\tilde{z})}) = 
\overline{U_l(z)}.$$ 
The connected set
$$X = \cap_{l \geq 1} \overline{\tilde{U}_l(\tilde{z})} \cap J(f)$$
is mapped by $f^{\circ k}$ onto $\{ z \}$. 
Thus, $X = \{ \tilde{z} \}$. The Lemma and parts (a), (b), (c) of 
Theorem~\ref{par-imp-th} follow.
\hfill
$\Box$

\bigskip
In order to prove part (d) of Theorem~\ref{par-imp-th}
 we consider the intersection $X(z)$ of all
 the unbounded puzzle pieces which contain a point $z$ with
 infinite forward orbit. 
Then we show that $X(z)$
 contains only finitely many external 
 rays. 

\medskip
This situation can be worded in terms of $\ratrel(f)$-unlinked
classes. Recall  that the $\ratrel(f)$-equivalence classes
are the rational types for $f$:

\begin{definition}
We say that $t, t^\prime \in \S \setminus \QS$ are 
$\ratrel(f)$-unlinked equivalent 
if for all rational types $A(z)$, 
$\{ t , t^\prime \}$ and $A(z)$ are unlinked.
\end{definition}

At the same time, we  obtain a result needed in the next 
Chapter and the proof of Proposition~\ref{lam-p}.

\begin{lemma}
\label{e-l}
Every $\ratrel(f)$-unlinked class $E$ is a finite set. 
Moreover, given $\epsilon >0$, if $E = \{ t_0 , \dots , t_{p-1} \}$
(subscripts respecting cyclic order and mod $p$) then
there exists rational types
 $A(w_0),\brkOK \dots ,\brkOK A(w_{p-1})$ such that $A(w_i)$ has elements
both in $(t_i, t_i + \epsilon)$ and $(t_{i+1} - \epsilon, t_{i+1})$.
\end{lemma}

Roughly, the idea is to realize some iterate $m_d^{\circ k}(E)$
 of $E$ as the type of 
 some Julia set element $\zeta$ of some polynomial $g$.
The point $\zeta$ will have infinite forward orbit.
This allow us to use the fact,
proved in Chapter~1, that the 
type of points with infinite forward orbit have finite
cardinality.
The polynomial $g$ will belong to the visible shift locus
where the pattern in which external rays land is 
completely prescribed by the critical portrait $\cp$ of
$g$ (see Lemma~\ref{lan-vshift-l}).
That is, we look for a critical portrait $\cp$
such that $m^{\circ n}_d (E)$ is contained in 
a $\cp$-unlinked class, for $n$ sufficiently large.
Hence, some iterate of $E$ will be contained in 
the type $A(\zeta)$
of a Julia set element $\zeta$ for
a polynomial $g$ in the visible shift
locus with critical portrait $\cp$.
Then, we can apply Theorem~\ref{con-th} to conclude
that $E$ is finite:

\bigskip
\noindent
{\sc Proof of Lemma~\ref{e-l} and Theorem~\ref{par-imp-th}~(d):} 
We saturate the 
rational types of $f$ by 
a sequence $G_l$ of collections 
of rational types.
Namely, for $l \geq 1$, let $G_l$ be the collection formed
by the rational types  $A(w)$ where $f^{\circ l}(w)$
is periodic of period less or equal than $l$.

For each $z \in J(f)$ with infinite forward orbit, let
$U_l(z)$ be the unbounded $G_l$-puzzle piece which
contains $z$. Let
$$X(z) = \cap_{l \geq 1} \overline{U}_l(z)$$
and 
$$E(z) = \cap_{l \geq 1} \overline{\pi_{\infty} U_l (z)}.$$
Observe that if $z \in Imp(t)$ then the external
ray $R^t_f$ must be contained in $X(z)$. Thus,
to prove part (d) of the Theorem it is enough 
to show that $X(z)$ contains finitely many rays
or equivalently that $E(z)$ is finite.
Also, notice that if $\{ t,t^{\prime} \} \subset \S \setminus \QS$ is unlinked
with all the rational types of $f$ then  
the external rays $R^t_f$ and $R^{t^{\prime}}_f$ must 
be contained in the same $G_l$-puzzle piece, for all $l \geq 1$.
Thus, to prove the Lemma, it is also enough to show that $E(z)$
is finite for all $z$ with infinite forward orbit.

Now $X(z)$ cannot contain a periodic or pre-periodic point $w$
because the proof of the previous Lemma shows that there
is a sequence of puzzle pieces whose intersection
with $J(f)$ shrinks to $w$. 
This implies that $E(z) \subset \S \setminus \QS$.
Also, for two points $z$
and $z^{\prime}$,  $X(z) = X(z^{\prime})$ or, 
$X(z)$ and $X(z^{\prime})$ are disjoint.
Thus, $E(z) = E(z^{\prime})$ or, $E(z)$ and
$E(z^{\prime})$ are unlinked.

Notice that $f(X(z)) = X(f(z))$. We claim that,
for all $n,k \geq 1$, 
$$X(f^{\circ n}(z)) \mbox{ and } X(f^{\circ n+k}(z))$$ are disjoint. 
Otherwise,
for $l$ large enough, 
$f^{\circ k}(\overline{U_l(f^{\circ n}(z))}) \supset 
\overline{U_l(f^{\circ n}(z))}$ and $f$ would have a periodic
point in $X(f^{\circ n}(z))$. Situation that we already
ruled out.
Hence, $\{ m^{\circ n}_d(E(z)) \}^{\infty}_{n=0}$ are disjoint and 
pairwise unlinked.

\smallskip
To capture the location of the critical
points by means of a critical portrait,
consider the rational types 
$$A(c_1), \dots, A(c_j)$$
of the critical points which are pre-periodic
(if any). 
Let $$\cp_1 \subset A(c_1), \dots, \cp_j \subset A(c_j)$$
be such that the cardinality of $\cp_i$ agrees with the local
degree of $f$ at $c_i$ and $m_d(\cp_i)$ is a single argument.

For the critical points $c$ with infinite forward orbit,
list without repetition the sets $X(c)$:
$$X_{j+1}, \dots, X_{m} \subset \C$$
Denote by $E_i \subset \S$ the arguments of the external
rays contained in $X_i$. 
Let $k_i$ be the number of critical points (counted 
with multiplicities) which belong to $X_i$.
Pick a subset $\cp_i \subset E_i$ 
with $k_i + 1$ elements such that $m_d(\cp_i)$ is 
a single argument.

By construction, $\cp_1, \dots , \cp_m$ are pairwise
unlinked. Counting multiplicities, conclude
that $\cp = \{ \cp_1, \dots, \cp_m \}$ is a critical
portrait.

\smallskip
Consider a polynomial $g$ in the visible shift
locus $\vshift$
with critical portrait $\cp(g) = \cp$ (Chapter 2).
For $n$ large enough, $E(f^{\circ n}(z))$ is contained in 
a $\cp$-unlinked component. Thus, $E(f^{\circ n}(z))$ is
contained in the type $A(\zeta)$
of a point $ \zeta \in J(g)$ with infinite
forward orbit (Lemma~\ref{lan-vshift-l}).
 By Theorem~\ref{con-th}, $E(f^{\circ n}(z))$ is finite.
It follows that $E(z)$ is also finite.
\hfill
$\Box$

\medskip
We leave record, for later reference, 
of the critical portrait $\cp$ that we found
in the proof above. 
This critical portrait $\cp$ abstractly
captures the location of the critical points
in the Julia set $J(f)$ of a polynomial $f$ with all
cycles repelling. Observe that, a posteriori, there
are only finitely many choices of critical portraits in the 
construction of the previous proof:

\begin{corollary}
\label{exist-c}
Let $f$ be polynomial with all  cycles repelling 
and connected Julia set $J(f)$. Then there exists at least one
and at most finitely many  
critical portraits $\cp = \{ \cp_1 , \dots , \cp_m \}$
such that $\cp_i$ is either contained in a rational
type or $\cp_i$ is unlinked with every rational type.
\end{corollary}

We finish this Chapter by proving the basic properties
of the rational lamination of polynomials with all 
cycles repelling:

\smallskip
\noindent
{\sc Proof of Proposition~\ref{lam-p}:}
Properties (R2) to (R5) follow from the 
Lemmas~\ref{inv-por},~\ref{mon-por} and ~\ref{unl-por}. 
To show that $\ratrel(f)$ is closed, let 
$t_n,t^\prime_n \in \QS$ be $\ratrel(f)$-equivalent
and suppose that $t_n \rightarrow t \in \QS$ 
and $t^\prime_n \rightarrow t^\prime$.
Consider the $\ratrel(f)$-class $A$ of $t$ and observe that, for $n$
sufficiently large,  
Lemma~\ref{approx-l} implies that $t^{\prime}_n$ is trapped in an
$\epsilon$-neighborhood of $A$. Thus, $t^{\prime} \in A$ and
$\ratrel(f)$ is closed.

To prove the existence of $\lambda_{\R}$ just consider the equivalence
relation in $\S$ such that each equivalence class $B$ is either
a $\ratrel(f)$-class or a $\ratrel(f)$-unlinked class.
The same argument used above to prove that
$\ratrel(f)$ is closed shows that $\lambda_{\R}$ is
closed.

To prove the uniqueness of $\lambda_{\R}$ observe that 
Lemma~\ref{e-l} leaves us no other choice.
\hfill
$\Box$

\newpage
\part*{Chapter 4: Combinatorial Continuity}
\vspace{-0.5cm}
\indent
\section{Introduction}
In this Chapter we describe which equivalence
 relations in $\QS$ appear as the rational 
 lamination of polynomials with connected
 Julia set and all cycles repelling.
Conjecturally, every polynomial with
 all cycles repelling and connected Julia set
 lies in the set $\partial \shift \cap \conec$
 where the shift locus $\shift$ and the connectedness locus
 $\conec$ meet. 
Here, we also give a  description of 
 where in $\partial \shift \cap \conec$
 a polynomial with all cycles repelling 
 and a given rational lamination
 can be found.

Our descriptions will be in terms of critical portraits.
On one hand we will show that each critical portrait $\cp$
 gives rise to an equivalence relation $\ratcri(\cp)$ in $\QS$
 which is a natural candidate to be the rational 
 lamination of a polynomial.
On the other hand,
 in Chapter 2, we have already seen 
 that critical portraits determine directions 
 to go from the shift locus $\shift$ to the connectedness 
 locus $\conec$.
More precisely, each critical portrait $\cp$ determines
 a non-empty connected subset of $\partial \shift \cap \conec$
 called the impression $I_{\conec}(\cp)$ of $\cp$.
Thus, a location in $\partial \shift \cap \conec$ will be given 
 in terms of critical portrait impressions.

\medskip
In order to be more precise, 
 recall that a critical portrait $\cp$ partitions the circle $\S$
 into $d$ $\cp$-unlinked classes $L_1, \dots, L_d$
 (see Definition~\ref{unl-d}).
Symbolic dynamics of $m_d : t \rightarrow dt \mbox{ mod 1 }$
 with respect to this partition give us the right
 and left itineraries $\mathrm{itin}^{\pm}_{\cp}$
 (see Definition~\ref{itin-d}).
That is, we let:
 $$\begin{array}{lccl}
 \mathrm{itin}^{\pm}_{\Theta} : & \S & \rightarrow & 
        \{ 1, \dots , d \}^{{\N}\cup \{0\}} \\
                    & t  & \mapsto     &
                        ( j_0 , j_1 , \dots )\\
 \end{array}$$
 if, for each $n \geq 0$, there exists $\epsilon >0$ such that
 $(d^n t, d^n t\pm \epsilon) ) \subset L_{j_n}$.

\smallskip
Now, each critical portrait generates an equivalence 
relation in $\QS$:

\begin{definition}
Given a critical portrait
$\Theta$ we say that two arguments 
$t, t^\prime$ in $\QS$ are $\Lambda_{\Q}(\Theta)$-equivalent
if and only if there exists $t=t_1, \dots ,t_n=t^{\prime}$ such that
one of the two 
itineraries $\mathrm{itin}^{\pm}_{\Theta} (t_i)$ 
coincides with one of the two itineraries
$\mathrm{itin}^{\pm}_{\Theta} (t_{i+1})$ for $i =1 , \dots , n-1$.
\end{definition}

Notice that the above equivalence relation $\ratcri(\cp)$ is
closely related
to the landing pattern of rational external rays 
for a polynomial $g$ in the visible shift locus 
with critical portrait $\cp(g) = \cp$ 
(see Lemma~\ref{lan-vshift-l}).

\smallskip
Let us illustrate the above definition with some 
examples:

\noindent
{\bf Example:} 
Consider the cubic critical portrait 
$\cp = \{ \{ 1/3 , 2/3 \} , \{ 1/9 , 7/9 \} \}$.
The $\cp$-unlinked classes are 
$L_1 = (7/9, 1/9)$, $L_2 = (1/9,1/3) \cup (2/3,7/9)$
and, $L_3 = (1/3, 2/3)$. Observe that 
$\{ 1/9,2/9,7/9,8/9 \}$ is a $\ratcri(\cp)$-equivalence
class. In fact, 
$$\mathrm{itin}^{+}_{\cp} (7/9) = 
   \mathrm{itin}^{-}_{\cp}(8/9) = 13111...$$ 
$$\mathrm{itin}^{+}_{\cp} (2/9)= 
  \mathrm{itin}^{-}_{\cp}(7/9) = 22111...$$
$$\mathrm{itin}^{+}_{\cp}(8/9) = 
  \mathrm{itin}^{-}_{\cp}(1/9) = 12111...$$

\bigskip
\noindent
{\bf Example:}
Consider the cubic critical portrait 
$$\cp = \{ \{ 11/216 , 83/216 \}, \{ 89/216 , 161/216 \} \}.$$
The $\cp$-unlinked classes are $L_1 = (11/216,83/216)$,
$L_2 = (83/216,89/216) \cup (161/216,11/216)$ and
$L_3 = (89/216,161/216)$. 
It is not difficult
to see that   $$\{ 11/216, 17/216 , 83/216 ,
89/216 ,155/216 , 161/216 \}$$ is a 
$\ratcri(\cp)$-equivalence class
(compare with Example 3 in Section~\ref{dyn-s}).

\bigskip
A main distinction needs to be made according to whether
an argument which participates in $\cp$ has
a periodic itinerary or not. 

\smallskip
\begin{definition}
Consider a critical portrait $\cp = \{ \cp_1, \dots, \cp_m \}$.
We say that $\cp$ has {\bf periodic kneading} if for 
 some $\theta \in \cp_1 \cup \dots \cup \cp_m$ one of
 the itineraries $\mathrm{itin}^{\pm}_{\cp}(\theta)$ 
 is periodic under the one sided shift. Otherwise,
 we say that $\cp$ has {\bf aperiodic kneading}.
\end{definition}

The equivalence relations that arise from
 critical portraits with aperiodic kneading
 are exactly those that appear as the rational lamination
 of polynomials with all cycles repelling:

\begin{theorem}
\label{corres-th}
Consider an equivalence relation $\ratrel$ in $\QS$.
$\ratrel$ is the rational lamination 
 $\ratrel(f)$ of  some polynomial $f$
 with connected Julia set and all cycles repelling if and only if
 $\ratrel = \ratcri(\cp)$ for some critical portrait
 $\cp$ with aperiodic kneading.

Moreover, when the above holds, there are at most finitely many
critical portraits $\cp$ such that $\ratrel = \ratcri (\cp)$.
\end{theorem}

Given a polynomial $f$ with all cycles repelling the
 existence of a critical portrait $\cp$ such that
 $\ratrel(f) = \ratcri(\cp)$ is shown in 
 Section~\ref{from-rat-s}. 
In Section~\ref{from-rat-s}  we also give a necessary 
 and sufficient condition for a critical portrait
 $\cp$ to generate the rational lamination of $f$
 (Proposition~\ref{from-p}). 

Given a critical portrait $\cp$ with aperiodic kneading
 we find, in $\partial \shift \cap \conec$,
 a polynomial $f$ with  rational lamination $\ratcri(\cp)$.
More precisely, we show that the rational lamination
 of polynomials in the critical portrait impression 
 $I_{\conec}(\cp)$ 
 is exactly $\ratcri(\cp)$.
In particular, $\cp$ completely determines the rational
 lamination of the polynomials in $I_{\conec}(\cp)$:

\begin{theorem}
\label{cri-imp-th}
Consider a map $f$ in the impression 
 $I_{\conec}(\cp)$ of a critical portrait $\cp$.

If $\cp$ has aperiodic kneading then $ \ratrel(f) = \ratcri(\cp)$
and all the cycles of $f$ are repelling.

If $\cp$ has periodic kneading then  at least one
cycle of $f$ is non-repelling.
\end{theorem}

From the Theorems above,
 we conclude that a polynomial $f \in \partial \shift \cap \conec $
 with all cycles repelling
 must lie in at least one of the finitely many 
 impressions of critical
 portraits $\cp$ such that $\ratrel(f) = \ratcri(\cp)$.

\smallskip
A case of particular interest is when the critical
 portrait $\cp$ is 
 formed by strictly pre-periodic arguments.
Under this assumption, the impression $I_{\conec}(\cp)$
is the unique critically pre-repelling polynomial
$f$ such that, for each $\cp_i$, the external rays with arguments
in $\cp_i$ land at a common critical point of $f$.
In fact, this is a direct consequence of the above Theorem
and the ``combinatorial rigidity'' of critically
pre-repelling maps (see Corollary~\ref{extra-c}). 
By ``combinatorial rigidity''
we mean that the rational lamination  of a critically
pre-repelling map
uniquely determines the polynomial
(see~\cite{jones-91,bielefield-92}).

It is also worth mentioning that one obtains
a proof of the 
Bielefield-Fisher-Hubbard~\cite{bielefield-92}
 realization Theorem which bypasses the application of 
 Thurston's characterization of post-critically
 finite maps~\cite{douady-93a}.
That is,  given a critical
portrait $\cp = \{ \cp_1 , \dots , \cp_m \}$ formed by strictly
pre-periodic arguments there {\bf exists} 
a critically pre-repelling
map $f$ such that, for each $\cp_i$, 
the external rays with arguments
in $\cp_i$ land at a common critical point of $f$
(see Corollary~\ref{extra-c}).

\begin{figure}[htb]
\centerline{
\psfig{file=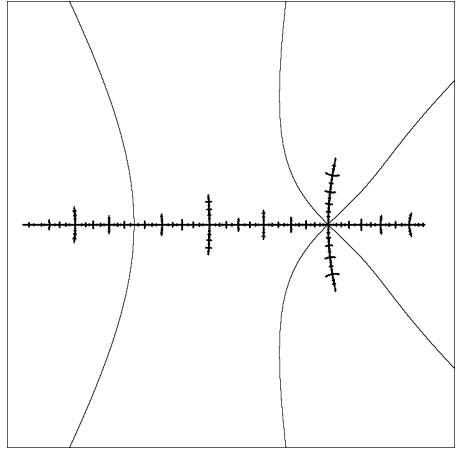,width=170pt}
}
\caption{The Julia set of the cubic polynomial 
$f(z) = z^3 - 9/4z + \sqrt{3}/4$
and the external rays landing at the critical points.}
\label{one-third}
\end{figure}

\smallskip
\noindent
{\bf Example:} The critically pre-repelling cubic polynomial
$f(z) = z^3 - 9/4z + \sqrt{3}/4$ has two critical points.
One critical point is the landing point of 
the external rays with arguments $1/3$ and $2/3$.
The other critical point is the landing point 
of the rays with arguments $1/9, 2/9, 7/9, 8/9$
(see Figure~\ref{one-third}).
Thus, Proposition~\ref{from-p} implies that 
the cubic critical portraits 
$$\cp = \{ \{ 1/3, 2/3 \} , \{ 1/9 , 7/9 \} \}
\mbox{ and }
\cp^\prime = \{ \{ 1/3, 2/3 \} , \{ 2/9 , 8/9 \} \}$$
are the only ones that 
generate the rational lamination of $f$. 
Moreover, putting together the fact that $f$
is uniquely determined by its rational lamination
with Theorem~\ref{cri-imp-th}, we have that
the  impressions $I_{\conec}(\cp)$ and $I_{\conec}(\cp^\prime)$
consists of the polynomial $f$.

\smallskip
\noindent
{\bf Example:} 
The polynomial $f(z) = z^3 + 0.2203 +1.1863I$ has a unique critical
point which is the landing point of the 
external rays with arguments
$$11/216, 17/216 , 83/216 ,
89/216 ,155/216 , 161/216.$$
Arguing as in the previous example, 
the cubic critical portraits
$$ \{ 11/216 , 83/216 , 155/216 \}$$
$$ \{ 17/216 , 89/216 , 161/216 \}$$
$$ \{ \{ 17/216 , 89/216 \} , \{ 11/216 , 155/216 \}\}$$
$$ \{ \{ 89/216 , 161/216 \} , \{ 11/216 , 83/216 \}\}$$
$$ \{ \{ 17/216 , 161/216 \} , \{ 83/216 , 155/216 \}\}$$
have as impression the polynomial 
$f(z) = z^3 + 0.2203 +1.1863I$.
Notice that although $f$ has a unique multiple critical point
there are critical portraits formed by two sets that 
generate the rational lamination of $f$.

\begin{figure}[htb]
\centerline{
\psfig{file=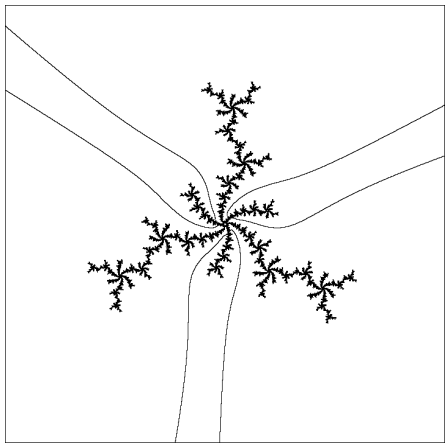,width=170pt}
}
\caption{The Julia set of the cubic polynomial 
$f(z) = z^3 + 0.2203 +1.1863I$
and the external rays landing at the critical point.}
\end{figure}

\smallskip
Note that the Mandelbrot local 
connectivity conjecture says that 
impressions of quadratic critical portraits are a single map.
We would like to stress that, for higher degrees,
 we do not expect this to be true.
That is, there might be non-trivial impressions of critical 
portraits with aperiodic kneading.

\section{From $\ratrel(f)$ to $\ratcri(\cp)$}
\label{from-rat-s}
\noindent \indent
In this section we show that every rational lamination 
$\ratrel(f)$
of a polynomial $f$ with all cycles repelling can be realized
as $\ratcri(\cp)$, for some critical portrait $\cp$.
Recall that,
under the assumption that all cycles of
$f$ are repelling, the rational lamination
 $\ratrel(f)$ is a maximal 
equivalence relation 
with finite and unlinked classes (Proposition~\ref{lam-p}).
 Thus, the strategy is first to
show that, for an arbitrary critical
portrait $\cp$, the equivalence relation $\ratcri(\cp)$
has finite and unlinked
classes (Lemma~\ref{ratcri-l} below).
Then we proceed to prove that 
$\ratrel(f) \subset \ratcri(\cp)$ for some critical portrait $\cp$
(Proposition~\ref{from-p}) and, by maximality of $\ratrel(f)$
we  conclude that $\ratrel(f) = \ratcri(\cp)$. 
 
\begin{lemma}
\label{ratcri-l}
Let $\Theta$ be a critical portrait.

Every $\ratcri(\cp)$-equivalence class 
$A$ is a finite set.

If $A_1$ and $A_2$ are distinct 
$\ratcri(\cp)$-equivalence classes 
then $A_1$ and $A_2$ are unlinked.
\end{lemma}

These properties of $\ratcri(\cp)$ can be 
proven abstractly.
Nevertheless, the intuition behind it is that there
exists a polynomial $g$ in the visible shift locus whose 
rational external
rays land in a pattern very closely related to that  
given by $\ratcri(\cp)$ (Lemma~\ref{lan-vshift-l}).
 Our proof will make use of this
fact:

\medskip
\noindent
{\sc Proof:}
Choose a polynomial $g$ in the visible shift 
locus with critical portrait $\Theta(g) = \Theta$.

To show that $A$ is finite we pick $t_1 \in A$
and make a distinction according to whether $t_1$
is periodic or pre-periodic. 
In the case that $t_1$ is periodic of period $p$,
it is enough to verify that all the elements of $A$ are 
periodic of period $p$. 
In fact, take $t_2 \in A$ such that
 $\mathrm{itin}^{\epsilon_1}_{\cp}(t_1) = 
 \mathrm{itin}^{\epsilon_2}_{\cp}(t_2)$ where 
 $\epsilon_1, \epsilon_2 \in \{ +, - \}$.
In view of Lemma~\ref{lan-vshift-l},  
 the external rays $R^{t^\epsilon_1}_g$ 
 and $R^{t^{\epsilon_2}_2}_g$ land at a periodic
 point $z$. 
By Theorem~\ref{dis-lan-th}, both $t_1$ and $t_2$
 have period $p$. 
It follows that all the elements
 of $A$ have the same period $p$. 
Now, in the case that $t_1$ is pre-periodic, say
 that $d^l t_1$ is periodic of period $p$.
For each element $t \in A$, a similar argument
 shows that $d^l t$ is periodic of period $p$. 
Thus, $A$ is a finite set of pre-periodic arguments. 

\medskip
Now we must show that the $\ratcri(\cp)$-equivalence classes 
 $A_1$ and $A_2$ are unlinked.
That is, $A_2$ is contained in a connected 
 component of $\S \setminus A_1$.
In fact, consider the union $\Gamma_1$ (resp. $\Gamma_2$)
 of all the external rays 
 $R^{t^\pm}_g$ with arguments $t \in A_1$ (resp. $t \in A_2$)
 and their landing
 points. 
Observe that $\Gamma_1$ and $\Gamma_2$ are disjoint connected
 sets.
Thus, $\Gamma_2$ is contained in a connected component of 
 $\C \setminus \Gamma_1$.
In the circle at infinity, it follows that $A_2$ is contained 
in a connected component of $\S \setminus A_1$.
\hfill
$\Box$

\bigskip
Now we characterize the critical portraits that generate
a given rational lamination:

\newpage
\begin{proposition}
\label{from-p}
Consider a  polynomial $f$ with connected Julia set
 $J(f)$ and all cycles repelling. 
Let $\ratrel(f)$ be its rational lamination. 
Then there exists at least
 one and at most finitely many critical portraits
 $\Theta $ such that
 $$\Lambda_{\Q}(\Theta) = \ratrel (f).$$
Moreover, all such critical portraits have aperiodic
 kneading. 
Furthermore, $\Theta = \{ \Theta_1 , \dots , \Theta_m \}$ 
 satisfies the identity above if and only if each
 $\cp_i$ is either contained in a $\ratrel(f)$-equivalence
 class or $\cp_i$ is unlinked
 with all $\ratrel(f)$-equivalence classes.
\end{proposition}

\noindent
{\sc Proof of Proposition:}
Corollary~\ref{exist-c} provides us with a critical portrait
$\cp = \{ \cp_1, \dots, \cp_m \}$ such that 

(a) If $\cp_i \subset \QS$ then $\cp_i$ is
 contained in a $\ratrel(f)$-equivalence class,

(b) If $\cp_i \subset \S \setminus \QS$ 
then $\cp_i$ is unlinked with all $\ratrel(f)$-equivalence classes. 

We show that these are sufficient conditions to
conclude that $\ratrel(f) = \ratcri(\cp)$.

Consider an arbitrary  $\ratrel(f)$-equivalence class
$A$ and let $t_1, t_2$ be two consecutive elements of $A$.
That is, $(t_1, t_2)$ is a connected component of
$\S \setminus A$. 
Observe  that (a) and (b) guarantee that the first symbol of 
$\mathrm{itin}^{+}_{\cp}(t_1)$ coincides
with the first symbol of  $\mathrm{itin}^{-}_{\cp}(t_2)$.
In view of property (R4) of $\ratrel(f)$ in Proposition~\ref{lam-p},
  $d^n t_1$ and $d^n t_2$ are consecutive
elements of $m^{\circ n}_d (A)$. 
Thus, $\mathrm{itin}^{+}_{\cp}(t_1) = \mathrm{itin}^{-}_{\cp}(t_2)$
and $A$ is contained  in a $\ratcri(\cp)$-equivalence class.
That is, $\ratrel(f) \subset \ratcri(\cp)$.
Now Proposition~\ref{lam-p} says that $\ratrel(f)$ 
is a maximal equivalence relation with finite and 
unlinked classes. By Lemma~\ref{ratcri-l}, 
we conclude $\ratrel(f) = \ratcri(\cp)$

\smallskip
Assume that $\ratrel(f) = \ratcri(\cp)$, 
we must prove that the conditions (a) and (b) hold. 
By contradiction, suppose that  $\cp_i$ lies 
in $\QS$ and it is not contained in a $\ratrel(f)$-class.
By Lemma~\ref{approx-l}, $\cp_i$ is linked with infinitely many
$\ratrel(f)$-classes. A class 
of $\ratcri(\cp)$ which is linked with $\cp_i$ must contain
an element of $\cp_i$. Therefore, $\ratcri(\cp)$ 
would have an infinite class. This contradicts Lemma~\ref{ratcri-l}
and implies that (a) holds.
If $\cp_i$ lies in $\S \setminus \QS$,  Lemma~\ref{e-l} allows us to
apply a similar reasoning to show that $\cp_i$ is unlinked
with every $\ratrel(f)$-class.

\smallskip
To show that a critical portrait $\cp$ such that 
 $\ratrel(f) = \ratcri(\cp)$
 must have aperiodic kneading, consider a polynomial
 $g$ in the visible shift locus such that 
 $\cp = \cp(g) = \{ \cp_1 , \dots , \cp_m \}$.  
No periodic $t \in \QS$
 participates in $\cp$, otherwise some $\cp_i$  contains 
 periodic and pre-periodic arguments and cannot 
 be contained in a $\ratrel(f)$-class as proved above.
Thus, all the periodic rays of $g$ are smooth.
Since $\ratrel(f) = \ratcri(\cp)$, 
 if $A \subset \QS$ is 
 a periodic point type of $f$ then $A$ is  
 a periodic point type of $g$.
Thus, $g$ has exactly $d^p$ rational periodic
 point types of period dividing $p$. 
This matches with
 the $d^p$ periodic points of $g$ of  period dividing $p$.
Therefore, no ray that bounces off a critical point
 of $g$ can land at a periodic point of $g$. 
It follows that $\cp(g) = \cp$ has aperiodic kneading.
\hfill $\Box$

\begin{lemma}
\label{extra-l}
Consider a critical portrait $\cp = \{ \cp_1 , \dots , \cp_m \}$
such that all the arguments in $\cp_1 \cup \dots \cup \cp_m$
are strictly pre-periodic. 
If $f$ is a polynomial such that $\ratrel(f) = \ratcri(\cp)$
then all the critical points of $f$ are strictly pre-periodic.
\end{lemma}

\noindent
{\sc Proof: } It is not difficult to check that $\cp$
has aperiodic kneading. 
A counting argument as above
shows that $f$ has all cycles repelling. 
Proposition~\ref{from-p} implies that the
external rays with arguments in $\cp_i$ land at a common
critical point. 
Counting multiplicities, it follows that all the critical
points of $f$ are the landing point of some $\cp_i$.
Thus, $f$ is a critically pre-repelling map.
\hfill
$\Box$

\section{Critical Portraits  with \\ aperiodic kneading}

\begin{figure}[htb]
\label{sep-f}
\centerline{
\psfig{file=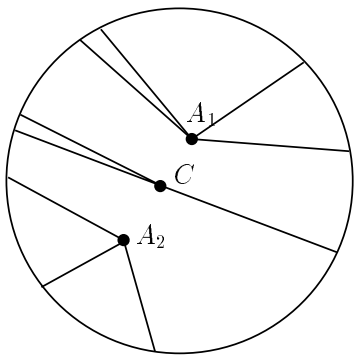,width=136pt}
}
\caption{The diagram illustrates the 
situation of Lemma~\ref{sep-l}.}
\end{figure}

\noindent \indent
In the next section we are going to show that 
 polynomials $f$ in the impression $I_{\conec}(\cp)$ of a critical 
 portrait $\cp$ with aperiodic kneading has all cycles
 repelling. 
A key step in the proof is to rule out 
 parabolic cycles of $f$. 
To do so we will show that 
 two periodic points cannot coalesce as one goes 
 from the shift locus to the connectedness locus 
 in the ``direction'' determined by $\cp$. 
Roughly, the obstruction
 for this collision to occur is that, for polynomials
 in $\vshift$ with critical portrait $\cp$, any
 two periodic points are separated by the external rays
 landing at a pre-periodic
 point (see Figure~\ref{sep-f}): 

\begin{lemma}[Separation]
\label{sep-l}
Consider a critical portrait $\Theta$ with aperiodic
 kneading and let $A_1$ and $A_2$ be two distinct 
 periodic $\ratcri(\Theta)$-equivalence classes.
Then there exists a strictly pre-periodic 
 $\ratcri(\cp)$-equivalence class $C$
 such that $A_1$ and $A_2$ lie in different
 connected components of $\S \setminus C$.
Moreover, $C$ can be chosen such that, for all 
 $n \geq 0$, $m^{\circ n}_d(C)$ is contained in
 a $\cp$-unlinked class.
\end{lemma}

Note that for $\cp^\prime$ close to $\cp$ each
 of the sets $A_1$, $A_2$ and $C$ is also contained
 in a $\ratcri(\cp^\prime)$-equivalence class.

\smallskip
The idea of the proof is similar to that of 
Lemma~\ref{approx-l}:

\noindent
{\sc Proof:}
Pick a polynomial $f$ in the visible shift locus
 with critical portrait $\cp = \cp(f)$.
Observe that since $\cp$ has aperiodic kneading
 all the external rays landing at periodic 
 points are rational and smooth. 
Moreover, the periodic point types of $f$ are exactly
 the periodic equivalence classes of $\ratcri(\cp)$.
Let $z_1, z_2$ be such that $A_1 = A(z_1)$
 and $A_2 = A(z_2)$. 
The idea is to construct, with smooth rays, a puzzle piece
 around $z_1$ which does not contain $z_2$.
For this, consider the collection $\mathcal{G}$ 
 formed by the rational typess $A(w)$ such that:

(i) All the external rays landing at a point in the
grand orbit of $w$ are smooth.

(ii) $w$ is not in the grand orbit of $z_1$.

We pass to an iterate of $f$ so that every periodic
 point $w$ whose type does not lie in $\mathcal{G}$
 is a fixed point and all the rays landing at $w$ are 
 fixed rays. 
In particular, now $z_1$ is a fixed
 point which is the landing point of fixed rays.

We saturate $\mathcal{G}$ by an increasing sequence of 
 finite collections $G_l$.
That is, let $G_l$ be the collection of types $A(w)$
 in $\mathcal{G}$ such that $f^{\circ l}(w)$ is periodic of 
 period less or equal to $l$.
Let $U_l(z_1)$ be the $G_l$-puzzle piece which contains
 $z_1$. 

It is enough to show that, for some $l$, $f$ maps $U_l(z_1)$ 
 homeomorphically onto its image. 
We suppose that this is not the case and, after some
 work, we arrive to a contradiction.
That is, suppose that for $l \geq l_0$ the map 
 $f_{|U_l}(z_1)$ has degree $k >1$. 
We can pick $l_0$ large enough so that every fixed 
 point in $\overline{U}_{l_0}(z_1)$ is the landing
 point of fixed external rays.

\noindent
{\it Claim 1:} For $l \geq l_0$, 
$f$ has $k$ fixed points in $\overline{U}_l(z_1)$.

\smallskip
\noindent
{\it Proof of Claim 1:} 
Let $\rho >0$ be large enough such that the equipotential
$g_f = \rho$ is a topological circle (i.e. all the critical
points of $f$ are contained inside $g_f = \rho$).
Consider the portion $U$ of $U_l(z_1)$ contained inside
this equipotential. Let $U^{\prime} = f(U)$.
Observe that $U$ and $U^{\prime}$ are puzzle pieces
that satisfy the conditions of the thickening
Lemma~\ref{thick-l}. The same construction shows that, 
after extracting a polynomial  like map of degree $k$,
we must have $k$ fixed points of $f$ in 
$\overline{U}_l(z)$.

\smallskip
Now let $V_n$ be the connected component of $f^{-n}(U_{l_0}(z_1))$
 that contains $z_1$ and 
 $$X = \cap_{n \geq 0} f^{-n}(\overline{V_n}).$$
In the circle at infinity, 
 let $\pi_{\infty} V_n = \{ t \in \S : R^{*t}_f \subset V_n \}$
 and
 $$E = \cap_{n \geq 0} \overline{\pi_{\infty} {V_n}}.$$

Observe that $m_d$ is $k$ to $1$ on $E$.
Moreover, $E$ contains elements of $k$ distinct
 fixed point types formed by fixed rays.
Our aim is to show that $E$ can only intersect $k-1$
 fixed point types which are formed by fixed arguments.
Roughly speaking, we will construct a semiconjugacy between
 $m_d|_E$ and a degree $k$ selfcovering of a circle.
For this, we need some basic facts about $E$.

\noindent
\smallskip
{\it Claim 2:}  If $(t_1, t_2)$ is a connected component of
 $\S \setminus E$ then $(dt_1, dt_2)$ is also a connected
 component of $\S \setminus E$. 
Moreover, the external rays
 $R^{t^+_1}_f$ and $R^{t^-_2}_f$ land at the same point.

\noindent
{\it Proof of Claim 2:}
Notice $(t_1,t_2) = \cup_{n \geq 1} (a_n,b_n)$ 
where $(a_n, b_n)$ is a connected component of
$\S \setminus \overline{\pi_{\infty}V}_n$.
Since $a_n,b_n$ are rational and 
 $f$ acts preserving the cyclic order of the rational 
rays that bound $V_n$ we have that
$(da_n, db_n)$ is a connected component of 
$\S \setminus \overline{\pi_{\infty}V}_{n-1}$.
Thus, $(t_1,t_2 )$ is a connected component of
$\S \setminus E$. 
Since the $\cp$-itineraries of $a_n$ and $b_n$
agree, $R^{t^+_1}_f$ and $R^{t^-_2}_f$ land at the same point.

\smallskip
Also, we need some control over the (possibly) isolated points
of $E$:

\smallskip
\noindent
{\it Claim 3:} Let $\{ t_i \} \subset E$ such that 
$(t_i, t_{i+1})$ is a connected component of $\S \setminus E$.
Then $\{ t_i \}$ is finite.

\noindent
{\it Proof of Claim 3:} 
There are two possibilities. In the case that 
there is a rational $t_{i_1} \in \{ t_i \}$, then
$R^{t^+_{i_1}}$ lands at a pre-periodic or periodic
point $z$. The previous Claim implies that the external
ray $R^{t^-_{i_1 +1}}_f$ also lands at $z$.
Hence, $t_{i_1 +1}$ is rational and it is
 $\ratcri(\cp)$-equivalent to $t_{i_1}$. It follows
that all the elements of $\{ t_i \}$ are 
$\ratcri(\cp)$-equivalent. Since  
$\ratcri(\cp)$-classes are finite, we conclude
that $\{ t_i \}$ is finite.

If all the elements of $\{ t_i \}$ are irrational
and this set is infinite then there exists two elements
$t_{j}, t_{j+l}$ such that
$d^{k_1} t_{j} = d^{k_2} t_{j + l} = \theta$   where 
$\theta \in \cp_1 \cup \cdots \cup \cp_m$ and $k_1 \neq k_2$.
For $n$ large enough, all the rays with arguments
$d^n t_{j}, d^n t_{j +1},  \dots, d^n t_{j +l}$ are
smooth. By the previous Claim, all these rays must land at the same
point $z$. It follows that the external rays
with arguments $d^{n -k_1} \theta$ and $d^{n -k_2} \theta$ land
at the same point $z$. Since $k_1 \neq k_2$ we have that
$z$ must be periodic. Thus, $\theta \in \QS$
and  $t_{j}$ must also be rational. Which puts us in the first 
case.

\medskip
It follows that the set $\tilde{E}$ obtained by removing from
$E$ its isolated points is a Cantor set.
Moreover, if $(t_1, t_2)$ is a connected component
of $\S \setminus \tilde{E}$ then $(dt_1, dt_2)$ is
also a connected component of $\S \setminus \tilde{E}$.
Every fixed point type that had an element in $E$
also has an element in $\tilde{E}$.
Furthermore, ${m_d}_{|\tilde{E}}$ is $k$ to $1$.

Now consider the quotient $\T$ obtained from $\S$ by
identifying $[t_1,t_2]$ to a point if and only if 
$(t_1,t_2)$ is a connected component of $\S \setminus \tilde{E}$.
Let $$h : \S \rightarrow {\T}$$
 be the quotient map.
It follows that $m_d$ projects to a degree $k$ selfcovering $g$
of the topological circle $\T$. 
Moreover, each fixed
point of $g$ is either the image of a fixed point $t$ of $m_d$ or 
of an interval $[t_1, t_2]$ whose endpoints are
fixed points of $m_d$. 
Observe that, in the latter case, 
$R^{t_1}_f$ and $R^{t_2}_f$ land at the same fixed point of $f$.
Therefore, $h$ does not identify
arguments of distinct fixed point types 
that have elements in $\tilde{E}$. 
Recall that there are
$k$ fixed point types with rotation number zero 
that have elements in $\tilde{E}$. 
Thus, $g$ has at least $k$ fixed points.
But, every fixed
point of $g$ is topologically repelling. Thus, $g$ has $k-1$
fixed points. Contradiction.
\hfill $\Box$

\smallskip
The previous Lemma says that periodic points are
separated by smooth pre-periodic rays. 
The next Lemma will allow us to show that this 
separation persists in the limit when we go from the
shift locus $\shift$ to the connectedness locus $\conec$. 

In order to make this precise, we consider
a sequence of external rays $R^t_{f_n}$ and
say that 
$$\limsup R^t_{f_n}$$
is the set of points $z \in \C$ such that every neighbourhood
of $z$ intersects infinitely many $R^t_{f_n}$.
This coincides with the usual definition of the 
Hausdorff metric on compact subsets of the Riemann sphere.

\begin{lemma}
\label{limsup-l}
Let $\cp$ be a critical portrait with aperiodic 
kneading. Consider  a sequence $f_n \in \vshift$ such that
$\cp(f_n) \rightarrow \cp$ and
$f_n \rightarrow f \in \conec$.

If $t \in \QS$ is periodic and
$z \in \limsup R^t_{f_n} \cap J(f)$ then 
$z$ is periodic under $f$.

If $t \in \QS$ is such that $t$ is
not periodic and $z \in \limsup R^{t^\pm}_{f_n} \cap J(f)$ then 
$z$ is not periodic under $f$.
\end{lemma}

\noindent
{\sc Proof:} 
Assume $t$ is periodic of period $p$. Since $\cp$
has aperiodic kneading, for $n$ sufficiently large,
$R^t_{f_n}$ is smooth.

\noindent
{\it Claim 1:}
Consider $z_n \in R^t_{f_n}$,
we claim that, for $n$ large enough,
$$\rho_{\Omega(f_n)} (z_n , f^{\circ p}_n ( z_n)) < 2 p \log d$$
where $\rho_{\Omega(f_n)}$ is the hyperbolic 
metric in $\Omega(f_n)$.

\smallskip 
Let us postpone this estimate and proceed with the proof.
If $$z \in \limsup R^t_{f_n} \cap J(f)$$ then consider
a sequence $ z_n \in R^t_{f_n} $ which converges to $z$.
Since repelling cycles of $f$ are dense in $J(f)$,
the euclidean distance between $z_n$ and $J(f_n) = \partial \Omega(f_n)$
goes to zero. In the other hand, the hyperbolic distance 
between $z_n$ and $f^{\circ p}_n (z_n)$ stays bounded.
The standard comparison between the hyperbolic metric
and the euclidean metric yields that $z_n$ and $f^{\circ p}_n(z_n)$
must converge to the same point $z$. It follows that
$f^{\circ p}(z) = z$, i.e. $z$ is periodic.

Now, in the case that $t \in \QS$ is strictly pre-periodic 
and we first consider the case in which $dt$ is a periodic 
argument. If $w \in \limsup R^{t^\pm}_{f_n} \cap J(f)$
then we proved that $f(w)$ must be periodic. 
There are two
possibilities, either $w$ is the unique periodic
preimage of $f(w)$ or $w$ is strictly pre-periodic.
We claim that only the latter occurs. 
By contradiction, suppose that $w$ is periodic.
Let $w_n \in R^{t^\pm}_{f_n}$ be such that 
$w_n \rightarrow w$ and 
$t^{\prime} \in \QS$ be the unique periodic 
preimage of $t$. Consider $w^{\prime}_n \in R^{t^\prime}_f$
such that $f_n (w^{\prime}_n) = f_n (w_n)$. Since $t^{\prime}$
is periodic,  
$w^{\prime}_n$ must converge to the unique periodic 
pre-image $w$ of the
periodic point $f(w)$. Thus, $f$ is not locally injective
at $w$, but $w$ is a periodic point in the Julia set $J(f)$.
Contradiction. 
The general case, for an arbitrary strictly pre-periodic $t$,
follows.

\smallskip
\noindent
{\it Proof of Claim 1:}
The estimate will follow from Lemma~\ref{ls-l}. Let 
$H = \{ z = x+ iy : x >0 \}$ and consider
the region $\tilde{V} = 2 \pi i t + \{ z \in H : | \arg z | < \pi/4 \}$.
By Lemma~\ref{ls-l}, for $n$ sufficiently
large,  $exp ( \tilde{V} ) \subset U_{f_n}$, where $U_{f_n}$
is the image of the B\"ottcher map $\phi_{f_n} : \Omega^*(f_n) 
\rightarrow U_{f_n}$. Observe that:
$$z_n = \phi^{-1}_{f_n} \circ exp (g_{f_n} (z_n) + 2 \pi i t )$$
$$f^{\circ p}_n (z_n) = 
\phi^{-1}_{f_n} \circ exp ( d^p g_{f_n} (z_n) + 2 \pi i t ).$$
Also, 
$$\rho_{\tilde{V}}( g_{f_n} (z_n) + 2 \pi i t , 
d^p g_{f_n} (z_n) + 2 \pi i t ) = 2 p \log d$$
where $\rho_{\tilde{V}}$ is the hyperbolic metric in $\tilde{V}$.
Since, $\phi^{-1}_{f_n} \circ exp$ is a contraction the Claim
follows.
\hfill
$\Box$

\section{Combinatorial Continuity}

\indent \noindent
Before we prove Theorem~\ref{cri-imp-th} we need
the following Lemma due to Goldberg and Milnor~\cite{goldberg-93}:

\begin{lemma}
\label{gm-l}
Consider $f_0 \in \param$ such that $J(f_0)$ is
connected.
Assume that $R^t_{f_0}$ is a smooth external 
ray which lands at a pre-repelling
or repelling periodic point $z_{f_0}$.
Also, assume that $z_{f_0}$ is not a pre-critical 
point.
Then, for any $f$ sufficiently close to  
$f_0$, the external ray $R^t_{f}$ is smooth and 
lands at the analytic continuation $z_f$ of $z_{f_0}$.
\end{lemma}

\noindent
{\sc Proof of Theorem~\ref{cri-imp-th}:}
Consider a critical portrait $\cp$ with aperiodic
kneading and a polynomial $f$ in its impression $I_{\conec}(\cp)$.
So let $f_n$ be a sequence in the visible shift locus
such that:
$$f_n \rightarrow f$$
$$\cp(f_n) \rightarrow \cp.$$

The next two claims  show that $f$ must have all cycles 
repelling:

\smallskip
\noindent
{\it Claim 1:} $f$ does not have a parabolic cycle.

\noindent
{\it Proof of Claim 1:}
By contradiction, suppose that $z$ is a parabolic 
 periodic point of $f$. 
It follows that there
 exists distinct repelling periodic points $z_1(n)$ and $z_2(n)$ 
 of $f_n$ that converge to $z$. 
Also, 
 the periods of $z_1(n)$ and $z_2(n)$ divide
 a fixed number $p$. 
Since $\cp$ has aperiodic kneading the itinerary of
 periodic elements of $\QS$ vary continuously. 
In particular,
 for $\cp^{\prime}$ in a sufficiently
 small neighbourhood of $\cp$
 the periodic classes of $\ratcri(\cp)$ with period
 dividing $p$ coincide with those of $\ratcri(\cp^{\prime})$.
After passing to a subsequence, we may assume that
 the type of $z_1(n)$ is $A_1$ and the type
 of $z_2(n)$ is $A_2$, where $A_1$ and $A_2$ are distinct
 equivalence classes of $\ratcri(\cp)$. 
By Lemma~\ref{sep-l} we know that there exists a strictly pre-periodic
 class $C$ of $\ratcri(\cp)$ such that it separates $A_1$
 and $A_2$ (i.e. $A_1$ and $A_2$ lie in different connected
 components of $\S \setminus C$).
Recall that none of the elements of $C$ and its forward
 orbit participate of $\cp$.
Hence,
 for $\cp^{\prime}$
 in a sufficiently small neighborhood of $\cp$,
 the elements of $C$ are also identified by $\ratcri(\cp^{\prime})$.
Thus, for $n$ large the external rays of $f_n$ with arguments
 in $C$ together with their landing points form a connected
 set $\Gamma_n$. 
Now the periodic
 points $z_1(n)$ and $z_2(n)$ lie in different
 connected components of $\C \setminus \Gamma_n$.
Passing to the limit, $z \in \limsup \Gamma_n$
 and therefore $z \in \limsup R^t_{f_n}$, for some
 $t \in C$. 
This contradicts Lemma~\ref{limsup-l} and shows
 that $f$ cannot have parabolic cycles.

\medskip
\noindent
{\it Claim 2:} $f$ does not have irrationally neutral cycles.

\noindent
{\it Proof of Claim 2:} Again by contradiction we suppose
that $z$ is an irrationally neutral periodic point of 
$f$ with period $p$. There exists a sequence $z(n)$ of
periodic points of $f_n$ that converge to $z$. As above,
after passing to a subsequence, we may assume that
$z(n)$ has type $A \subset \QS$, for all $n$. Pick an element
of $t \in A$ and notice that $R^t_f$ must land at
a repelling periodic point $w$. Hence, for $n$ large enough,
$R^t_{f_n}$ has to land in the analytic continuation
of $w$ which is not $z(n)$.

\smallskip

\noindent
{\it Claim 3:} If $f$ is a polynomial
with all cycles repelling  and $f \in I_{\conec}(\cp)$
for some critical portrait $\cp \in \ang$ then
$\ratcri(\cp) = \ratrel(f)$.

\smallskip
It follows from Proposition~\ref{from-p} that $\cp$ has aperiodic kneading.

\noindent
{\it Proof of Claim 3:} Consider
$f_n$ in the visible shift locus
such that $$f_n \rightarrow f$$
and $$\cp(f_n) \rightarrow \cp = \{\cp_1 , \dots , \cp_m \}$$
Now we show that each $\cp_i$ is either contained
in a $\ratrel(f)$-class or it is unlinked with all $\ratrel(f)$-classes.
If $\cp_i$ has elements in two distinct
equivalence classes of $\ratrel(f)$
then $\cp_i$ is linked with the rational type $A(w)$
of a point $w$ which is not in the grand orbit of a critical
point. For $\cp^{\prime}$ close to $\cp$,
we also have that $C = A(w)$ has points in two different
$\cp^{\prime}$-unlinked classes.
 In view of Lemma~\ref{gm-l}, for $n$ large enough, the external
rays of $f_n$ with arguments in $C$ are smooth and land at a 
common point. Therefore, they intersect the external radii
of $f_n$ terminating at some critical point.
This is impossible.
A similar situation occurs when $\cp_i \subset \S \setminus \QS$
 is linked with a $\ratrel(f)$-class.
By Proposition~\ref{from-p}, $\ratcri(\cp) = \ratrel(f)$.
\hfill $\Box$

\medskip
Theorem~\ref{corres-th} follows from Theorem~\ref{cri-imp-th}
and Proposition~\ref{from-p}.

\begin{corollary}
\label{extra-c}
Assume $\cp= \{ \cp_1, \dots, \cp_m \}$ 
is a critical portrait formed by strictly
pre-periodic arguments. Then the critical portrait impression 
$I_{\conec}(\cp)$ is formed by the 
unique critically pre-repelling 
polynomial $f$ such that, for each $\cp_i$,
the external rays with arguments in $\cp_i$ land at
a common critical point.
\end{corollary}

\noindent
{\sc Proof:} 
By Theorem~\ref{cri-imp-th} and Lemma~\ref{extra-l}, 
each polynomial $f$ in the impression $I_{\conec}(\cp)$
is critically  pre-repelling and 
such that $\ratrel(f) = \ratcri(\cp)$.
According to Proposition~\ref{from-p},
 this occurs if and only if the external rays with 
 arguments in $\cp_i$ land
 at same point.
As mentioned above it follows from the 
work of Jones~\cite{jones-91} or the work of Bielefield,
Fisher and Hubbard~\cite{bielefield-92} 
that $f$ is uniquely determined
by its rational lamination $\ratcri(\cp)$.
\hfill
$\Box$

\bigskip\bigskip

This paper is a copy of the author's Ph. D. Thesis
produced at the Mathematics Deparment, S.U.N.Y. at
Stony Brook, U.S.A.

\smallskip
\noindent
Current Address:

{\sc 
\noindent
Facultad de Matem\'aticas,

\noindent
Pontitificia Universidad Cat\'olica de Chile.

\noindent
Casilla 306, Correo 22, Santiago,

\noindent
Chile. }
  
\noindent
{\tt jkiwi@riemann.mat.puc.cl}

\vfill\eject
\end{document}